\input amssym
\font\goth=eufm10
\def\H{{\cal H}}
\def\K{{\cal K}}
\def\P{{\cal P}}
\def\g{\hbox{\goth g}}
\def\G{\widehat{\hbox{\goth g}}}
\def\s{\hbox{{\goth sl}$_2$}}
\def\u{\hbox{\goth u}}
\def\a{\hbox{\goth a}} 
\def\v{\hbox{\goth vir}}
\def\w{\hbox{\goth witt}}
\bf \centerline{DIRECT PROOFS OF THE FEIGIN--FUCHS CHARACTER FORMULA}
\centerline{FOR UNITARY REPRESENTATIONS OF THE VIRASORO ALGEBRA}
\vskip .1in
\centerline{Antony Wassermann, CNRS, Institut des Math\'ematiques de Luminy}
\vskip .2in
\noindent \bf ABSTRACT. \it Previously we gave a proof of the Feigin--Fuchs character formula for the irreducible
unitary discrete series of the Virasoro algebra with $0<c<1$. The proof showed directly that the mutliplicity 
space arising in the coset construction of Goddard, Kent and Olive was irreducible, using the elementary part of the
unitarity criterion of Friedan, Qiu and Shenker, giving restrictions on $h$ for $c=1-6/m(m+1)$ with $m\ge 3$. In this
paper we consider the same problem in the limiting case of the coset construction for $c=1$. Using primary fields, we
directly establish that the Virasoro algebra acts irreducibly on the multiplicity spaces of irreducible 
representations of $SU(2)$ in the two level one irreducible representations of the corresponding affine 
Kac--Moody algebra. This gives a direct proof that the only singular vectors in these representations are those given by
Goldstone's formulas, which also play an important part in the proof. For this proof, the theory is 
developed from scratch in a self--contained semi--expository way. Using the Jantzen filtration and the Kac determinant formula, 
we give an additional independent proof for the case $c=1$ which generalises 
to the case $0<c<1$, where it provides an alternative approach to that of Astashkevich. 
\vskip .2in
\bf\centerline{CONTENTS}
\vskip .05in
\rm
\noindent 1. Introduction.

\noindent 2. Constructions of positive energy representations.

\noindent 3. Character formulas for affine Lie algebras.

\noindent 4. Existence of singular vectors in the oscillator representations.

\noindent 5. Uniqueness of singular vectors in the oscillator representations.

\noindent 6. Density modules, primary fields and the Feigin--Fuchs product formula.

\noindent 7. Multiplicity one theorem.

\noindent 8. Proof of the Feigin--Fuchs product formula using primary fields.

\noindent 9. Proof of the character formula for $L(1,j^2)$ using the Jantzen filtration.

\noindent 10. Proof of the character formula for the discrete series $0<c<1$ using the Jantzen filtration.

\noindent Appendix A: Alternative proofs of the Fubini--Veneziano relations.

\noindent Appendix B: Explicit construction of singular vectors and asymptotic formulas of Feigin--Fuchs.

\noindent Appendix C: Proof of Feigin--Fuchs product formula using explicit formula for singular vectors.

\noindent Appendix D: Holomorphic vector bundles and flat connections.

\noindent References. 
\vskip .2in
\noindent \bf 1. Introduction. \rm The Witt algebra $\w$ is the Lie algebra of complex vector fields 
$f(\theta) \cdot d/d\theta$ on the unit circle with $f(\theta)$ a trigonometric polynomial. It has a basis 
$$\ell_n=ie^{in\theta}{d\over d\theta}$$
for $n\in {\Bbb Z}$ and commutation relations
$$[\ell_m,\ell_n]=(m-n)\ell_{m+n}.$$
This Lie algebra has a central extension by ${\Bbb C}$, called the Virasoro algebra $\v$,
with basis $L_n$ ($n\in {\Bbb Z}$), $C$ and the commutation relations
$$[L_m,L_n]=(m-n)L_{m+n} + \delta_{m+n,0} {m^3-m\over 12} C.$$
In this article we will be interested in positive energy unitary representations of the Virasoro algebra.
By this we mean an inner product space $\H$ on which the Lie algebra $\v$ acts by operators $\pi(L_n)$ satisfying
$\pi(L_n)^*=\pi(L_{-n})$, $\pi(C)=\pi(C)^*$, with $\pi(C)$ acting as a scalar $cI$ and
such that $\H$ has a decomposition as an algebraic orthogonal direct sum of egenspaces of $\pi(L_0)$
$$\H=\oplus_{n\ge 0} \H(n)$$
with $\H(n)$ finite--dimensional and $\pi(L_0)\xi=(h+n)\xi$ for $\xi\in \H(n)$. Vectors in $\H(n)$ 
are said to have energy $h+n$. Necessarily $h$ and $c$ are real; and it is
easy to see that $c\ge 0$ and $h\ge 0$. Every positive energy unitary representation can be written as a direct sum of
irreducible representations which are uniquely determined by the lowest eigenvalue $h$ of $\pi(L_0)$. In this case
$\H(0)$ is one--dimensional and generated by a unit vector $\xi_h$. The character of a 
positive energy representation on $\H$ is defined as the formal sum
$$\chi_\H(q)= q^h \sum_{n\ge 0}d(n)\, q^{n},$$
where $d(n)={\rm dim}\, \H(n)$. We shall often write $L_n$ in place of $\pi(L_n)$, regarding $\H$ as a module over 
the Virasoro algebra. The above character is then formally
$$\sum_{n\ge 0} {\rm Tr}_{\H(n)}\, q^{L_0}= {\rm Tr}_\H \,q^{L_0}.$$ 
The work of Friedan--Qiu--Shenker {\bf [12]} 
shows that the irreducible positive energy unitary
representations are given precisely by a continuous series $c\ge 1$ and $h\ge 0$ together with a discrete series
$$c=1-{6\over m(m+1)},\,\, h=h_{p,q}(m)\equiv {(p(m+1)-qm)^2 -1\over 4m(m+1)},\,\,m\ge 2,\,\,1\le q\le p\le m-1.$$
The case $m=2$ corresponds to $c=0$ and had already been analysed by Gomes: 
only the trivial one--dimensional representation
can occur ($h=0$). Otherwise the representations are all on infinite--dimensional inner product spaces $L(c,h)$. The 
problem of computing the character $\chi_{c,h}(q)=\chi_\H(q)$ when $\H$ is the irreducible representation with parameters 
$(c,h)$ was first solved by Feigin and Fuchs {\bf [8]}. We outline their proof and some of its subsequent simplifications below.
It relies on a detailed knowledge of singular vectors and non--unitary representations.

The vector $\xi=\xi_h$ satisfies $L_n\xi=0$ for $n>0$ and $L_0\xi=h\xi$. 
The space $\H$ is therefore spanned by all
monomials $L_{-k}^{n_k} \cdots L_{-1}^{n_1}\xi$ with $n_i\ge 0$. On the other hand there is a universal 
module for the Virasoro algebra, the Verma module $M(c,h)$, which has as basis such vectors. Formally it can be defined 
as the the module induced from the one--dimensional representation $L_n=0$ ($n>0$), $L_0=h$ and $C=c$ of the 
Lie subalgebra spanned by $L_n$ ($n\ge 0$) and $C$. Its existence is standard and follows from the 
Poincar\'e--Birkhoff--Witt theorem applied to the universal enveloping algebra of the Virasoro algebra. (By that result the universal enveloping algebra has a basis
consisting of monomials of the above form post--multiplied by monomials 
$C^k L_0^{m_0} L_1^{m_1}\cdots$ 
with $k\ge 0$ and $n_i\ge 0$.) By definition $M=M(c,h)$ has positive energy and $d(n)={\rm dim}\,M(n)$
is given by the number of partitions $\P(n)$ of $n$. Thus
$$\chi_{M(c,h)}(q)= q^h\prod_{n\ge 1} (1-q^n)^{-1}.\eqno{(1)}$$
By the definition 
of $M(h,c)$, the irreducible representation $L(c,h)$ is the quotient of $M(c,h)$ by its unique maximal submodule.

The Verma module is in general not an inner product space, but for $c,h$ real, it does carry a unique Hermitian form
$(\xi,\eta)$ such that $(\xi_h,\xi_h)=1$ and 
$(L_n\xi,\eta)=(\xi,L_{-n}\eta)$, the so--called {\it Shapovalov form}. 
It is easy to see that the maximal submodule of $M(c,h)$ 
coincides with the kernel of the Shapovalov form. Thus for all real values of $c$ and $h$, the Shapovalov form defines
a non--degenerate Hermitian form on $\H=L(c,h)$. The representation is unitary precisely when this 
Hermitian form is positive--definite on $\H$, or equivalently, since the subspaces $\H(n)$ are orthogonal, on each $\H(n)$.

In {\bf [36]} we gave a direct method of determining the 
characters of the discrete series $0<c<1$ involving only unitary representations, which we summarise below. 
When $c>1$ the character formula 
is particularly simple because the Verma module is itself irreducible so the character formula is given by (1). This
remains true in the limiting case $c=1$ provided $h\ne m^2/4$ with $m\in {\Bbb Z}$. When $h=m^2/4$ or equivalently $j^2$ 
with $j$ a non--negative half--integer, the Verma module is no longer irreducible and the character formula takes the form
$$\chi_{L(1,j^2)}(q) =(q^{j^2}-q^{(j+1)^2})  \prod_{n\ge 1} (1-q^n)^{-1}.\eqno{(2)}$$
The main result of this article is a direct proof of this formula using only unitary representations; 
it is the counterpart of our previous proof for the discrete series
and relies only on a detailed knowledge of the boson--fermion
correspondence in conformal field theory. 

Before explaining the proof for $c=1$, it will be helpful to recall the method used for the discrete series $0<c<1$ 
which relies on the coset construction of Goddard--Kent--Olive {\bf [15]}. The discrete series at $c=1$ can be regarded as 
a limiting case of this construction. The complexification $\g$ of the Lie algebra of $SU(2)$ can be identified with
the $2\times 2$ complex matrices of trace zero. It is closed under taking adjoints and the Lie algebra of $SU(2)$ can be 
identified with the skew--adjoint matrices of trace zero. The affine Lie algebra $L\g$ is defined as the 
space of trigonometric polynomial maps into $\g$ with the pointwise bracket. 
It has generators $X_n=e^{in\theta}X$ ($X\in \g$) with commutation relations $[X_m,Y_n]=[X,Y]_{m+n}$. 
The Witt algebra acts naturally by differentiation, so that $[\ell_n,X_k]=-kX_{k+n}$. In particular $\ell_0$ defines a 
derivation $[d,X_n]=-nX_n$. The affine Kac--Moody algebra $\G$ is defined as a 
central extension of $L\g\rtimes {\Bbb C}d$ by ${\Bbb C}$. It has generators $X(n)$, $D$ and $C$ (central) 
satisfying the commutation relations 
$$[X(m),Y(n)]=[X,Y](m+n) +m\delta_{n+m,0} {\rm Tr}(XY)\cdot C,\,\,[D,X(n)]=-nX(n).$$
The Lie subalgebra generated by the $X(n)$ and $C$ is denoted ${\cal L}\g$: 
it is a central extension of $L\g$ by ${\Bbb C}$.
Again we can define positive energy unitary representations on an inner product space $\H$ by 
requiring that $C^*=C$ acts as a scalar $\ell I$, 
$D^*=D$. $X(n)^*=X^*(-n)$ and that $\H=\oplus_{n\ge 0} \H(n)$, where the $\H(n)$'s are finite--dimensional eigenspaces of $D$ corresponding to 
the eiegnvalue $n$. Any such representation decomposes as a direct sum of irreducible 
positive energy representations of the same form (possibly after adjusting $D$ by subtracting a positive scalar). 
For infinite--dimensional irreducible representations, $\ell$ must be a positive integer, called the {\it level}, 
and $\H(0)$ an irreducible 
representation of the Lie subalgebra $\g(0)=\g$. Each such $V_j$ is specified by a non--negative half--integer spin $j$
with ${\rm dim}\, V_j =2j+1$. At level $\ell$, the only spins that occur are those satisfying $0\le j\le \ell/2$.
Yhe Segal--Sugawara construction shows that in each irreducible positive energy representation operators
$L_n$ can be defined with 
$$L_0=D+ h\cdot I,\,\,\,\,h= {j^2+j\over 2(\ell+2)}$$
such that $[L_n,X(m)]=-mX(n+m)$, $L_n^*=L_{-n}$ and
$$[L_m,L_n]=(m-n)L_{n+m} +{m^3-m\over 12}c_\ell\cdot I$$
where $c_\ell=3\ell/(\ell+2)$.
The character of the representation is defined as the formal sum
$$\chi_{\ell,j}(q,g)=q^h \sum_{n\ge 0} q^n {\rm Tr}_{\H(n)} (g),$$
where $g\in SU(2)$. 
Since this only depends on the conjugacy class of $g$, as usual we may take 
$$g=\pmatrix{\zeta & 0\cr 0 & \zeta^{-1}},$$
with $\zeta|=1$. Note that formally 
$$\chi_{\ell,j}(q,g)={\rm Tr}_\H(gq^{L_0}).$$
These characters are given by the Weyl--Kac character formula (see {\bf [23], [31]})). 
In {\bf [36]} this formula was proved purely in terms of unitary representations together using the supersymmetric operators of Kazama--Suzuki; see also {\bf [26]}.

Let $\H_{\ell,j}$ denote the irreducible representation of level $\ell$ and spin $j$. 
The algebraic tensor product $\H=\H_{\ell,r}\otimes \H_{m,s}$ can be expanded as a direct sum
$$\H=\H_{\ell,r}\otimes \H_{m,s}=\bigoplus \H_{m+\ell,t}\otimes M_t,$$ 
where $t-r-s$ is an integer and $M_t$ is a multiplicity space
that can be 
identified with ${\rm Hom}_{{\cal L}\g} (\H_{m+\ell,t},\H)$. Now if $T\in   {\rm Hom}_{{\cal L}\g} (\H_{m+\ell,t},\H)$
so too is $L_n T - TL_n$ where the action on the right is by the Segal--Sugawara construction on $H_{m+\ell,t}$
and the action on the left on  $\H=\H_{\ell,r}\otimes \H_{m,s}$ is given 
the tensor product of the Segal--Sugawara constructions on $\H_{\ell,r}$ and $\H_{m,s}$. Hence the 
operator defined on the tensor product as 
$\pi_\ell(L_n)\otimes I + I \otimes \pi_m(L_n)-\pi_{\ell +m}(L_n)$ acts on the multiplicity space $M_t$. 
This is the Goddard--Kent--Olive construction. It gives a representation of the Virasoro algebra with central charge
$c=c_\ell+c_m - c_{\ell+m}$. In particular when $m=1$, we have a decomposition
$$\H=\H_{\ell,r}\otimes \H_{1,s}=\oplus \H_{\ell+1,t}\otimes M_t,$$
with $s=0$ or $1/2$, $r-t\in s +{\Bbb Z}$ and  
$$c={3 \ell\over \ell +2} + 1 - {3(\ell+1)\ell\over \ell+3}=1-{6\over (\ell+2)(\ell+3)}.$$
The character of the multiplicity space is determined directly by the Weyl--Kac characters at level $1$, $\ell$ and 
$\ell+1$ and this gives the lowest eigenvalue of $L_0$ in $M_t$, which turns out to be $h=h_{p,q}(m)$ with 
$p=2r+1$ and $q=2t+1$. The representation $L(c,h)$
therefore occurs as a summand of $M_t$. This shows that for this value of $c$, all the $h$ appearing in the list
of Friedan--Qiu--Shenker are indeed unitary. To complete the proof we just have to show that $M_t$ is irreducible.

This is accomplished by invoking the easy part of the Friedan--Qiu--Shenker theorem. 
Indeed as observed in {\bf [36]}, their proof of the classification theorem splits naturally into two parts. 
The first part shows that for $c$ fixed at $1-6/m(m+1)$, the values of $h$ have to be one of those in their list.
The proof requires the Kac determinant formula for $M=M(c,h)$ (see {\bf [22]}). This computes the determinant of 
$(v_i,v_j)$ where the $v_i$'s run across the canonical monomial basis of $M(n)$. If $c$ is fixed, 
this is a constant times a polynomial in $h$. On the other hand, 
since $\chi_{M_t}$ gives an upper bound for $\chi_{L(c,h)}$, it follows that $M=M(c,h)$ 
has a singular vector in $M(pq)$. The Kac determinant for 
$M(n)$ therefore must vanish when $n\ge pq$. 
This implies
that $h-h_{p,q}(c)$ divides the Kac determinant if $m$ is given by $c=1-6/m(m+1)$. By the structure of the Verma modules, 
it occurs with multiplicity at least ${\cal P}(n-pq)$ where ${\cal P}(n)$ is the partition function. Thus 
the Kac determinant is divisible by 
$$\prod_{p,q; 1\le pq \le n} (h-h_{p,q}(c))^{{\cal P}(n-pq)}.$$
Since it can be checked that this product has the same degree as the Kac determinant and they agree up to a constant 
in the highest power, the product above is proportional to the Kac determinant. We will use the same strategy of 
proof to establish a similar product formula, due to Feigin and Fuchs, for a 
polynomial arising elsewhere in the representation theory of the Virasoro algebra. 

Finally to prove that $M_t$ is irreducible, on the one had hand the Kac determinant formula 
gives a lower bound for the character of $\H=L(c,h)$. It follows that the $M_t(n)$ agrees with $
\H(n)$ for $n<m(m+1) -(m+1)p +mq=M$. On the other hand if $M_t$ is reducible it is a sum of representation
$L(c,h^\prime)$ from the list of Friedan--Qiu--Shenker. But it is elementary to verify that 
$h^\prime<M+h$ for all such $h^\prime$. It follows that $M_t$ must be irreducible.
\vskip .1in
The proof of the character formula for $c=1$ and $h=j^2$ proceeds similarly by indentifying a multiplicity 
space on which the Virasoro algebra acts and proving that it is irreducible. Graeme Segal {\bf [32]} gave a 
direct proof of the first case when $h=0$. Our proof is different in this case. It avoids the 
Segal--Sugawara construction and the Kac determinant formula, relying instead on various well known aspects 
of the fermion--boson correspondence as tools, including elementary parts of the theory of vertex algebras (see {\bf [6], [12], [14], [24]}).  

Indeed in this case a single complex fermion field is given by a set of operators $(e_n)$, $d$ 
acting on an inner product space
$\H$ subject to the anticommutation relations
$$e_me_n+e_ne_m=0,\,\, e_me_n^* + e_n^*e_m=\delta_{m+n,0}\cdot I, \,\, d=d^*,\,\, [d,e_n]=-(n+{1\over 2})e_n.$$
Slightly relaxing the positive energy condition, we require that $\H$ has a decomposition $\H=\oplus \H(r)$ where
$r$ runs over non--negative half integers with $\H(r)=\{\xi:d\xi=r\xi\}$ finite--dimensional. There is essentially
only one such irreducible representation on an inner product space 
${\cal F}_f$ generated by a vacuum vector $\Omega\in {\cal F}_f(0)$.
${\cal F}_f$ is called {\it fermionic Fock space} and can be identified with an exterior algebra. 

A single boson field is formed
by operators $(a_n)$, $d$ acting on an inner product space $\K$ and satisfying the commutation relations
$$a_ma_n-a_na_m=m\delta{m+n,0}\cdot I, \,\, a_n^*=a_{-n},\,\, d=d^*,\,\, [d,a_n]=-n a_n.$$
Again there is an essentially unique irreducible positive energy representation, this time with non--negative integer
eigenvalues for $d$. It is generated by $\Omega\in \K(0)$ satisfying $a_0\Omega=\mu \Omega$.A unitary charge operator
$U$ can be added to this bosonic system satisfying the additional relations
$$Ua_nU^*=a_{n}+\delta_{n,0}I,\,\,\,\, UdU^*= d+a_0 + {1\over 2} I.$$   
The ``discrete'' system 
$U$, $a_0$ and $d$ has an essentially unique positive energy representation with $d$ having non--negative
half integer eigenvalues. From this it follows that the charged boson system $(a_n,U,d)$ has an essentially unique
positive energy representation on an inner product space ${\cal F}_b$, obtained as a tensor product. ${\cal F}_b$ is
called {\it bosonic Fock space} and can be identified with a symmetric algebra tensored with the algebraic group algebra
of ${\Bbb Z}$. 

Boson--fermion duality is the statement that there is a natural unitary isomorphism between 
${\cal F}_f$ and ${\cal F}_b$ compatible with the operators, carrying $\Omega_f$ onto $\Omega_b$. Identifying the spaces,
the compatibility conditions state that $d$ should be preserved and that
$$[a_n, e_m] = e_{n+m},\,\, Ue_nU^*=e_{n+1}.$$
On the other hand it is straightforward 
to see that there are essentially unique operators on ${\cal F}_f$ satisfying these relations and that they act
irreducibly. In fact $U$ arises as an explicit shift operator on ${\cal F}_f$ and the $a_n$'s 
can be written as linear combination of operators of the form $e_ie_j^*$. These operators act irreducibly on ${\cal F}_f$,
so in this sense charged bosons can be constructed from fermions. To proceed in the other direction, i.e. to construct 
fermions from charged bosons, requires the introduction of vertex operators. For $m\in {\Bbb Z}$, these are defined as formal power series
in $z$ and $z^{-1}$ by the formula
$$\Phi_m(z)=U^m z^{-ma_0} \exp (\sum_{n<0} {mz^na_n\over n})\exp(\sum_{n>0} {mz^na_n\over n}).$$
The fundamental identities are then 
$$\Phi_1(z)=\sum e_n z^{-n-1},\,\,\Phi_{-1}(z)=\sum e_{-n}^* z^{-n}.$$
Identifying ${\cal F}_f$ and ${\cal F}_b$, which we write simply as ${\cal F}$, 
there are Virasoro operators $L_n$ with $c=1$ satisfying the 
covariance relations  
$$[L_n,e_m]=-(m+{1\over 2}(m+1))e_{m+n},\,\, [L_n,a_m]=-ma_{m+n}.\eqno{(3)}$$
The operators $L_n$ can be written either as linear combination of operators of the form $e_ie_j^*$ or of the form
$a_ia_j$. The identity (3) is a special case of the more general Fubini--Veneziano relation
$$[L_n,\Phi_m(z)] = z^{n+1}\Phi_m^\prime(z) +{m^2\over 2} (n+1)z^n\Phi_m(z).
\eqno{(4)}$$ 
Note that the operators $L_n$ commute with $a_0$ and so leave invariant its eigenspaces. The action 
on the non--zero eigenspaces turns out to be irreducible. The action on the zero eigenspace, however, 
is not irreducible. We will see that it is a direct sum of representations $L(1,m^2)$ with $m$ a non--negative
integer. Similarly the Virasoro algebra acts on the irreducible representation of $(a_n)$, $d$ with 
$a_0={1\over 4}I$, which it will turn out decomposes as a direct sum of representations $L(1,(m+{1\over2})^2)$.

A natural way to understand and study these decompositions is by passing to 
${\cal F}^{\otimes 2}={\cal F}\otimes {\cal F}$. The description of ${\cal F}$ as an exterior algebra shows that there is
a natural action on ${\cal F}^{\otimes 2}$ of the group $SU(2)$ and thus its Lie algebra. This action commutes
with the operator $d\otimes I + I \otimes d$. 
Conjugating the action of $a_n\otimes I$ by elements of $SU(2)$ then leads to operators $E_{ij}(n)$ with 
$E_{11}(n)=a_n$ satisfying $E_{ij}(n)^*=E_{ji}(-n)$ and
$$[X(m),Y(n)]=[X,Y](n+m) +m\delta_{m+n,0} {\rm Tr}\, XY,$$
if $X=\sum x_{ij}E_{ij}$ and $Y=\sum y_{ij}E_{ij}$. In particular taking $X\in \s$, the matrices with zero trace, we get
a positive energy representation of the Kac--Moody Lie algebra at level $1$. These operators commute with
the operators $A_n={1\over 2}(a_n\otimes I + I\otimes a_n)$. In particular the $A_n$'s commute with the operators
$B_m=H(m)$, $E(m)$ and $F(m)$ where $H={1\over 2}(E_{11}-E_{22})$, $E=E_{12}$ and $F=E_{21}$ is the standard basis of $\s$.
These operators are related to the vertex operators by the formulas
$$\Phi_1(z)\otimes \Phi_{-1}(z)=\sum E(n) z^{-n-1}, \,\, \Phi_{-1}(z)\otimes \Phi_1(z) =\sum F(n)z^{-n-1},\eqno{(5)}$$
identities originally due to Frenkel--Kac {\bf [10]} and Segal {\bf [32]}.

The space ${\cal F}^{\otimes 2}$ decomposes as a sum of two components 
$${\cal F}^{\otimes 2}=\H_{0}\otimes M_0\oplus \H_{1/2}\otimes M_{1/2},$$
with the affine Kac--Moody algebra acting irreducibly on $\H_j$ and the bosonic operator irreducibly on the multiplicity
spaces $M_j$. The action of the Virasoro algebra on ${\cal F}^{\otimes 2}$ preserves this decomposition. Its actions on the
tensor factors coincides with the natural action associated with $(B_n)$ and $(A_n)$. The action on $H_j$ 
commutes with the action of $SU(2)$ and satisifies $[L_n,X(m)]=-mX(n+m)$. Each space $\H_j$ may be further
decomposed according to the action of $SU(2)$
$$\H_j=\bigoplus_{i\ge 0,\, i-j\in {\Bbb Z}} V_i\otimes \K_{ij}.$$
The operators $L_n$ commute with $SU(2)$, so preserve the multiplicity spaces $\K_{ij}$, each of which contains a copy of 
$L(1,i^2)$ generated by a singular vector in the lowest energy space.  This action on the multiplicity spaces can be
regarded as a limiting case of the coset construction of Goddard--Kent--Olive. 
The singular vectors were described explicitly by Goldstone {\bf [21]} 
and we shall refer to them as
{\it Goldstone vectors}. Goldstone's
formulas were shown by Segal {\bf [32]} to be a direct consequence of the vertex operator formulas in (5). Although the formulas
can be developed without any advanced theory of symmetric functions {\bf [27], [29]}, the simplest way to describe them involves the
combinatorics of the Weyl character formula for $U(n)$ {\bf [40]}. In fact setting $X_n=a_{-n}$ 
and regarding these as the symmetric functions with signature $(n,0,0,0,\dots)$, the Goldstone 
vectors correspond to the symmetric functions $X_f$ with signature $(k+m,k+m,\dots,k+m,0,0,\dots)$, where $k+m$ appears
$k$ times applied to a singular vector for the $a_n$ with $n>0$. (If we took $X^\prime_k=\sum_i z_i^k$, then $X^\prime_f$ would just be $\det\, z_j^{f_i+n-i}/ \det\, z_j^{n-i}$, the character of the irreducible representation $V_f$ of $U(n)$ with signature $f$.)  Goldstone conjectured that
these vectors were the only singular vectors in $\H_j$ so that consequently the Virasoro algebra acted irreducibly
on the multiplicity spaces. Since ${\rm ch}\,\K_{ij}=(q^{i^2}-q^{(i+1)^2})\cdot \varphi(q)$, this establishes the character
formulas for $L(1,i^2)$. 

To prove this directly, it suffices to show that in $\H_j$ the Goldstone vectors are the only singular vectors 
in a fixed eigenspace $\K$ of $H(0)$. 
It is easy to check that $\K$ is an irreducible representation of the $H(n)$'s
and that, up to scalar multiples, there is at most one singular vector at any fixed energy level in $\K$. Moreover as module over $\v$, $\K$ is a direct sum of 
irreducible representations, with one component for each singular vector. 
If there were a component $\K^\prime$ isomorphic to $L(1,p)$ with $p$ not of the form $(j+k)^2$, then by irreducibility, for some component 
$\K^{\prime\prime}$ generated by a Goldstone vector, $H(a)\K^{\prime\prime}$ would have to have a non--zero projection on $\K^\prime$. But 
if $P^\prime$ and $P^{\prime\prime}$ denotes the orthogonal projections 
onto $\K^\prime$ and $\K^{\prime\prime}$ , then 
$$\Psi(z)=P^\prime H(z)P^{\prime\prime}$$ 
would give a 
formal power series of operators  from $\K^{\prime\prime}$ to $\K^{\prime}$ satisfying
$$[L_k,\Psi(z)] = z^{k+1}d\Psi^\prime(z).\eqno{(6)}$$
Such family of operators is a special case of what is called a {\it primary field}. 
Writing $\Psi(z) =\sum \psi(n)z^{-n}$, (6) can be rewritten
$$[L_k,\psi(n)]=-n \psi(n+k). \eqno{(7)}$$ 
This prompts the introduction of the {\it density modules} $V_{\lambda,\mu}$ giving the natural representations of
the Witt algebra on expressions of the form $f(\theta)e^{i\mu\theta}(d\theta)^\lambda$, with $f$ a trigonometric 
polynomial. Multiplication and integration over the circle
gives a natural pairing with $V_{1-\lambda,-\mu}$. A primary field, for appropriate $\lambda$ and $\mu$, then defined to be
a linear map $\K^{\prime\prime}\otimes V_{\lambda,\mu}\rightarrow \K^\prime$ 
commuting with the action of the Virasoro algebra. By duality
it is the same as an equivariant map $\K^{\prime\prime}\rightarrow \K^\prime\otimes V_{\lambda^\prime,\mu^\prime}$, 
where $\lambda^\prime=1-\lambda$ 
and $\mu^\prime=-\mu$. The formula (6) is then replaced by the more general relation
$$[L_k,\Psi(z)] = z^{k+1}d\Psi^\prime(z)+\lambda z^k\Psi(z).\eqno{(8)}$$
Returning to (6) and (7), 
the action there is simply on functions, so that $\lambda=0=\mu$. But then
$V^\prime_{0,0}=V_{1,0}$, the space of differentials $f(\theta)\,d\theta$. 
On the other hand if $\K^{\prime\prime}$ is isomorphic to $L(1,i^2)$ then it is the quotient of a Verma module
$M(1,i^2)$ by a submodule containing at least one singular vector $w$ of energy $(i+1)^2$. If $v$ is the
lowest energy vector of $M(1,i^2)$, then $w$ has the form $Pv$ for some non--commuting polynomial in $L_{-k}$ ($k\ge 1$).
Taking the component of the lowest energy vector of $K_1$, it is then easy to see that a necessary condition for the
existence of a primary field $\Psi(z)$ is that the action of $P$ on one of the 
standard basis vectors of $V_{\lambda^\prime,\mu^\prime}$ must be zero. But the action of the Witt algebra is given in this basis
as shift operators weighted by polynomials, so the condition is equivalent to the vanishing of a polynomial. A
formula for this polynomial was given without proof by Feigin--Fuchs {\bf [8]}. In this particular case where $c=1$, 
sufficiently many factors of the polynomial can be produced by explicitly exhibiting enough
 primary fields to determine completely the polynomial. (These are constructed as operators between different $L(1,i^2)$ by compressing $\Phi_a(z)\otimes \Phi_b(z)$.) 
Irreducibility then follows by noting that this polynomial does not vanish for the 
hypothetical field $\Psi(z)$ constructed above. 

The proof that there are sufficiently many primary fields is essentially the first step in showing that under fusion 
the representations $L(1,i^2)$ for $i$ a non--negative half integer behave exactly like the irreducible representations
$V_i$ of $SU(2)$ under tensor product (see {\bf [38]}). Now if $\Psi(z)$ is a primary field
and $L_{-1}w=0$, we see that $F(z)=\Psi(z)w$ satisfies the so--called ``equation of motion''
$${dF\over dz}=L_{-1}F,$$ 
so that 
$$F(z)=e^{zL_{-1}} \Psi(0)\Omega.$$
This identity allows the existence of sufficiently many
primary fields to be reduced to checking that \break $L_{1}^{|f|}X_fv\ne 0$ where $v$ is a singular vector of the $H(n)$'s 
for $n>0$ with $H(0)v=pv$ with $p$ a poitive integer.  But it is easy to see that $(|f|!)^{-1}L_{1}^{|f|}X_fv=a v$ for some constant $a$. Like $X_f$, the value of $a$ is given 
by a determinant, 
but this time of binomial coefficients. It can be calculated as an explicit product using either the Weyl dimension 
formula or an elementary matrix computation.  In particular, in the case of interest, it does not vanish: indeed remarkably $a$ equals the dimension of the irreducible representation
$V_f$ of $U(p)$ when this makes sense and vanishes otherwise. This can be summarised by saying that, if the Goldstone vectors are given formally by formulas relating to the Weyl character formula, then  their non--degeneracy properties follow from the 
corresponding Weyl dimension formula. 

Sections 2 to 8 give a self--contained step--by--step account 
of this method of establishing the character formula for $L(1,j^2)$.
Appendix A contains a alternative direct verification of the Fubini--Veneziano relations for vertex operators. 
In Appendices B and C, we give an alternative approach to the product formula of Feigin--Fuchs using
an elegant algorithm of Bauer, Di Francesco, Itzykson and Zuber ({\bf [3], [7]}) for determining the polynomial $P$, related to an earlier formula of Benoit and St Aubin {\bf [4]}: 
we give a simplification of their treatment in Appendix~A. The most general Feigin--Fuchs formula can be proved 
using the coset construction for primary fields for the 
unitary discrete series with $0<c<1$. Indeed the primary fields for the affine Kac--Moody algebra $\s$ can be constructed
explicitly and then used to contruct primary fields for the Virasoro algebra 
{\bf [28], [39]}. In this case
the irreducible representations for $c=1-6/m(m+1)$ with $m$ large are indexed by non--negative half integers $(i,j)$ with
$i-j$ and integer. The corresponding representation has $h=h_{p,q}(m)$ with $p=2i+1$ and $q=2j+1$. The singular vector
at level $pq$ acts on the density module as a multiple of 
$\prod(h- h_{p^\prime,q^\prime}(m))$ where $V_{i}\otimes V_{i_1}=\bigoplus V_{i^\prime}$ and 
$V_{j}\otimes V_{j_1}=\bigoplus V_{j^\prime}$. Other methods of proving the Feigin--Fuchs formula are given in {\bf [7]} and  {\bf [13]}. 

In Section 9 we desribe a simple method for determining the character 
formula for $L(1,j^2)$ 
using only the Kac determinant formula and the Jantzen filtration {\bf [20]}.
It is generalised in Section~10 to the case $0<c<1$ and gives a different method of 
proof to that of Astashkevich {\bf [1]}.  Let $A(x)=\sum_{i\ge 0} A_i x^i$ be an 
analytic family of non--negative self--adjoint matrices defined 
for $x$ real and small on $V={\Bbb R}^n$ or ${\Bbb C}^n$. We define a filtration $V^{(i)}$ ($i\ge 0$) by
$V^{(0)}=V$ and for $m\ge 1$ 
$$V^{(m)}=\bigcap_{i=0}^{m-1}{\rm ker}\, A_i.$$
In the examples $A(x)$ will be invertible for $x\ne 0$ in a neighbourhood of $0$.  It follows 
that $V^{(m)}=(0)$ for $m$ sufficiently large. The fundamental identity, which is easy to prove, is that the order
of $0$ as a root of ${\rm det}\, A(x)$ equals
$$\sum_{i\ge 1}{\rm dim}\, V^{(i)}.$$
Taking the standard inner product on $V$, we can define a Hermitian form on $V$ by $(v_1,v_2)_x=(A(x)v_1,v_2)$.
This filtration has an important functorial property. Suppose we are given two such space $(V, A(x))$ and $(W,B(x))$ and in addition analytic maps $X(x):V\rightarrow W$, $Y(x):W\rightarrow V$
such that 
$$(X(x)v,w)_x=(v,Y(x)w)_x.$$
It then follows easily that $X(0)$ and $Y(0)$ induce maps between $V^{(i)}$ and $W^{(i)}$ preserving the canonical Hermitian 
forms induced by $A_i$ and $B_i$. 

Now let $M=M(1+x,j^2)$. Then $V=M(n)$ is independent of $x$ with a canonical basis, 
so can be identified with 
${\Bbb C}^{{\cal P}(n)}$.  We consider the Jantzen filtration for $M=M(1,j^2)$ associated with the parameter $x$ in 
$c=1+x$. The last property above implies that each $M^{(i)}\equiv\oplus M(n)^{(i)}$ is a submodule for the 
Virasoro algebra and comes with an invariant Hermitian form. By the Kac determinant formula we can explicitly 
compute $\sum_{i\ge 1} {\rm ch}\,M^{(i)}$. Now the determinant formula or the coset construction 
above show that $M(1,j^2)$ has a filtration by Verma modules
$M(1,(j+k)^2)$ for $k\ge 1$. Let $M^{(n_k)}$ be the last term in the Jantzen filtration containing $M(1,(j+k)^2)$; it turns out that $(n_k)$ is strictly increasing (we set $n_0=0$). But a lower bound for the character sum above is given by
$\sum_{k\ge 1} (n_k-n_{k-1}) {\rm ch}\, M(1,(j+k)^2$ and hence $\sum_{k\ge 1} {\rm ch}\, M(1,(j+k)^2$. 
However this sum actually equals the character sum. Hence $M^{(k)}=M(1,(j+k)^2)$ and hence
$L(1,j^2)=\varphi(q)\cdot (q^{j^2} - q^{(j+1)^2})$.

The proof for the discrete series is similar but slightly but more intricate. 
Let $c=1-6/m(m+1)$ and $M=M(c(m),h_{p,q}(m))$ with $1\le q\le p\le m-1$. 
Again the Kac determinant formula 
produces a sequence of singular vectors in $M$ which generate a series of Verma submodules $A_i,B_i$ for $i\ge 1$ such that
$A_i$ and $B_j$ both contain $A_k$ and $B_k$ for $k>i$. Thus $A_k\cap B_k \supseteq A_{k+1} + B_{k+1}$. 
Again, this time using $x=h-h_{p,q}(m)$ as a parameter, we can compute the character sum 
$\sum_{i\ge 1} {\rm ch}\,M^{(i)}$ as well as that for $A_1$ corresponding to the
singular vector with energy $h+pq$. Analogous inequalities to those above prove that 
$A_k\cap B_k =A_{k+1} +B_{k+1}$ for $k\ge 1$ and that $M^{(k)}=A_k+B_k$ for $k\ge 1$. Using the isomorphisms
$A_k\oplus B_k/(A_{k+1}+B_{k+1})=A_k+B_k$, the character formula for $L(c(m),h_{p,q}(m))=M/(A_1+B_1)$ follows.  

Although we have not checked this, 
it seems likely that variants of the two proofs above can be used to 
simplify substantially the proof of Feigin and Fuchs. Such a simplification had been given by Astashkevich {\bf [1]}
for all cases except III${}^{00}_\pm$. The method of Feigin--Fuchs in this case was considerably 
more complicated than what is 
required for the discrete series $c<1$ (the case III${}_-$) 
and invokes the Riemann--Roch theorem for holomorphic vector bundles on the sphere. 
\vskip .1in
We end this overview by briefly describing the Feigin and Fuchs' original proof of the character formula for
$L(1,j^2)$ in {\bf [8]}---their case III$_-^{00}$.  In addition to a knowledge of the discrete 
series characters and their Jantzen filtrations, their proof requires the fact that the formula for the singular vector 
is a polynomial in $t$ of degree $2j+1$, with constant term $L_{-1}^{2j+1}$ and leading coefficient $[(2j)!]^2 L_{-2j-1}$. 
This was proved indirectly by Astashkevich and Fuchs {\bf [2]} for $M(1,j^2)$. In this case we show in Appendix~A that it is an easy 
consequence of the formula of Bauer, Di Francesco, Itzyksohn and Zuber {\bf [3]} for the singular vector of $M(h_{2j+1,1}(t),
13 -6t -6t^{-1})$. 
 
Instead of the Riemann--Roch theorem, all that is needed is the classical theory of holomorphic vector bundles 
on the Riemann sphere. Recall
that according to the Grothendieck--Birkhoff theorem, any 
such bundle can be written uniquely as a sum of line bundles, each classified by their degree (the sum with 
multiplicitied of a the degrees of poles and zeros of any meromorphic section). Grothendieck's 
original proof {\bf [19]} was algebraic--geometric using the language of divisors, although it was soon realized that the theorem 
followed from the Birkhoff factorization theorem (see {\bf [16], [31], [38]}). Indeed it is known {\bf [9]}
that a rank $n$ holomorphic vector bundle is trivial when restricted to either of the open discs 
$z<R$ and $|z|>r$ with $r<1<R$ (see Appendix D for a short analytic proof). But then it is specified by a holomorphic map
$f$ from the annulus $r<|z|<R$ into $GL_n({\Bbb C})$. By Birkhoff's factorization theorem the 
restriction of $f$ to the circle $|z|=1$ can be written
$f(z)=f_-(z) D(z) f_+(z)$ where $D(z)$ is a matrix with entries $z^{n_i}$ on the diagonal, $f_\pm(z)$ 
is a function holomorphic on $D_\pm(z)=\{z:|z|^{\pm 1} <1\}$ and smooth on its closure. 
The matrix $D(z)$ is uniquely determined up to the order of its entries. Since $f$ extends to $r<|z|<R$ 
it follows immediately from the Schwarz reflection principle that $f_+(z)$ and $f_-(z)$ extend to 
homolorphic function on $|z|<R$ and $|z|>r$. But then the clutching function on the annulus can be replaced by $D(z)$
which corresponds to a direct sum of holomorphic line bundles. (Each of these can in turn
be understood in 
terms of holomorphic sections of a homogeneous line bundle for the group $G=SU(2)$, 
identifying the Riemann sphere $S$ with $G/T$, where $T$ 
is the subgroup of diagonal matrices.)

Given a representation such as $L(1,j^2)$, the method of Feigin and Fuchs proceeds as follows. Let $r=2j+1$. 
At each point of the curve $C$ given by $\varphi_{r,1}(c,h)=0$, there is a singular vector of energy $2j+1$, 
in the Verma module $M=M(c,h)$. The action of the operators $L_{-n}$ for $n>0$ is independent of $c$ and $h$ so the
Verma module can be identified along the curve. The curve on the other hand can be parametrised by
$c=13-6t -6t^{-1}$ and $h=(j^2+j)t-j$ with $t\in {\Bbb C}^*$. The singular vector 
$w(t)\in M(2j+1)$ depends polynomially on $t$ (see Appendix B) 
and thus defines a line bundle on the Riemann sphere: 
its degree $d$ is determined by the degree of the polynomial, 
since the constant term is non--vanishing. Let $W_t$ be the 
Verma submodule of $M$ generated by $w(t)$. At any fixed energy level $k$ it is a sum of line bundles of degree $d$.
So $M(k)/W_t(k)$ defines a holomorphic vector bundle $E$ over the Riemann sphere, 
which is a sum of line bundles with the sum of the degrees determined. 
But if the rank of $E$ is $m$, 
$F=\lambda^m E$ is a line bundle of known degree. Now the Kac determinant formula implies that 
generically $M/W_t$ is irreducible on the curve $\varphi_{r,1}(c,h)=0$. 
The Kac determinant in the quotient produces a meromorphic section of $F^*$, 
the line bundle dual to $F$.  The degree of $F^*$ is known. It can also be 
computed in terms of the order of poles and zeros of the 
meromorphic section. The contributions from $0$ and $\infty$ can be computed directly, since there are two explicit bases 
independent of $t$ in $V$ complementary to $W_t$ for $t\ne 0,\infty$.  The 
computation reduces to that of a prinicipal minor where the diagonal terms dominate. Only real points given
by the intersection of the curve with other curves $\varphi_{p,q}(c,h)=0$ with $pq\le N$ give other contributions. 
For $c\ne 1$, these correspond to Verma modules of type III$_-$ 
having the same structure of 
singular vectors as the discrete series
with $0<c<1$. These  will have a singular vector of type $(r,1)$ amongst the other 
singular vectors (possibly of lower energy).
Since the structure of these modules is known explicitly, 
the Jantzen filtration on $M/W_t$ corresponding to the (real part of the) 
curve $C$ can be computed and hence the contribution from these points. Adding up these contributions, it 
then follows that the contribution from $(1,j^2)$ is $0$, as required.     
 
\vskip .1in
\noindent \bf 2. Constructions of positive energy representations. \rm We define the bosonic 
algebra $\a$ to be the Lie algebra with basis $a_n$ ($n\in {\Bbb Z}$), $c$ and $d$, with non--zero brackets 
$$[a_m,a_n]=m\delta_{m+n,0}\cdot c,\,\,\,[d,a_n]=-n a_n.$$ 
(Thus $c$ is central.) A unitary positive energy representation of 
$\a$ consists of an inner product space ${\cal H}$ which can be written 
as the orthogonal algebraic direct sum of finite--dimensional spaces ${\cal H}(n)$ ($n\ge 0$) 
$$\H=\bigoplus_{n\ge 0} {\cal H}(n)$$
on which $\a$ acts such that
$d\xi=(n+h)\xi$ for $\xi\in {\cal H}(n)$, with $h\in {\Bbb R}$, $c=I$ and $a_n^*=a_{-n}$. 
$d$ is called the {\it energy operator} and an element of ${\cal H}(n)$ is said to have energy $h+n$. Note that replacing 
$d$ by $d+tI$ gives another essentially equivalent representation. 

Note that if an operator $a$ satisfies the canonical commutation relations 
$a^*a-aa^*=\mu I$ and $a^*v=0$, then an easy induction argument shows that 
$$(a^m v,a^n v)=\delta_{m,n} m! \mu^m.\eqno{(1)}$$   
Now let $V={\Bbb C}[x_1,x_2,\dots]$ and for $n>0$ define operators $a_{-n}p=x_n p$, $a_np=n\partial_{x_n}p$ and 
$a_0 p=\lambda p$. We define ${\rm deg}\,x_n=n$ and let $D$ be the operator which multiplies a monomial by its degree, 
so that 
$$D=\sum_{n\ge 1} n \cdot x_n\partial_{x_n}.\eqno{(2)}$$
This gives a representation of the Lie algebra $\a$. Since $p(x_1,x_2,\dots)=p(a_1,a_2,\dots)1$, repeated applications
of (1) show that it is unitary for the inner product making monomials orthogonal with
$$\|x_1^{m_1}x_2^{m_2}\dots \|^2=\prod m_n! \,n^{m_n}.$$
It is irreducible because any non--trivial $\a$--submodule must contain a non--zero vector $p$ of 
lowest energy. Thus $a_n p=0$ for all $n>0$, i.e. $\partial_{x_n}p=0$ for all $n>0$. But then $p$ is a constant and the
submodule must be the whole of $V$. Conversely if $W$ is any positive energy unitary representation, 
then by (1) any lowest energy 
vector unit $w$ will satisfy $(p(a_{-1},a_{-2},\dots)w,q(a_{-1},a_{-2},\dots)w)=(p,q)$, 
the inner product on $V$ constructed above. Let $W_1$ 
be the $\a$--module generated by $w$. Adjusting $a_0$ and $d$ by constants if necessary, it follows from (1) that 
$W_1$ is unitarily equivalent to $V$ by the map $Up=p(a)w$. The orthogonal complement of $W_1$ is also a positive energy 
unitary representation so likewise contains an irreducible submodule $W_2$ of the same form. Continuing in this way, 
the positive energy condition shows that we can write $W=\bigoplus_{m\ge 1} W_m$, an orthogonal direct sum where each $W_m$ 
is isomorphic to $V$ but with $a_0$ acting as $\lambda_m I$ and $d$ as $D+h_mI$. 
Thus $\lambda_m$ and $h_m$ are the eigenvalues of $a_0$ and $d$ on the lowest energy vector $w_m$: 
clearly this determines the representation uniquely. Representations of the Virasoro algebra can be define on these
spaces by the generalising the formula (1) as follows:
$$L_0={1\over 2} a_0^2 +\sum a_{-n}a_n,\,\, L_n={1\over 2}\sum_{r+s=n}a_ra_s \,(n\ne 0).$$
Note that in the second formula the operators $a_r$ and $a_s$ commute. 
By definition $L_n$ changes energy by $-n$. It is easy to check that $[L_n,a_m]=-ma_{m+n}$ and hence that
$A_{m,n}=[L_m,L_n]-(m-n)L_{m+n}$ commutes with all the $a_j's$. It leaves each $W_k$ invariant. By uniqueness of the lowest
energy vector in $W_k$, it must act as a scalar $\nu_{m,n}$ and therefore preserve energy. Thus $A_{n,m}=0$ unless $n=-m$.
Otherwise for $n>0$
$$\eqalign{(([L_n,L_{-n}]-2nL_0)w_k,w_k)&=-n(a_0^2w_k.w_k) +{1\over 2}(L_n\sum_{r+s=-n} a_ra_s w_k,w_k) \cr
&=-n(a_0^2w_k,w_k)+{1\over 2}\sum_{r+s=-n} -r(a_{n+r}a_sw_k,w_k) 
-s(a_ra_{n+s}w_k,w_k)\cr
&={1\over 2}\sum_{t=1}^n t(n-t)\cr
&=(n^3-n)/6.\cr}$$
Thus 
$$[L_m,L_n]=(m-n)L_{m+n} + \delta_{m+n,0} {m^3-m\over 12}I.$$
\vskip .1in
We now make similar constructions using the canonical anticommutation relations. In this case the notion of positive
energy representation has to be slightly generalised to allow the eigenvalues of the energy operator to have
the form $n+h$ with $n\in {1\over 2}{\Bbb Z}$. We define $V$ to be the inner product space 
spanned by the mutually orthogonal unit vectors 
$$e_{i_1}\wedge e_{i_2} \wedge e_{i_2} \wedge \cdots,\eqno{(3)}$$ 
with $i_k\in{\Bbb Z}$, $i_1<i_2<i_3<\cdots$ and $i_{k+1}=i_k+1$ for $k$ 
sufficient large. We define the vacuum vector by $\Omega=e_{0}\wedge e_{1}\wedge e_{2} \wedge \cdots$ and make the 
$e_i$'s act by exterior mutiplication. It is easily checked that the $e_j^*$'s act by interior multiplication and that
$$e_me_n^*+e_n^*e_m=\delta_{m,n}I,\,\, e_me_n+e_ne_m=0.$$
The space $V$ can be canonically identified with the exterior algebra in the anticommuting variables 
$e_{-n}$ for $n>0$ and $e_n^*$ for $n\ge 0$. 
We can define the energy and charge of a basis vector in (3) as follows. Firstly it can be regarded as a monomial 
in the exterior algebra in the $e_{-n}$'s ($n>0$) and the $e_{n}^*$ ($n\ge 0$). A term $e_{-n}$ or $e_{n}^*$ 
contributes $n-1/2$ to the energy and $+1$ and $-1$ to the charge. The vector 
$$\Omega=\Omega_0 =e_0\wedge e_1 \wedge e_2 \wedge \cdots$$
has energy $0$ and charge $0$. Let $d$ and $a$ be the operators multiplying the basis vectors by their energy and 
charge respectively. Thus $a^*=a$ and $d^*=d$. Morevoer $[d,e_i]=-(i+1/2)e_i$ and $[a,e_i]=e_i$. Taking adjoints we get
$[d,e_i^*]=(i+1/2)e_i^*$ and $[a,e_i^*]=-e_i^*$. Alternatively a basis vector can be written
$$v=e_{m_1} \wedge \cdots \wedge e_{m_j} \wedge e_{j+k}\wedge e_{j+k+1}\wedge \cdots,$$
where $m_1<m_2<\cdots <m_j< k+j$. Thus $m_i<k+i$. The vector will have charge $-k$ and energy $\sum_{i=1}^j (k+1/2-i-m_i)$.
Indeed the vector
$$\Omega_k=e_k\wedge e_{k+1}\wedge \cdots$$
has charge $-k$ and energy $k^2/2$. The operators $e_ae_b^*$ commute with $a$ so preserve charge.
Since $v$ is obtained from $\Omega_k$ by successively applying the operators $e_{m_i}e_{k-i}^*$, which raises the energy
by $k-i-m_i$, it follows that $v$ has energy $\sum_{i=1}^j (k-i-m_i)$. 
We define
$$a_k=\sum_{n-m=k} e_n e_m^*,\,\, a_0=\sum_{n\ge 0} e_{-n}e_{-n}^* -\sum_{n>0} e_n^*e_n$$
and
$$L_k^\prime=\sum_{n-m=k} -(m+1/2+k/2) e_ne_m^*,\,\, L_0^\prime= \sum_{n> 0}(n-{1\over 2}) e_{-n}e^*_{-n} + \sum_{n\ge 0} 
(n+{1\over 2})e_n^*e_n.$$
It is straightforward to check that $a_n^*=a_{-n}$, $(L^\prime_n)^*=L^\prime_{-n}$ and
$$[L_n^\prime,e_m]=-(m+{n\over 2}+{1\over 2})e_{n+m},\,\, [L_n^\prime,e_m^*]=(m-{n\over 2}-{1\over 2})e_{m-n}^*,\,\, 
[a_n,e_m]=e_{n+m},\,\,[a_n,e_m^*]=-e_{m-n}^*.$$
The operators $d-L_0^\prime$ and $a-a_0$ commute with the $e_n$'s and $e_n^*$'s and annihilate $\Omega$. By cyclicity 
of $\Omega$, it follows that $d=L_0^\prime$ and $a=a_0$. The operator $[a_n,a_m]$ commute with  
the $e_k$'s and $e_k^*$'s
and raise energy by $-n-m$. If $n+m<0$ it must annihilate $\Omega$ and hence vanish. Since $[a_n,a_m]^*=[a_{-m},a_{-n}]$.
it follows that $[a_n,a_m]=0$ if $n+m\ne 0$. Now up to a scalar multiple $\Omega$ is the unique vector of energy $0$, 
so that if $n>0$, $[a_n,a_{-n}]\Omega=\lambda_n\Omega$; by cyclicity $[a_m,a_n]=\lambda_n\delta_{m+n,0}I$. To evaluate
$\lambda_n$ we have
$$\lambda_n=([a_n,a_{-n}\Omega,\Omega)=-\|a_{-n}\Omega\|^2=n.$$
Thus 
$$[a_m,a_n]=m\delta_{n+m,0}I.$$
Similarly 
$[L^\prime_m,L^\prime_n]-(m-n)L^\prime_{m+n}$ commutes with  the $e_k$'s and $e_k^*$'s
and raise energy by $-n-m$. Using adjoints as before, it must be $0$ if $n+m\ne 0$ and a scalar $\mu_m$ if $m+n=0$.
Clearly $\lambda_{-m}=-\lambda_m$ using adjoint and if $n>0$
$$\eqalign{\lambda_n&=([L^\prime_n,L^\prime_{-n}]\Omega,\Omega)\cr
&=\|L^\prime_{-n}\Omega\|^2 \cr
&= \sum_{m=1}^n (m-n/2 -1/2)^2\cr 
&={n(n+1)(2n+1)\over 6} +{n(n+1)^2\over 4}-{(n+1)^2n\over 2}\cr
&={1\over 12}n(n+1)[4n + 2 + 3(n+1) - 6(n+1)]={1\over 12} (n^3-n).\cr}$$
Hence we have
$$[L^\prime_m,L^\prime_n]=(m-n)L^\prime_{m+n} + \delta_{m+n,0} {m^3-m\over 12} I.$$
\vskip .1in

\noindent \bf Lemma (Gomes). \it 
The only irreducible positive energy unitary representation of the Virasoro algebra with $c=0$ is the
trivial representation.
\vskip .05in
\noindent \bf Proof. \rm Let $\Omega$ be a lowest
energy unit vector with $L_0\Omega=h\Omega$ and $L_n\Omega=0$ for $n>0$. Let
$\xi_1=L_{-2n}\Omega$ and $\xi_2=L_{-n}^2\Omega)$. Then $\|\xi_1\|^2= 4nh$, $\|\xi_2\|^2=4n^2h(2h +n)$ and 
$(\xi_1,\xi_2)= 6n^2h$. Hence 
$$\det (\xi_i,\xi_j) = 16n^3h^2(2h +n) -36n^4h^2=4n^3h^2(8h-5n).$$
To be non--negative for all $n$ we must have $h=0$. 
But then unitarity implies that $L_{-n}\Omega=0$ for all $n>0$. 
\vskip .1in
We have constructed the positive energy representation of the fermionic operators 
$e_n$ and $e_n^*$ on ${\cal F}$ with vacuum vectors $\Omega$. We have constructed operators, bilinear in fermions,
$a_n$ and $L^\prime_n$ satisfying $[a_m,a_n]=m\delta_{m+n,0} I$, $[L^\prime_n,e_m]=-m e_{n+m}$ and
$[L^\prime_n,a_m]=-ma_{m+n}$. The $L_n^\prime$'s define a representation of the Virasoro algebra with $c=1$. 
On the other hand we may also define operators $L_n$, bilinear in the bosons $a_i$, which give another representation
of the Virasoro algebra with $c=1$ satisfying $[L_n,a_m]=-m a_{m+n}$. The coset operators $L^{\prime\prime}_n=L^\prime_n-L_n$ 
define a positive energy representation of the Virasoro algebra with $c=0$ on the multiplicity spaces
$M_i=\{\xi|a_0\xi=i\xi,\,\, a_n\xi=0 \,\,(n>0)\}$. So $L_n=L^\prime_n$ on $M_i$ and hence everywhere.

\vskip .1in
\noindent \bf Corollary. \it $L_n^\prime= L_n$. 
\vskip .1in

\rm We now generalize this contruction to produce representations of the Virasoro algebra with $c=1$ and $h=j^2$ 
with $j$ a non--negative half--integer. The space ${\cal F}$ is canonically 
${\Bbb Z}_2$--graded and the grading is compatible with its identification with an exterior algebra. The operators $e_i$ and $e_i^*$ act as 
odd operators. On the ${\Bbb Z}_2$--graded tensor product ${\cal F}_2={\cal F}\otimes {\cal F}$ 
there are operators $e_i^{(1)}=e_i\otimes I$ and $e_i^{(2)}=I\otimes e_i$. The space ${\cal F}_2$ can be identified with
 the exterior algebra on the anticommuting variables $e_{-n}^{(j)}$ ($n>0$) and $(e_n^{(j)})^*$ ($n\ge 0$). 
If for $v=(\alpha,\beta)\in {\Bbb C}^2$ we set $v_i=\alpha e_n^{(1)}+\beta e_n^{(2)}$ and $v_i^*= 
\overline{\alpha} (e_n^{(1)})^*+\overline{\beta} (e_n^{(2)})^*$, then the natural actions of $G=U(2)$ on $V$ and $V^*$ 
(the conjugate space) extend to a unitary action $\pi$ on the exterior algebra and hence ${\cal F}_2$. Evidently 
$$v_nw_m+w_mv_n=0,\,\, v_nw_m^*+w_m^*v_n=\delta_{n,m} (v,w)I,\,\,\pi(g)v_i\pi(g)^*=(gv)_i.$$
For $i,j=1,2$ we define
$$E_{ij}(k) =\sum_{n-m=k} e_n^{(i)}(e_m^{(j)})^*,\,\, 
E_{ij}(0)=\sum_{n\ge 0} e_{-n}^{(i)}(e_{-n}^{(j)})^* -\sum_{n>0} (e_n^{(j)})*e_n^{(i)}.$$
Identifying $E_{ij}$ with the canonical basis of ${\rm End}\,V=M_2{(\Bbb C})$, we can extend this  
definition by linearity to define operators $X(n)$ for $X\in {\rm End}\, V$ . Note that these operators are even. We
define $D=L_0\otimes I + I \otimes L_0$, so that $[D,v_i]=-iv_i$ and $[D,X(n)]=-nX(n)$. More generally set
${\cal L}_n=L_n\otimes I + I\otimes L_n$. Thus 
$$[{\cal L}_m,{\cal L}_n]=(m-n){\cal L}_{m+n}+\delta_{m+n,0}{m^3-m\over 6} I,\,\,[{\cal L}_m,v_n]=-nv_{n+m}.$$  

\vskip .1in
\noindent \bf Lemma. \it (a) $\pi(g)X(n)\pi(g)^*=(gXg^*)(n)$ for $g\in U(2)$.

\noindent (b) $E_{11}(n)=a_n\otimes I$, $E_{22}(n)=I\otimes a_n$.

\noindent (c) $X(n)^*=X^*(-n)$.

\noindent (d) $[X(n),v_i]=(Xv)_{i+n}$.

\noindent (e) $[{\cal L}_m,X(n)]=-nX(n+m)$.

\noindent (f) $[X(m),Y(n)]= [X,Y](m+n) +m\delta_{m+n,0}\cdot {\rm Tr}(XY)I$.

\vskip .05in
\noindent \bf Proof. \rm The first three identities are immediate. By complex
linearity it suffices to check the remaining relations for $X,Y$ skew--adjoint, i.e. in the Lie algebra of $U(2)$.  
Relation (d) is evident for diagonal operators using (b) and the relation $[a_m,e_n]=e_{n+m}$. Using (a) and the relation
$\pi(g)v_i\pi(g)^*=(gv)_i$, it follows for skew--adjoint $X$, since we can choose $g\in U(2)$ so that $gXg^*$ is diagonal.
Similarly (e) is evident for diagonal $X$ and follows for skew--adjoint $X$ from the fact that the 
operators $\pi(g)$ and ${\cal L}_n$ commute, because $L_n=L^\prime_n$.  
From (d), it follows that $[X(m),Y(n)]-[X,Y](m+n)$ commutes with the $v_i$'s. Taking adjoints, it also commutes with the
$v_i^*$'s. It changes energy by $-n-m$, so vanishes if $n+m\ne 0$. Otherwise, if $n=-m$, it preserves energy and therefore
takes the vacuum vector $\Omega\otimes \Omega$ to a scalar multiple of itself. Since the vacuum vector is cyclic, it acts
as this scalar everywhere, so that
$$[X(m),Y(n)]-[X,Y](m+n) =\delta_{m+n,0}\cdot \lambda_m(X,Y).$$
Evidently $\lambda_m(X,Y)$ is bilinear in $X$ and $Y$ and, conjugating the left hand side by $\pi(g)$, 
satisfies 
$$\lambda_m(gXg^*,gYg^*)=\lambda_m(X,Y).$$
Thus $\lambda_m(X,Y)={\rm Tr}((X\otimes Y) Z)$ where
$Z\in {\rm End}\, V\otimes V$ commutes with $U(2)$, or  equivalently, since $Z$ automatically commutes scalars, 
with $SU(2)$. But as a representation of $SU(2)$, 
$V=V_{1/2}$ and $V_{1/2}\otimes V_{1/2}=V_0\oplus V_1$. Thus by Schur's lemma, the space of 
invariant bilinear forms is two--dimensional, with an easily identified basis: 
$$\lambda_m(X,Y)=\alpha_m \cdot {\rm Tr}(X)\cdot {\rm Tr}(Y)+\beta_m\cdot{\rm Tr}(XY).$$ 
Since $E_{11}(m)=a_m\otimes I$ and $E_{22}(n)=I\otimes a_{n}$ commute and $[a_m\otimes I ,a_n\otimes I]=m\delta_{m+n,0}I$,
it follows that $\alpha_m=0$ and $\beta_m=m$, as claimed.
\vskip .1in

\rm We have thus constructed an action of $\widehat{\u_2}$ on ${\cal F}_2$ by operators $X(n)$. Let $L^{(1)}_n=L_n\otimes I$ and $L^{(2)}_n=I\otimes L_n$. 
From the above, if $H_i=E_{ii}$, then $L^{(i)}_n$ is the Virasoro algebra 
associated to the bosonic system 
$H_i(k)$ and ${\cal L}_n=L^{(1)}_n + L^{(2)}_n$ satisfies $[L_n,X(m)]=-mX(m+n)$. 
Set $H={1\over 2}(H_1-H_2)$, $K={1\over 2}(H_1+H_2)$, 
$E=E_{12}$ and $F=E_{21}$. Then the operators $K(n)$ commutes with the 
operators $E(m)$, $F(m)$ and $H(m)$.  Let $L^{SU(2)}_n$ and $L^{U(1)}_n$
be the Virasoro operators constructed using the two oscillator algebras $(\sqrt{2} H(n))$ and $(\sqrt{2} K(n))$.
Since $H_1^2 + H_2^2= 2(H^2+K^2)$, we have $L_n=L^{SU(2)}_n + L^{U(1)}_n$.
Explicitly we have
$$L_0^{SU(2)}=H(0)^2 +2\sum_{n>0} H(-n)H(n),\,\,\,\,
L_0^{U(1)}=K(0)^2 +2\sum_{n>0} K(-n)K(n).$$  
Since $K(n)$ commutes with $X(m)$ for $X\in {\rm Lie}\, SU(2)$, 
it follows that $[L^{SU(2)}_n,X(m)]=-m X(m+n)$. 
By construction, since the operators $H(m)$ leave any $\widehat{\s}$--submodule invariant, 
so do the Virasoro operators $L^{SU(2)}_n$; moreover $L^{U(1)}_n$ leaves invariant the multiplicity spaces 
of the distinct level $1$ representations of $\widehat{\s}$. 
 
There are at most two such representations, classified by their lowest energy spaces $\H(0)$, 
isomorphic to either $V_0$ or $V_{1/2}$. In fact let $V$ be an 
irreducible $SU(2)$--submodule of ${\cal H}(0)$ and let ${\cal K}$ be the $\widehat{\s}$--submodule 
generated by $V$.  By the commutation relations any monomial in the 
$X(n)$'s can be written as a sum of monomials of the form $RDL$
with $R$ a monomial in energy raising operators $X(-n)$ ($n>0$), $D$ a monomial in constant energy operators $X(0)$, 
and $L$ a monomial in energy lowering operators $X(n)$ ($n>0$). Hence ${\cal K}$ is spanned by products $Rv$. 
But then clearly ${\cal K}(0)=V$. By irreducibility ${\cal H}={\cal K}$ and hence 
${\cal H}(0)=V$. We claim that if $V=V_j$, then $j=0$ or $1/2$. Indeed let $e=F(1)$, $f=E(-1)$ 
and $2h=[e,f]=(E_{22}(0)-E_{11}(0)) + I= -2H +I$. 
Thus $h^*=h$, $e^*=f$, $[h,e]=e$ and $[h,f]=-f$. 
Suppose that ${\cal H}(0)\cong V_j$ and that $v\in V_j$ satisfies $Hv=jv$. 
So $hv=(1-2j)v$ and $ev=0$. By standard $\s$ theory, it follows that $1-2j\ge 0$, i.e. $j=0$ or $j=1/2$. 
Note that if ${\cal H}(0)\cong V_j$ with $j=0$ or $j=1/2$,  
then each ${\cal H}(n)$ decomposes as a sum of integer ($j=0$) or half--integer spin ($j=1/2$).
Since the decomposition of ${\cal F}_2$ contains both integer and half--integer spin representations of $SU(2)$,
it follows that there are irreducible positive energy representations with ${\cal H}(0)=V_j$ ($j=0,1/2$). 

To prove uniqueness, note that any monomial $A$ in operators from $\widehat{\s}$ is a sum of monomials $RDL$ 
with $R$ a monomial in energy raising operators $X(-n)$ ($n>0$), $D$ a monomial in constant energy operators $X(0)$, 
and $L$ a monomial in energy lowering operators $X(n)$ ($n>0$). If $v,w\in {\cal H}(0)$, then the inner products 
$(A_1v,A_2w)$ are uniquely determined by $v,w$ and the monomials $A_i$. Indeed $A_2^*A_1$ is a sum of terms $RDL$
and $(RDLv,w)=(DLv,R^*w)$, with $R^*$ an energy lowering operator. 
Hence if $\H^\prime$ is another irreducible positive energy representation
with ${\cal H}(0)\cong {\cal H}^\prime(0)$ by a unitary isomorphism $v\mapsto v^\prime$, then $U(Av)=Av^\prime$ defines 
a unitary map of ${\cal H}$ onto ${\cal H}^\prime$ intertwining the action of $\widehat{\s}$.  

\vskip .1in
\noindent \bf 3. Character formulas for affine Lie algebras. \rm The character of the spin $j$ 
representation $V_j$ of $SU(2)$, and hence $\s$, is given by
$$\chi_j(\zeta)={\zeta^{2j+1}-\zeta^{-2j-1}\over \zeta -\zeta^{-1}},$$
writing $\zeta$ for the element $\pmatrix{\zeta & 0\cr 0 & \zeta^{-1}\cr}$ of $SU(2)$.
The character of the level $1$ spin $j$ representation $\H_j$ ($j=0,1/2$) is given by 
$$X_j(\zeta,q)={\rm ch}\, \H_j=
\sum_{k=0}^\infty \chi_j(\zeta) (q^{(j+k)^2}-q^{(j+k+1)^2})\varphi(q)=\sum_{n\in j+{\Bbb Z}}\zeta^{2n} q^{n^2} \varphi(q),$$
where 
$$\varphi(q)=\prod_{n\ge 1} (1-q^n)^{-1}.$$
This formula can be deduced directly from first prinicples or by decomposing the 
fermionic representation on ${\cal F}^{\otimes 2}$. In both cases we use the fact that
there are just two level one representations of 
${\cal L}\s$, $\H_j$ with $j=0,1/2$, $\H_0$ 
involving only integer spin representations of $SU(2)$ and $\H_{1/2}$ representations 
only non--integer spins. We use the fact that the operators $L^{SU(2)}_n$ are given by the construction of the Virasoro algebra associated with the bosonic system $(H(k))$. 
\vskip .1in
\noindent \bf Method I. Direct treatment. \rm  We consider first the It is generated by a 
lowest energy vector $\Omega$ with $L_0\Omega=0$. Let $\H[i]=\{\xi\in \H_0:H(0)\xi=i\xi\}$. 
Since $\Omega$ is cyclic, 
these eigenspaces are invariant under the Virasoro algebra and are non--zero only if 
$i$ is an integer. We know that $\H[0]$ is non--zero. We claim that $\H[i]$ is non--zero for every $i\in {\Bbb Z}$.
Indeed $E(-n)$, $F(n)$ and $X={1\over 2}[E(-n),F(n)]$ give a copy of $\s$ with $X\Omega=-n\Omega$. 
Hence $\xi_n=E(-n)^n\Omega\ne 0$ is a non--zero vector in $H[n]$. 
Similarly using $E(n)$, $F(-n)$ and $Y={1\over 2}[E(n), F(-n)]$ as the copy of $\s$, 
it follows that  $\xi_{-n}=F(-n)^n\Omega\ne 0$ is a non--zero vector in $H[-n]$. 
Now $\H[0]$ has a unique lowest energy vector
and is an irreducible oscillator module. For any singular vector would satisfy $L^{SU(2)}_0 \xi=0$ 
and hence be a multiple of the vacuum vector. 

We claim that the $\xi_n$ is a lowest energy vector in $\H[n]$ 
and generates $\H[n]$ as an oscillator module. Clearly $[L_0,E(-n)^n]=n^2E(-n)^n$ and $[L_0,F(-n)^n]=n^2F(-n)^n$ 
so that $L_0\xi_n=n^2$.  But with the present normalisation $L_0=H(0)^2 +2\sum_{n>0} H(-n)H(n)$, so that 
$H(m)\xi_n=0$ for $m>0$. If there were other singular vectors in $\H[n]$ with $H(0)\eta=n\eta$ and 
$H(m)\eta=0$ for $m>0$, then $L_0\eta=n^2\eta$ so that $F(n)^n\eta$ or $E(n)^n\eta$ would 
be a non--zero vector in $\H[0]$. But then it would be be proportional to the vacuum vector $\Omega$ and, 
by $\s$ theory, $\eta$ would be proportional to $\xi_n$. Consequently the character of the vacuum representation
is $\sum_{n\in \Bbb Z} q^{n^2} z^n \varphi(q)$.

We repeat this argument for the spin $1/2$ representation $\H_{1/2}$. In this case the lowest energy representation gives 
non-zero lowest energy vectors in $\xi_{\pm 1/2}=\H[\pm 1/2]$ of energy $1/4$ with respect to the operator $L_0^{SU(2)}$. 
This time the vectors $\xi_{n+1/2}=E(-n-1)^n\xi_{1/2}$ and $\xi_{-n-1/2}=F(-n-1)^n\xi_{-1/2}$ 
are non--zero vectors in $H[\pm(n+1/2)]$. They have energy $(n+1)n + 1/4=(n+1/2)^2$. 
As before an $\s$ 
argument shows that, up to scalar multiples, these are the only singular vectors in $\H[\pm(n+1/2)]$. Hence
the character of the spin $1/2$ representation is $\sum_{n\in 1/2 + {\Bbb Z}} q^{n^2} z^n \varphi(q)$.       
\vskip .1in
\noindent \bf Method II. Treatment using fermions. \rm   
Using the direct sum 
decomposition from the previous section, as a representation of $U(1)$,
${\cal F}$, as a representation of $U(1)$, has character
$$\theta(z,q)=\sum_{n\in {\Bbb Z}} z^n q^{n^2/2}\varphi(q),$$
for $z\in U(1)$. We now use the 
identity:
$$\eqalign{\theta(z,q)\theta(\zeta^{-1},q)&=\sum_{m,n\in {\Bbb Z}} (z\zeta^{-1})^m(z\zeta)^nq^{(m^2 +n^2)/2} \varphi(q)^2\cr
&=\sum_{j=0,1/2}\sum_{a,b\in j+{\Bbb Z}} z^{2a}\zeta^{2b}q^{a^2+b^2}\varphi(q)^2\cr 
&=\sum_{j=0,1/2}X_j(\zeta,q)\Psi_j(z,q),\cr}$$
where 
$$\Psi_j(z,q)=\sum_{m\in j+ {\Bbb Z}} z^{2m} q^{m^2}\varphi(q).$$
We have already exhibited explicit singular vectors 
$\Omega_p\in {\cal F}$  
with $p\in {\Bbb Z}$. By boson--fermion duality these exhaust the singular vectors. 
This yields singular vectors $\Omega_p\otimes \Omega_q$ in ${\cal F}^{\otimes 2}$. If we split up the 
character of ${\cal F}^{\otimes 2}$ as the sum of characters of 
integer and non--integer spin $\Theta_j(z,\zeta,q)$, 
then $\Theta_j(z,\zeta,q)=X_j^\prime(\zeta,q)\cdot \Psi_j^\prime(z,q)$,
with $X^\prime_j(\zeta,q)$ the character of $\H_j$ and $\Psi_j^\prime(z,q)$ the character
of the mutliplicity space of $\H_j$.
If we fix a spin $j=0,1/2$, then there are certainly singular vectors of charge $k+j$ for each $k\in {\Bbb Z}$ for the $H(n)$'s.
Thus 
$$X_j^\prime(\zeta,q) \ge X_j(\zeta,q)=\sum_{n\in j+{\Bbb Z}} \zeta^{2n} q^{n^2}\varphi(q),\eqno{(1)}$$
where the inequality is to be taken coefficient by coefficient.
Similarly, for the character of the multiplicity space,
$$\Psi_j^\prime(z,q) \ge \Psi(z,q)=\sum_{n\in j+{\Bbb Z}} z^{2n} q^{n^2}\varphi(q).\eqno{(2)}$$
On the other hand the identity above implies that
$$\sum_j X_j^\prime(\zeta,q)\Psi_j^\prime(z,q) =\sum X_j(\zeta,q)\Psi_j(z,q),$$
so that equality must hold in (1) and (2). 
\vskip .1in
\noindent \bf Remark. \rm In Section~8 we will introduce shift operators which together with the $H(n)$'s and the $K(n)$'s act irreducibly on the $\H_j$'s and their multiplicity spaces. This gives a more conceptual operator--theoretic explanation of the second proof. 
\vskip .1in
\noindent \bf 4. Existence of singular vectors in the oscillator representations. \rm According to the 
character formula $\H_j$ has a series of vectors $\xi_m$ ($m\in j+ {\Bbb Z}$) with $H(0)\xi_m=m\xi_m$ 
and $H(n)\xi_m=0$ for $m>0$. Moreover $\H_j$ is the direct sum of the cyclic 
oscillator modules $\K_{m}$ generated by these vectors;  note that, if $H=H(0)$, then $\K_m$ is the $m$--eigenspace of $H$. 
Following Graeme Segal, we construct non-zero singular vectors in $\K_j\subset \H_j$. 
In fact the vector $\xi_m$ satisfies $H\xi_m=m\xi_m$. If $m\ge 0$ then $\eta=E\xi_m$ lies in $ \K_{m+1}$ and satisfies 
$L_0\eta={1\over 2} m^2\eta$. Since the lowest energy level in $\K_{m+1}$ is ${1\over 2}(m+1)^2$ it follows that 
$\eta=0$, i.e. $E\xi_m=0$. Similarly if $m\le 0$, $F\xi_m=0$. Thus $\xi_m$ generates an irreducible $\s$--representation
of dimension $2|m|+1$. If $m>0$, then $F^{[m]}\xi_m$ is a non--zero singular vector in $\K_j$ of energy $m^2/2$. 
If $m<0$ then $E^{[m]+1} \xi_m$ is a non--zero singular vector in $\K_j$ of energy $m^2/2$. This proves:       
\vskip .1in
\noindent \bf Proposition. \it The oscillator representation with $H(0)=jI$ has non--singular singular 
vectors for the Virasoro algebra with energy $(j+n)^2/2$ for $n\ge 0$. The non--zero singular vector constructed above 
in $\K_j$ of energy $m^2/2$ ($m\ge 0$) generates a copy of $V_m$ as an $\s$-module.
\vskip .1in 
\noindent \rm We shall call these singular vectors {\it Goldstone vectors}. In Section 8, following Segal, 
we will use vertex operators to deduce the explicit formulas of Goldstone for these vectors.
\vskip .1in
\noindent \bf 5. Uniqueness of singular vectors in the oscillator representations. \rm Let $(a_m)$ 
satisfy $[a_m,a_n]=m\delta_{m+n,0}I$. These act as operators on $V={\Bbb C}[x_1,x_2,\dots]$ via $a_0=\sqrt{2}\mu$, $a_n=n\partial_{x_n}$ 
and $a_{-n}=x_n$ for $n>0$. Set $L_0=a_0^2/2 + \sum_{n>0} a_{-n}a_n$ and $L_k={1\over 2} \sum_{p+q=k} a_pa_q$ for $k\ne 0$.
Thus $[L_{k},a_n]=-na_{n+k}$ and $[L_m,L_n] =(m-n)L_{m+n} + \delta_{m+n,0} (m^3-m)/12$. 
\vskip .1in
\noindent \bf Proposition. \it For $n\ge 0$ let $V(n)=\{\xi: L_0\xi=(\mu^2+n)\xi\}$ with $\mu\in {\Bbb R}$. Then 
${\rm dim}\, V(n)=\P(n)$ the partition function and the subspace of singular vectors, i.e. sRolutions of 
$L_1\xi=0=L_2\xi$, is at most one-dimensional.

\vskip .05in
\noindent \bf Proof. \rm If we define the degree of $x_i$ to be $i$, then $V(\mu^2+n)$ is spanned by all monomials 
of total degree $n$, so that ${\rm dim}\, V(n)=\P(n)$. We shall treat a more general case which evidently implies the proposition:
\vskip .1in
\noindent \bf Lemma. \it Let $W={\Bbb C}[x_1,x_2,\dots]$ with ${\rm deg}\, x_i=i$ and let $W(n)$ be the 
space spanned by monomials of total degree $n$, so thet ${\rm dim}\, W(n)=\P(n)$. For $a,b\in {\Bbb C}$ let
$$A=a \partial_{x_1}+\sum_{n> 0} (n+1)x_{n}\partial_{x_{n+1}},\,\, 
B= {1\over 2}\partial_{x_1}^2 + b\partial_{x_2}+ \sum_{n\ge 1}(n+2)x_n\partial_{x_{n+2}}.$$ 
Then the space of solutions of $W_0(n)=\{p\in W(n):Ap=0,\,\, Bp=0\}$ is at most one--dimensional.
\vskip .05in
\noindent \bf Proof. \rm For $i\ge 0$, let $W(n,i)=\{x_1^{n-i}q(x_2,\dots): {\rm deg}\,q=i\}$. Evidently 
$$W(n)=\bigoplus_{i\ge 0} W(n,i)$$
so that any element $p\in W(n)$ has a decomposition $p=\sum p_i$ with $p_i\in W(n,i)$. Thus 
$p_i=x_1^{n-i}q_i(x_2,x_3,\dots)$ with ${\rm deg}\, q_i=i$. We assume that the leading coefficient 
$p_0$ (a multiple of $x_1^n$) is zero and prove by induction on $i$ that $q_i=0$. Each $q_i $ 
is a sum of linear combination of monomials $x^\alpha=x_2^{\alpha_2} \dots x_n^{\alpha_n}$ with 
$\sum_{k\ge 2} k \alpha_k= i$. Suppose that $\alpha_2>0$. Then 
$$Ax_1^{n-i}x^\alpha= 
a(n-i)x_1^{n-i-1}x^\alpha + 2x_1^{n-i+1}\alpha_2 x_2^{\alpha_2-1} 
x_3^{\alpha_3}\cdots  + 3\alpha_3 x_1^{n-i} x_2^{\alpha_2+1} x_3^{\alpha_3-1}x_4^{\alpha_4} \cdots + \cdots .$$
Because of the inductive hypothesis, this is the only monomial which 
produces a term $x_1^{n-i+1}x_2^{\alpha_2-1}x_3^{\alpha_3}\cdots$ 
when $A$ is applied. It follows that its coefficient must be zero. So only the coefficients of monomials $x_1^{n-i}
 x_3^{\alpha_3}x_4^{\alpha_4}\cdots$ The image under $B$ of such a monomial is
$$B x_1^{n-i} x_3^{\alpha_3}x_4^{\alpha_4}\cdots = 
{1\over 2} (n-i)(n-i-1)x_1^{n-i-2}  x_3^{\alpha_3}x_4^{\alpha_4}\cdots + 3\alpha_3x_1^{n-i+1} x_3^{\alpha_3-1} x_4^{\alpha_4}
\cdots.$$
Because of the inductive hypothesis this is the only monomial which produces a 
term $x_1^{n-i-2}  x_3^{\alpha_3}x_4^{\alpha_4}\cdots $ when $B$ is applied. 
It follows that its coefficient must be zero. Thus $p_i=0$ and hence $p=0$, as required. 
\vskip .1in
\noindent \bf 6. Density modules, primary fields and the Feigin--Fuchs product formula. \rm For 
$\lambda,\mu\in {\Bbb C}$, we define the density module $V_{\lambda,\mu}$ 
to be the vector space with basis $v_n$ ($n\in {\Bbb Z}$) and define
operators $\ell_k$ on $V_{\lambda,\mu}$ by
$$\ell_k v_n=-(n+ \lambda k +\mu)v_{n+k}.$$
It is easy to check that $[\ell_m,\ell_n]=(m-n)\ell_{n+m}$, so this defines a representation of the Virasoro algebra
with $c=0$. Clearly the change of basis by the shift operator $v_n\mapsto v_{n+k}$ 
gives a natural isomorphim
between $V_{\lambda,\mu}$ and $V_{\lambda,\mu+k}$.  

Recall that if $\pi:\g\rightarrow {\rm End}\, W$ is a representation of a Lie algebra $\g$, then the {\it dual
representation} on the dual space $W^\prime$ is given by 
$(\pi^\prime(X)\xi,v)=-(\xi,Xv)$ for $v\in W$. Let $w_n$ be the canonical basis of $V_{1-\lambda,-\mu}$. We identify
$U$ with a subspace of the dual $V_{\lambda,\mu}$ via the pairing $(v_n,w_m)=\delta_{n+m,0}$. It is immediately 
verified that the dual representation restricts to the natural representation on $V_{1-\lambda,-\mu}$. With 
an obvious abuse of notation, we shall write $V_{1-\lambda,-\mu}\cong V_{\lambda,\mu}^\prime$. 
These representations are the infinitesimal version of the action of ${\rm Diff}\, S^1$, 
or more properly a certain central extension, on the densities $e^{i\mu\theta}f(\theta)(d\theta)^\lambda$. The 
duality between $V_{\lambda,\mu}$ and $V_{1-\lambda,-\mu}$ is given by the pairing into 1--forms
$$(f(\theta)e^{i\mu\theta}(d\theta)^\lambda,g(\theta)e^{-i\mu\theta}(d\theta)^{1-\lambda}) \mapsto fg d\theta$$
followed by integration over the circle, i.e. ${1\over 2\pi} \int_0^{2\pi} fg\,d\theta$. Identifying $v_n$ with
$e^{in\theta}e^{i\mu\theta}(d\theta)^\lambda$, we see that
$$\ell_k(fe^{i\mu\theta}(d\theta)^\lambda))=-e^{ik\theta}(i{d\over d\theta} + k\lambda)(f e^{i\mu\theta}) (d\theta)^\lambda.$$
Thus it is the action on ${\Bbb C}[e^{i\theta},e^{-i\theta}]e^{i\mu\theta}$ given by the operators
$$\ell_k=-e^{ik\theta}(i{d\over d\theta} + k\lambda).$$
Changing variables to $z=e^{i\theta}$, regarded as a formal variable, and introducing an extra factor $z^\lambda$, it may also be identified with the action
on ${\Bbb C}[z,z^{-1}]z^{\mu+\lambda}$ given by
$$\ell_k= -z^{k+1}{d\over dz} - (k+1)\lambda z^k.$$ 
For applications here $\lambda$ and $\mu$ will be in ${\Bbb Z}$ (for applications to fusion 
we require $\lambda,\mu\in {1\over 4}{\Bbb Z}$).

For $i=1,2$, let $\H_i=L(1,h_i)$ with lowest energy vectors $\xi_i=\xi_{h_i}$. 
We define a primary field of type $(\lambda,\mu)$ between 
$\H_1$ and $\H_2$ to be a map 
$\Phi:\H_1\otimes V_{\lambda,\mu}\rightarrow \H_2$ commuting with the
action of the Virasoro algebra. Let $\Phi(n)\xi=\phi(\xi\otimes v_n)$. Then
$$[L_k,\Phi(n)]=-(n+k\lambda +\mu)\Phi(n+k).$$
We have 
the following uniqueness result:
\vskip .1in 
\noindent \bf Lemma~A. \it If $(\Phi(n)\xi_1,\xi_2)=0$ then $\Phi(n)=0$ for all $n\in {\Bbb Z}$.  If $\Phi$ is non--zero, then $m=h_1-h_2-\mu$ is an integer and 
$\Phi(m)\xi_1$ is a non--zero multiple of $\xi_2$. 
\vskip.05in
\noindent \bf Proof. \rm Suppose that $(\Phi(n)\xi_1,\xi_2)=0$. From the commutation relations it follows that for $n_i>0$ 
$$(\Phi(n)L_{-{n_k}} \cdots L_{-n_1}\xi_1,\xi_2)=0.$$
Hence $(\Phi(n)\xi,\xi_2)=0$ for all $\xi$. From the commutation relations is then follows that for $n_i>0$ 
$$(\Phi(n)\xi,L_{-{n_k}} \cdots L_{-n_1}\xi_2)=0,$$
so that $(\Phi(n)\xi,\eta)=0$ for all $\xi$, $\eta$ and hence $\Phi(n)=0$ for all $n\in {\Bbb Z}$. Finally if $(\Phi(m)\xi_1,\xi_2)\ne 0$, then $\Phi(m)\xi_1$ is a non--zero multiple of $\xi_2$. Hence $h_2=h_1-m-\mu$.  
\vskip .1in 

\noindent \bf Remark. \rm This result allows us to normalise a primary field by shifting 
numbering of the modes $\Phi(n)$ to $\Phi(n+k)$. Thus, if a primary field exist, we
may assume that $\mu=h_1-h_2$ and that $\Phi(0)\xi_1=\alpha \xi_2$ with $\alpha\ne 0$. 
\vskip .1in
In the language of vertex operators and vertex algebras, primary fields are usually described in terms of generating 
functions in a formal variable $z$. We set
$$\Phi(z)=\sum_{n\in {\Bbb Z}} \Phi(n)z^{-n-\delta}.$$
Then 
$$[L_k,\Phi(z)]=z^{k+1} {d\Phi(z)\over dz} +(k+1)\Delta z^k\Phi(z).$$
Taking coefficients of $z^{-n-\delta}$, we get
$$[L_k,\Phi(n)] =-(n+k+\delta  -(k+1)\Delta)\Phi(n+k).$$
Thus $\lambda=-\Delta+1$ and $\mu=\delta-\Delta$. Hence $\Delta=1-\lambda$ and $\delta=\mu +\lambda-1$. 

Now the Verma module representation $M(1,j^2)$ ($j=0,1/2,1,3/2.\dots$) certainly has one singular vector at energy
$(j+1)^2=j^2+d$ with $d=2j+1$, since the representation $L(1,j^2)$ is a subrepresentation of the multiplicity space. 
This corresponds to a element of the universal enveloping algebra ${\cal U}_1$ generated by $L_{-k}$ ($k\ge 1$). Let
${\cal U}_k$ be the universal enveloping algebra generated by $L_{-i}$ with $i\ge k$. It is the universal 
enveloping algebra of the Lie algebra generated by $L_{i}$ for $i\ge k$. This is a Lie ideal in the Lie algebra generated 
by $L_{-j}$ with $j\ge 0$. As a consequence every element of the algebra ${\cal U}_1$ can be written uniquely in the form
$P=\sum_{k\ge 0} p_k L_{-1}^k$ and ${\rm ad}L_{-1}$ defines a derivation of ${\cal U}_2$. In particular the
singular vector  has the form $Pv_{j^2}$ where 
$$P\equiv P_d=\sum_{k =0}^d q_k L_{-1}^k.$$
Each term is a $q_k$ is a sum of monomials $L_{-i_1} \cdots L_{-i_r}$ where $i_1+\cdots +i_r=d-k$. In particular 
$q_d$ is a scalar.  The following uniqueness result was stated by Fuchs.
\vskip .1in 

\noindent \bf Lemma~B. \it The coefficient of $L_{-1}^d$ in $P_d$ is non--zero; in particular if a singular 
vector exists, it is unique up to a scalar multiple.
\vskip .05in

\noindent \bf Proof. \rm 
Suppose that $q_d=0$. We prove that $P_d=0$. 
Each $q_d$ can be written as a linear combination of monomials $L_{-s}^{n_s} \cdots L_{-2}^{n_2}$ with $n_i\ge 0$. 
Amongst all monomials with non--zero coefficients, we can find one $p\ge 2$ minimal with $n_p>0$. 
We can also choose this monomial so that $n=n_p$ is also minimal. Suppose that it occurs in $q_c$ with $0\le c<d$. 
We may in addition assume that $c$ is chosen to be maximal. Let $w=P_dv$ with $v$ the cyclic lowest energy 
vector in Verma module. Then $L_{p-1}w=0$. This can be written as
$$\sum_{k=0}^c L_{p-1}q_k L_{-1}^{k}v=\sum_{k=0}^c[L_{p-1},q_k]L_{-1}^kv+q_kL_{p-1}L_{-1}^kv.$$
We look for terms ending with $L_{-p}^{n-1}L_{-1}^{c+1}$ in this expression. We first note that it 
follows by induction on $k$ that if $j\ge 1$ then
$L_jL_{-1}^kv$ lies in ${\rm lin}\, \{AL_{-1}^iv:i<k,\,\, A\in {\cal U}_2\}$. Indeed
$$L_{j}L_{-1} L_{-1}^{k-1} v=(j+1)L_{j-1}L_{-1}^{k-1}v +L_{-1}L_jL_{-1}^{k-1}v,$$
which has the same form since $L_{-1}{\cal U}_2={\cal U}_2 L_1$ and $L_{-1}^{k-1}v$ is an eigenvector of $L_0$.
If $k\ne c$, then the minimality of $p$ forces $q_k$ to be in ${\cal U}_p$. So either $k> c$, in which case any monomial
$q_k$ could have a non--zero contribution from a monomial $L_{-s}^{m_s} \cdots L_{-p}^{m_p}$ with $m_p>n_p$. 
Clearly taking the Lie bracket with $L_{p-1}$ can diminish the exponent 
of $L_{-p}$ by at most one and at the same time must 
increase the exponent of $L_{-1}$ by at least one. 
So no monomials ending with $L_{-p}^{n-1} L_{-1}^{c+1}$ can appear with non--zero coefficient.
If $k<c$, there is no way to increase the power of $L_{-1}^k$ to $L_{-1}^{c+1}$ 
by taking the Lie bracket with $L_{p-1}$. For $k=c$ and the terms $[L_{p-1},AL_{-p}^{n}] $ with $A$ a 
monomial in ${\cal U}_{p+1}$, we have
$$[L_{p-1},AL_{-p}^{n}]=[L_{p-1},A]L_{-p}^{n}+A[L_{p-1},L_{-p}^{n_p}].$$
The first term lies in ${\cal U}_2$ while for the second
$$[L_{p-1},L_{-p}^{n}]v=(2p-1) \sum_{a+b=n-1} L_{-p}^a L_{-1} L_{-p}^bv=(2p-1)n L_{-p}^{n-1}L_{-1} + B,$$
where $B\in {\cal U}_2$. There could be several terms in $q_c$ with $n_p=n$, but on bracketing with $L_{p-1}$ 
and taking the term ending with $L_{-p}^{n-1}L_{-1}^{c+1}$ is equal to $A  L_{-p}^{n-1}L_{-1}^{c+1}$. Thus 
all these terms are distinct. But then there can be no  cancellation and 
the coefficient of $ L_{-r}^{n_r} \cdots L_{-p}^{n_p-1}L_{-1}^{c+1} v$ must therefore be non--zero, a contradiction. 
\vskip .1in
\noindent \bf Remark. \rm Fuchs used a slightly different ordering on monomials, which can
also be used to give prove the uniqueness theorem. For a monomial $L_{-s}^{n_s} \cdots L_{-2}^{n_2}L_k$, we set $n_1=d-k$ and use the lexicographic ordering determined by 
$(n_1,n_2,\dots)$. The proof is almost word--for-- word that used in the second proof of
Fuchs' algorithm in Section~10. 
\vskip .1in
Let $P_d\in {\cal U}_1$ be the element giving the singular vector normalized so that the coefficient of 
$L_{-1}^d$ equals~$1$. In the representation $V_{\lambda,\mu}$, we evidently have 
$$P_d v_0 = a_d(\lambda,\mu)v_{-d}$$
where $a_d(\lambda,\mu)$ is a inhomogenous polynomial of degree $d$ in $\lambda$ and $\mu$. Indeed since the 
$L_{-1}^d$ appears with coefficient $1$, the formula for the action shows that the coefficient of $\mu^d$ is $(-1)^d$.
Using the coset construction in Section~8 and an explicit formula for $P_d$ in Appendix~C, 
we give two methods
to establish the following product formula of Feigin--Fuchs, which we state first in the two simplest cases:
\vskip .1in
\noindent \bf Feigin--Fuchs product formula. \it Let $j=0,1/2,1,3/2,\dots$ and set $S=\{-j, -j+1, -j+2, \dots, j\}$. 

\noindent (a) $a_d(0,\mu) =(-1)^d \prod_{k\in S}(\mu+j^2 -k^2)$.

\noindent (b) $a_d(1,\mu)=(-1)^d \prod_{k\in S} (\mu+j^2 -(k+1)^2)$

\noindent (c) $a_d(p^2,\mu)=(-1)^d \prod_{k\in S} (\mu +j^2 -(k+p)^2)$ for $p\ge 0$.

\noindent (d) $a_d(\lambda,\mu)^2=\prod_{k\in S} ((\lambda -\mu -j^2+k^2)^2 -4\lambda k^2)$.
\vskip .1in
\rm 
We have the following necessary condition for the existence of a primary field:
\vskip .1in
\noindent \bf Lemma~C. \it If there is a non--zero primary field from $L(1,j^2)$ to $L(1,h)$ of type $(\lambda,\mu)$,
then 
$$a_d(1-\lambda,h-j^2)=0.$$
\vskip .05in
\noindent \bf Proof. \rm 
Let $\Phi$ be a non--zero primary field of type $(\lambda,\mu)$ 
from $L(1,j^2)$ to $L(1,h)$. Then we may assume that $\Phi$ is normalised so that
$\mu=j^2-h$ and 
$\Phi(0)\xi_{j^2}=\alpha \xi_h$ for some constant $\alpha\ne 0$. 
Let $\xi=\xi_{j^2}$ and $\eta=\xi_h$. 
Then $(\Phi(0)\xi,\eta)\ne 0$. On the other hand $P_d\xi=0$. Let $P_d=\sum a_m L_{-1}^{m_1} L_{-2}^{m_2} \cdots$.
For an operator we set $\pi(L_j)A=[L_j,A]$. This defines an action of the Virasoro algebra with $c=0$. Now
$$0=(\Phi(d)P_d\xi,\eta)= 
\sum a_m ( \pi(-L_{-n})^{m_n}\cdots \pi(-L_{-2})^{m_2}\cdot \pi(-L_{-1})^{m_1})\cdot\Phi(d)\xi,\eta).$$
On the other hand
$$\sum a_m \pi(-L_{-n})^{m_n}\cdots \pi(-L_{-2})^{m_2}\cdot \pi(-L_{-1})^{m_1})\cdot\Phi(d)=b\Phi(0),$$
where for the basis $(v_i)$ of $V_{\lambda,\mu}$ 
$$\sum a_m (-L_{-n})^{m_n}\cdots (-L_{-2})^{m_2}\cdot (-L_{-1})^{m_1})v_{d}=b v_0.$$
Thus we require that $b$ should vanish. To calculate $b$, we take $(w_j)$ to be a 
dual basis of $V_{1-\lambda,-\mu}$. Then
$$b=\sum a_m 
((-L_{-n})^{m_n}\cdots (-L_{-2})^{m_2}(-L_{-1})^{m_1}v_{d},w_{0})
=(v_d,P_dw_0).$$
Thus $P_d w_{0}=b w_{-d}$ in $V_{1-\lambda,-\mu}$
and hence
$$b=a_d(1-\lambda,-\mu)=a_d(1-\lambda, h-j^2).$$
\vskip .1in
\noindent \bf Remark. \rm This result can be seen a little more transparently 
if the dual of $V_{\lambda,\mu}$ is used.
Indeed, passing to a suitable completion of the tensor product, a primary field 
of type $(\lambda,\mu)$ from $L(c,h_1)$ to $L(c,h_2)$ can be regarded 
as a map $\H_1\rightarrow \H_2 \widehat{\otimes} V_{\lambda^\prime,\mu^\prime}$,
where $\lambda^\prime=1-\lambda$ and $\mu^\prime=-\mu$. 
Explicitly the map becomes 
$$\xi \mapsto \Phi(z)\xi=\sum \varphi(n)\xi z^{-n},$$ 
where $z^k$ here is identified with a basis element of $V_{\lambda^\prime,\mu^\prime}$. If $v_{h_2}=\phi(-k)v_{h_1}$, 
then, picking out the coefficient of $v_{h_2}$,
the condition $P(L_{-1},L_{-2}, \dots)v_{h_1}=0$ evidently implies that 
$P(\ell_{-1},\ell_{-2},\dots)z^k=0$, i.e. the condition stated in Lemma~C. It is therefore 
natural to refer to the primary field as having {\it dual type} $(\lambda^\prime,\mu^\prime)$. Thus if the primary field is normalised so that $\Phi(0)\xi_1=\alpha\xi_2$ with $\alpha\ne 0$, then $\mu^\prime=h-j^2$ and 
Lemma~C requires that $(a_d(\lambda^\prime,\mu^\prime)=0$.
 
\vskip .1in
\noindent \bf 7. Multiplicity one theorem. \rm We wish to prove the following theorem:
\vskip .1in
\noindent \bf Theorem~A. \it The action of the Virasoro algebra 
on the multiplicity space $M_m$ of 
$V_m$ ($m\ge 0$, $m\in j+{\Bbb Z}$) in $\H_j$ is irreducible for $j=0,1/2$.  
\vskip .1in
Let us start by showing that it suffices to prove:
\vskip .1in
\noindent \bf Theorem~B. \it The Goldstone vectors are the only singular vectors in
$\K_m=\{\xi\in \H_j:H\xi=m\xi\}$ for $m\in j+{\Bbb Z}$.
\vskip .05in
\noindent \bf Proof that Theorem~B implies Theorem~A. \rm  
The multiplicity space $M_m$ can be identified with 
$$M_m=\{\xi\in \H_j:H\xi=m\xi,\,E\xi=0\}\subset \K_m=\{\xi:H\xi=m\xi\}.$$
If the action on $M_m$ is not irreducible, then
$M_m$ would have to contain a singular vector of higher energy. But by 
Theorem~B any such vector would have to be a Goldstone vector in $\K_m$. 
But, by the proposition in Section~4, 
$\xi_m$ is the only Goldstone vector in $\K_m$ annihilated by $E$.

\vskip .1in
Thus we have only to prove Theorem~B. By the uniqueness theorem 
it suffices to show that there are non singular 
vectors in $\K_j$ with energy not of the form 
$m^2/2$ with $m\in j+{\Bbb Z}$. For this we require the 
notion of a particularly simple kind of primary field. Let 
$V_{0,0}={\Bbb C}[z,z^{-1}]$ with $L_n$ acting as the operator $\ell_n=-z^{n+1} d/dz$. 
A$(0,0)$--primary field is an operator 
$$\Phi:L(1,h_1)\otimes V_{0,0}\rightarrow L(1,h_2)$$
which intertwines the action of the Virasoro algebra. 

\vskip .1in
\noindent \bf Lemma~A. \it Let $\Phi(n)\xi=\Phi(\xi\otimes z^{n})$ 
for $\xi\in L(1,h_1)$. Then $[L_m,\Phi(n)]= -n\Phi(n+m)$.
Moreover any such assignment corresponds to a $(0,0)$--primary field via the above formula.
\vskip .05in
\noindent \bf Proof. \rm Equivariance implies that
$$L_m\Phi(n)\xi=L_m(\Phi(\xi\otimes z^n))=\Phi(L_m \xi\otimes z^n) + \Phi(\xi\otimes \ell_mz^n)
=\Phi(n)L_m\xi - m\Phi(x,m+n).$$
\vskip .1in
\noindent \bf Lemma~B. \it If $(\Phi(n))$ is a $(0,0)$--primary field so is $\Psi(n)=\Phi(-n)^*$.
\vskip .05in
\noindent \bf Proof. \rm We have
$$[L_m,\Psi(n)]=[L_m,\Phi(-n)^*]=-[L_{-m},\Phi(-n)]^*=-m\Phi(-n-m)^*=-m\Psi(n+m).$$
\vskip .1in
\noindent \bf Lemma~C. \it If $(\Phi(n)\xi_1,\xi_2)=0$ with $\xi_i$ the lowest energy vector in $L(1,h_i)$, then $\Phi(n)=0$ for all $n\in {\Bbb Z}$.  
\vskip.05in
\noindent \bf Proof. \rm Suppose that $(\Phi(n)\xi_1,\xi_2)=0$. From the commutation relations it follows that for $n_i>0$ 
$$(\Phi(n)L_{-{n_k}} \cdots L_{-n_1}\xi_1,\xi_2)=0.$$
Hence $(\Phi(n)\xi,\xi_2)=0$ for all $\xi$. From the commutation relations is then follows that for $n_i>0$ 
$$(\Phi(n)\xi,L_{-{n_k}} \cdots L_{-n_1}\xi_2)=0,$$
so that $(\Phi(n)\xi,\eta)=0$ for all $\xi$, $\eta$ and hence $\Phi(n)=0$ for all $n\in {\Bbb Z}$.
\vskip .1in 
\bf \noindent Remark. \rm Note that if $L_0\xi_i=h_i \xi_i$, then such a $\Phi$ can exist only if $h_1-h_2$ is an integer. 

\vskip .1in

\noindent \bf Proposition. \it If $j$ is a non--negative half--integer and $k\ne (j+m)^2$ with $m\in {\Bbb Z}$  
there are no non--zero $(0,0)$--primary fields 
from $L(1,j^2)$ to $L(1,k)$ or from $L(1,k)$ to $L(1,j^2)$.
\vskip .05in

\noindent \bf Proof. \rm Using adjoints it suffices to prove the first assertion. If a non--zero primary field exists, 
then by Lemma~C of Section~6, 
necessarily $a_d(1,j^2-k)$ would have to vanish. 
On the other hand by the Feigin--Fuchs product formula,
if $S=\{-j,-j+1,\cdots,j-1,j\}$, we know that
$$a_d(1,j^2-k)=\prod_{t\in S} ((t-1)^2-k)\ne 0.$$
\vskip .1in
\noindent \bf Proof of Theorem~B. \rm We have to prove that the Goldstone vectors are up to scalar multiples the only singular vectors in $\K_m$. Let $P$ be the projection onto the submodule generated by the Goldstone 
vectors in $\K_m$, and $P_0$ 
the projection onto the submodule generated by a Goldstone vector $\xi$.
If there is another singular vector $\eta$ not proportional to a Goldstone vector then its energy does not have the 
form $(j+k)^2$ with $k$ an integer, by the uniqueness theorem. 
Let $Q$ be be the projection on submodule generated by $\eta$.
Then $\Phi(n)=QH(n)P_0$ an $P_0\Phi(n)Q$ are $(0,0)$--primary fields since $L_m$ commutes 
with $P_0$ and $Q$ and $[L_m,H(n)]=-nH(n+m)$. By the proposition, we must have 
$QH(n)P_0=0$ and $P_0H(n)Q=0$. Hence
it follows that $QH(n)P=0=PH(n)Q$. But then $H(n)$ leaves the subspace $P\K_m$ invariant. Since the $H(n)$'s 
act irreducibly on $\K_m$ it follows that $P=I$, a contradiction. So there are no singular 
vectors other than the Goldstone vectors.  

\vskip .1in

\noindent \bf Corollary~1. \it The character of $L(1,m^2)$ for $m$ a non--negative half--integer is 
$(q^{m^2}-q^{(m+1)^2})\varphi(q)$.
\vskip .1in
\noindent \bf Corollary~2. \it Up to a scalar multiple $P_m\cdot \xi_m$ is the unique singular vector of energy $(m+1)^2$  in $M(1,m^2)$ and the quotient module is irreducible. The other singular vectors in $M(1,m^2)$
are given by \break $P_{m+k-1} \cdots P_{m+1} P_m\xi_m$ (of energy $(m+k)^2$).
\vskip .05in
\noindent \bf Proof. \rm By the character formula, there is only one singular vector $\eta$  of energy 
$(m+1)^2$ in $M=M(1,m^2)$  and no singular vectors of lower energy except the generating vector of energy $m^2$. 
The submodule $N$ generated by is isomorphic to $M(1,(m+1)^2)$. 
By the character formula $M(1,m^2)/N$ is irreducible. Any singular vector in $M$ not in 
$N$ its image would give a non--zero singular vector in $L(1,m^2)=M/N$ 
and would have to be a multiple of the generating vector. So it follows that all singular vectors of energy 
greater than or equal to $(m+1)^2$ lie in $N$, whence the result.
\vskip .1in
\noindent \bf 8. Proof of the Feigin--Fuchs product formula using primary fields. \rm We shall prove below 
that if $j$, $j_1$, $j_2$ are integers then there is a primary field of dual type $(j_1^2,j_2^2-j^2)$ 
from $L(1,j^2)$ to $L(1,j_2^2)$ provided that $V_{j_2}$ occurs as a component in $V_j\otimes V_{j_2}$. 
Thus  $a_d(j_1^2,j_2^2-j^2)=0$, by Lemma~C in Section~6. This is enough information to compute $a_d(\lambda,\mu)$ completely. 

We fix $j$, with $d=2j+1$, and freeze $j_2$, and hence $\mu=j_2^2-j^2$, with $j_2\ge j$. 
Then the possible values of $j_1$ are $j_2+k$ with $k\in
S=\{-j,-j+1,\cdots,j-1,j\}$. But $a_d(\lambda,j_2^2-j^2)$ has degree at most
$d$ in $\lambda$ and $d=|S|$. Hence
$$a_d=a_d(\lambda,j_2^2-j^2)=C \prod_{k\in S} (\lambda - (j_2+k)^2),$$
where the constant $C$ might depend on $\mu$. But then, 
using the fact that $-S=S$, we can write
$$\eqalign{a_d(\lambda,j_2^2-j^2)^2&=C^2 \prod_{k\in S} (\lambda -j_2^2 -k^2 -2j_2k)(\lambda -j^2 -k^2 +2j_2k)\cr
&=C^2\prod_{k\in S} ((\lambda-j_2^2-k^2)^2 -4j_2^2k^2)\cr
&=C^2\prod_{k\in S} ((\lambda-\mu-j^2-k^2)^2 -4(j^2+\mu)k^2).\cr}$$
Since $\mu$ can take infinitely many values, it follows that $C^2=1$ and
$$a_d(\lambda,\mu)^2=\prod_{k\in S} ((\lambda-\mu-j^2-k^2)^2 -4(j^2+\mu)k^2).$$
On the other hand
$$\eqalign{(\lambda-\mu-j^2-k^2)^2 -4(j^2+\mu)k^2
&=(\lambda+\mu)^2 -2(\lambda-\mu)(j^2+k^2) +(j^2+k^2)^2 -4j^2k^2-4\mu k^2  \cr
&=\lambda^2 -2\lambda \mu +\mu^2 -2\lambda(j^2+k^2) +2 \mu(j^2-k^2) +(j^2-k^2)^2\cr
&=(-\mu -j^2 +k^2)^2 + 2\lambda(-\mu -j^2+k^2) + \lambda^2 -4k^2\lambda\cr
&= (\lambda-\mu -j^2 +k^2)^2 -4k^2\lambda.\cr}$$
Thus we have
$$a_d(\lambda,\mu)^2=\prod_{k\in S}[ (\lambda-\mu -j^2 +k^2)^2 -4k^2\lambda].$$
When $k=p^2$, this can be rewritten
$$a_d(p^2,\mu)^2=\prod_{k\in S} (-\mu -j^2 +(k+p)^2)^2,$$
again using the fact that $S=-S$. Since the coefficient of $\mu^d$ in $a_d(\lambda,\mu)$ is $(-1)^d$, 
we deduce that
$$a_d(p^2,\mu)=(-1)^d \prod_{k\in S} (\mu+j^2-(k+p)^2),$$
as required.
\vskip .1in
We now explicitly construct all these primary fields. First we recall a version of the 
vertex operator formulation of the boson--fermion correspondence. We define the shift operator $U$ on fermionic Fock
space ${\cal F}$ by
$$U(e_{i_1}\wedge e_{i_2} \wedge e_{i3}\wedge \cdots)=e_{i_1+1}\wedge e_{i_2+1} \wedge e_{i_3+1} \wedge \cdots$$
By definition $U$ is unitary and $Ue_iU^*=e_{i+1}$. Thus uniquely specifies $U$ up to a scalar multiple,
because if $U^\prime$ is another such unitary, then $V=U^\prime U^*$ commutes with the $e_i$'s and their adjoints. 
Since the vacuum vector is characterised up to a scalar multiple by the conditions $e_i\Omega=0$ for $i\ge 0$ and
$e_i^*\Omega=0$ for $i<0$, $V\Omega=\zeta \Omega$ for some $\zeta\in{\Bbb T}$. Since $\Omega$ is cyclic for the
$e_i$'s and their adjoints, it follows that $V=\zeta I$. The shift operator has the following properties, which show in particular that it defines by conjugation an automorphism of $\a$:
\vskip .1in
\noindent \bf Lemma~A. \it (a) $Ue_iU^*=e_{i+1}$

\noindent (b) $UL_kU^*=L_k +a_k +{1\over 2}\delta_{k,0} I$

\noindent (c) $Ua_nU^*=a_n +\delta_{n,0}I$

\noindent (d) The system
$\a$, $U$ with the relations $Ua_nU^*=a_n +\delta_{n,0}c$,$UdU^*=d +a_0 +{1\over 2}c$
and  $d=a_0^2/2 +\sum_{n>0} a_{-n}a_n$
has a unique irreducible positive energy representation with $c$ acting as the scalar $I$ and $a_0$ having integer eigenvalues.

\noindent (e) If $m$ is a positive integer, the system
$\a$, $V$ with the relations $Va_nV^*=a_n +m\delta_{n,0}c$,
$VdV^*=d +ma_0 +{m^2\over 2}c$ and $d=a_0^2/2 +\sum_{n>0} a_{-n}a_n$
has $m$ inequivalent irreducible positive energy representations 
with $c$ acting as the scalar $I$ and $a_0$ having integer eigenvalues. These representations 
can be characterised as having a lowest energy vector of charge $i$ and energy $i^2/2$ with $i=0,\dots, m-1$. For each $i$, the other eigenvalues of $a_0$ are congruent to $i$ modulo $m$.
\vskip .05in
\noindent \bf Proof. \rm (a) has already been proven. For (b), note that 
$$[L_k-UL_kU^*,e_i]=-(i+(k+1)/2)e_{i+k} +(i+(k-1)/2)e_{i+k}=-e_{i+k}=[-a_k,e_i].$$
Hence $UL_kU^*-L_k-a_k$ commutes with the $e_i$'s; the same relation for $-k$ in place of $k$, shows that it also
commutes with their adjoints. Thus it carries $\Omega$ onto a scalar multiple of itself and by cyclicity is equal to 
mutliplication by this scalar. The scalar need only be computed for $k\ge 0$, by hermiticity. When $k>0$ we get $0$ while
for $k=0$, we get $UL_0U^*=L_0+a_0+1/2$. For (c), note that $Ua_nU^*-a_n$ commutes with the
$e_i$'s and their adjoints so is a scalar operator, which only need be calculated for $n\ge 0$ by looking at how it acts
on $\Omega$. For $n>0$ we get $0$, while for $n=0$, we get $Ua_0U^*-a_0= I$. Finally to prove (d), in the inner product space $\H$ take a vector $\xi_0$ of lowest energy which is an eigenvector of $a_0$. Thus
$d\xi=h\xi$ and $a_0\xi=\mu\xi$.  Then
the vectors $\xi_n =U^{-n} \xi$ 
satisfy $a_0\xi_n=(\mu+n)\xi_n$ and $d\xi_n= {1\over 2} n^2 \xi_n$. Each $\xi_n$ generates
an irreducible $\a$--module, necessarily mutually orthogonal since the eigenvalues of $a_0$ are different on each.  By irreducibility their direct sum gives the 
whole inner product space. But this decomposition shows that, after subtracting $hI$ from 
$d$ and $\mu I$ from $a_0$, $\H$ can be written 
as a tensor product $\K\otimes E$ 
where $\K$ is the standard irreducible representation of $(a_n)$ ($n\ne 0$) and $d$ and $E$ is the standard inner product space with orthonormal basis $e_n$ on which 
$Ue_n=e_{n+1}$, $de_n= {1\over 2} n^2 \xi_n$ and $a_0e_n=ne_n$. 
To prove (e), note that taking $V=U^m$, the space $E={\Bbb C}[{\Bbb Z}]$ splits up into
$m$ cosets for the subgroup $m{\Bbb Z}$ corresponding to the subgroup generated by $V$. Moreover each coset gives an irreducible summand with the vectors $\Omega_i$ ($i=0,\dots,m-1$) giving the corresponding lowest energy vectors. The representations are inequivalent because 
$i$ is specified as the smallest non--negative solution $\lambda$ of $a_0v=\lambda v$, $L_0v={1\over 2} \lambda^2v$. 
\vskip .1in 
\noindent \bf Corollary. \it In the setting of (e), the Virasoro operators assocyaed with the $a_n$'s satisfy $VL_kV^*=L_k + ma_k +{m^2\over 2}\delta_{k,0} I$.

\vskip .05in
\noindent \bf Proof. \rm It suffices to note that this holds for $V=U^m$. 
\vskip .1in
\rm
Now for $m\in {\Bbb Z}$ we define the formal power series $z$ and $z^{-1}$
$$\Phi_m(z)=U^{-m} z^{ma_0} \exp(\sum_{n>0} {m\cdot z^n a_{-n}\over n}) 
\exp( \sum_{n<0} {m\cdot z^n a_{-n}\over n})=\sum \Phi_m(n) z^{-n}.$$
If we define 
$$E_\pm^m(z)=\exp( \sum_{\pm n>0} {-m\cdot z^{-n} a_{n}\over n}),$$
then 
$$\Phi_m(z)=U^{-m} z^{ma_0} E_-^m(z) E_+^m(z).$$
To see that this is well defined, note that, when applied to a finite energy vector $\xi$, 
this becomes a formal power series in $z$ multiplied by a power of $z^{-1}$, 
since $E_+^m(z)\xi$ is a polynomial in $z^{-1}$ with vector coefficients.
We shall simply write $E_{\pm}(z)$ for $E_\pm^1(z)$. The notation has been chosen consistently so that, for $m>0$, 
$$E_\pm^m(z)=E_{\pm}(z)^m.$$
\vskip .1in
\noindent \bf Fubini--Veneziano relations. \it 
$[a_i,\Phi_m(z)]=mz^i \Phi_m(z)$ and $[L_k,\Phi_m(z)]=z^{k+1}\Phi^\prime_m(z) +{m^2\over 2}(k+1)z^k\Phi_m(z)$.
\vskip .1in

\rm This result can be proved in a number of ways, the most conceptual using operator product expansions, which is 
now usually formulated within the language of vertex algebras. The original proof was along those lines and 
is described in {\bf [14]}. Here we give a very elementary proof relying on the 
boson--fermion duality and no detailed computations. 
In Appendix~A we give two other proofs of the relations: the first is a simple
indirect proof based on a uniqueness result
for primary fields associated with the system $(a_n), U, L_0$; the second is a  direct verification by brute force 
along the lines sketched by Frenkel and Kac {\bf [10]}. 

\vskip .1in
\noindent \bf Lemma~B. \it Let $A$ be a formal power series in $z$ (or $z^{-1}$) with operator coefficients and let 
$D$ be an operator such that $[D,A]$ commutes with $A$. Then $[D,e^A]=[D,A]e^A$.
\vskip .05in
 \noindent \bf Proof. \rm We have $[D,A^N]=\sum_{p+q=N-1} A^p[D,A]A^q=N[D,A]A^{N-1}$.
\vskip .1in
\noindent \bf Corollary. \it (a) $[a_n,E_+(z)]=0=[a_{-n}, E_-(z)]$ if $n\ge 0$.

\noindent (b) $[a_n,E_-(z)]= z^n E_-(z),\,\, [a_{-n},E_+(z)=]z^{-n}E_+(z)$ if $n>0$.

\noindent (c) $[L_0,E_+(z)]=(\sum_{n>0} a_{n} z^{-n})E_+(z)$, $[L_0,E_-(z)]=(\sum_{n>0} a_{-n} z^{-n})E_-(z)$.

\vskip .05in
\noindent \bf Proof. \rm These formulas are straightforward consequences of the lemma, 
setting $D=a_n$ or $L_0$ and 
$$A=\sum_{\pm n>0} {m\over n} a_nz^{-n}.$$
\vskip .1in
\noindent \bf Proof of the Fubuni--Veneziano relations. \rm If $i\ne 0$, then
$a_i$ commutes with $U$ and $A_0$, while 
$$[a_i,E_-(z)E_+(z)]=mz^iE_-(z)E_+(z).$$
Hence $[a_i,\Phi_m(z)]=mz^i\Phi_m(z)$.
For $i=0$, we have $U^m a_0U^{-m} =a_0 + mI$ and hence $[a_0,U^{-m}]=mU^{-m}$, 
so that $[a_0,\Phi_m(z)]=m\Phi_m(z)$,
It follows that $\Phi_m(n)$ carries ${\cal F}[k]$ into ${\cal F}[k+m]$.  

If $k=0$, $UL_0U^*=L_0+a_0 + {1\over 2}I$, so that it follows by induction that 
$$U^mL_0U^{-m}=L_0 +m a_0 + {m^2\over 2}I.$$ 
Premultiplying by $U^{-m}$, we get
$$[L_0,U^{-m}]=U^{-m}(ma_0 +{m^2\over 2}I).$$ 
Hence
$$[L_0,\Phi_m(z)]= U^{-m} z^{ma_0}(ma_0 +{m^2\over 2}) E_-(z) E_+(z) + U^{-m} z^{ma_0} (B_-E_-E_+ +E_-B_+E_+),$$
where $B_\pm(z)=\sum_{\pm n>0} ma_{n}z^{-n}$. On the other hand
$$z{d\Phi_m(z)\over dz}=U^{-m} (ma_0) z^{ma_0}E_-E_+ +U^mz^{-ma_0}(B_-E_-E_++E_-B_+E_+).$$
Thus
$$[L_0,\Phi_m(z)]= z\Phi^\prime_m(z) +{m^2\over 2} \Phi_m(z).$$
We also have the relations
$$U\Phi_m(z)U^*=z^{m}\Phi_m(z),\,\,\,z^{-km}\Phi_m(z)\Omega_k|_{z=0}=\Omega_{k+m}.$$
The first is an immediate consequence of the relation $Ua_0U^*=a_0+I$. The second
follows from the first and the relation $\Phi_m(z)\Omega|_{z=0}=\Omega_m$. 

We shall say that
$\phi(n):{\cal F}\rightarrow {\cal F}$ is a primary field for the system of operators 
$(a_n)$, $U$ and $L_0$ if
$$[a_n,\phi(k)]=m\phi(k+n),\,\,[L_0,\phi(n)]=-(n+\mu)\phi(n),\,\, U\phi(n)U^*=\phi(n+m),$$
for some $m\in{\Bbb Z}$ and $\mu\in {\Bbb R}$.
In particular $\phi(n)=\Phi_m(n)$ satisfies these conditions with $\mu=-m^2/2$.
Note that if $(\phi(n))$ is a primary field, so too is 
the shifted field $\psi(n)=\phi(n+k)$ for any $k\in {\Bbb Z}$. We now show that after a suitable shift $\phi(n)$ is 
proportional to the modes of the field $\Phi_m(z)$. Since $U$, $a_n$ act 
cyclicly on $\Omega$, if $\phi(n)\Omega=0$, 
then $\phi$ is identically zero. Take $n$ minimal subject to $\xi \phi(n)\Omega\ne 0$: this is possible because 
${\cal F}$ is a positive energy representation. Shifting $\phi$ if necessary we many assume $n=0$.
Hence $a_i\xi=0$ for all $i>0$. Morevoer $a_0\xi=m\xi$. It follows
that $\xi$ is proportional to $\Omega_m$. Scaling if necessary, we may assume that $\xi=\Omega_m$. 
Taking $k=n$ in $[L_0,\phi(k)]=-(k +\mu)\phi(k)$, we get $\mu=m-m^2/2$. 
We have the following uniqueness result for primary fields:
\vskip .1in
\noindent \bf Lemma~C. \rm \it (a) Let $\K_i$ be irreducible representations of $\a$ with lowest energy vectors $\xi_i$ with $a_0\xi_i=\lambda_i\xi_i$. Let
$\phi(n):\K_1\rightarrow \K_2$ be operators such that
$[d,\phi(n)]=-(n+\mu)\phi(n)$ and $[a_k,\phi(n)]=\nu \phi(n+k)$ for $\nu\in {\Bbb R}$ fixed.
If $(\phi(n)\xi_1,\xi_2)=0$ for all $n$, then $\phi(n)\equiv 0$.

\noindent (b) Let $(a_n)$, $d$, $U$ be the bosonic system of operators acting on ${\cal F}$ and suppose that $\phi(n):{\cal F} \rightarrow {\cal F}$ satisfies
$[d,\phi(n)]=-(n+\mu)\phi(n)$, $[a_k,\phi(n)]=\nu \phi(n+k)$ and $U\phi(n)U^*=\phi(n+a)$ for some $\mu,\nu,a\in {\Bbb Z}$. 
Then if $\phi(n)\Omega$, we have $\phi(n)\equiv 0$. In particular if 
$(\phi(n)\Omega,\Omega_\nu)=0$, then $\phi(n)\equiv 0$.

\noindent (c) Let $(a_n)$, $d$, $U$ be the bosonic system of operators acting on ${\cal F}$
and define $a_n^{(1)}=a_n\otimes I$, $a_n^{(2)}=I\otimes a_n$, $U_1=U\otimes I$, $U_2=I\otimes U$ on (the ${\Bbb Z}_2$--graded tensor product) ${\cal F}^{\otimes 2}$. Suppose that
 $\phi(n):{\cal F}^{\otimes 2} \rightarrow {\cal F}^{\otimes 2}$ satisfies
$[d,\phi(n)]=-(n+\mu)\phi(n)$, $[a_k^{(i)},\phi(n)]=\nu_i \phi(n+k)$ 
and $U_i\phi(n)U_i^*=\phi(n+a)$ for some $\mu,\nu,a\in {\Bbb Z}$. 
Then if $\phi(n)\Omega\otimes \Omega=0$, 
we have $\phi(n)\equiv 0$.

\vskip .05in
\noindent \bf Proof. \rm (a) From the commutation conditions
it follows that $(\phi(n)p(a_{-1},a_{-2},\cdots))\xi_1,\xi_2)=0$ for any polynomial $p$.
Thus $(\phi(n)\xi,\xi_2)=0$ for all $\xi$. But then by the commutation relations
$(\phi(n)\xi,p(a_{-1},a_{-2},\cdots) \xi_2)=0$ for any polynomial $p$. But then
$(\phi(n)\xi,\eta)=0$ for all $\xi,\eta$ and hence $\phi(n)\equiv 0$.

\noindent (b) Note that by (a) it suffices 
to show that $(\phi(n)\Omega_r,\Omega_s)=0$ for all $r,s$,
But then
$$(\phi(n)U^r \Omega,U^s \Omega)= (\phi(n-ar)\Omega,U^{s-r}\Omega)=(\phi(n-ar)\Omega,\Omega_{r-s})=0.$$
Finally $\phi(n):{\cal F}_0\rightarrow {\cal F}_\nu$, so by (b), the second condition implies that $\Phi(n)\Omega=0$.

\noindent (c) Let $\Omega^\prime$ be an eigenvector of $d$ in ${\cal F}$ and define $\psi(n)$ by
$$(\psi(n)\xi,\eta)=(\phi(n)\Omega\otimes \xi, \Omega^\prime\otimes \eta).$$
Then by $\psi(n)$ satisfies the conditions in (b), so $\psi(n)=0$. Since $\Omega^\prime$ was
arbitrary, it follows that $\phi(n)\Omega\otimes \xi=0$ for all $\xi$. Fixing $\xi_2,\eta_2$ and defining
$\psi^(n)$ by
$$(\psi^\prime(n)\xi,\eta)=(\phi(n)\xi\otimes \xi_2,\eta\otimes \eta_2),$$
we deduce that $\psi^\prime(n)=0$. Hence $\phi(n)\equiv 0$.
\vskip .1in
This uniqueness result allows us to prove the boson--fermion correspondence in its vertex 
operator formulation. 
\vskip .1in
\noindent \bf Example~1. \rm $\Phi_1(z)=\sum e_n z^{-n-1}$ and $\Phi_{-1}(z)=\sum e_{-n}^*z^{-n}$. Indeed both sides of these identities satisfy the same covariance relations with
$\Phi_1(-1)\Omega=\Omega_1=e_{-1}\Omega$ and $\Phi_{-1}(0)\Omega=\Omega_{-1}=e_0^*\Omega$. 
\vskip .1in
As a consequence we get
$$[L_n,\Phi_{\pm 1}(z)]=z^{n+1}\Phi{\pm_1}^\prime (z) + {1\over 2}(n+1)z^n \Phi_{\pm 1}(z).\eqno{(*)}$$
We claim that this implies the Fubini--Veneziano relations for $\Phi_{\pm m}(z)$.
Indeed $(*)$ immediately implies the relations for the $m^2$--fold 
${\Bbb Z}_2$--graded tensor product 
$$\Theta_{\pm}(z)
=\Phi_{\pm 1}(z) \otimes \Phi_{\pm 1} (z)\otimes \cdots \otimes \Phi_{\pm 1}(z)$$
on ${\cal F}^{\otimes m^2}$. But this is just the vertex operator constructed with
the operators 
$$A_n=a_n\otimes I \otimes I\otimes \cdots \otimes I +\cdots + I \otimes \cdots \otimes I \otimes a_n$$
and
$$V=U\otimes U \otimes \cdots \otimes U$$
on ${\cal F}^{\otimes m^2}$. Now
$$[A_p,A_q]=m^2 p\delta_{p+q,0} I, \,\,\,\, VA_0V^*=A_0 + m^2 \cdot I.$$
Thus letting $a_n^\prime =m^{-1}A_n$, we have
$$[a_p^\prime,a_q^\prime]
=p\delta_{p+q,0} I, \,\,\,\,Va_0^\prime V^*=a_0^\prime + m\cdot I.$$
The operators $L_k$ on  ${\cal F}^{\otimes m^2}$ satisfy
$$[L_k,a_n^\prime]=-na_{n+k}^\prime,\,\,\,\,VL_kV^*=L_k + ma_k^\prime +{m^2 \over 2}\delta_{k,0}I.\,\,\,\, [L_k,\Theta_{\pm}(z)]=z^{k+1}{d\over dz}\Theta_{\pm}(z) + {m^2\over 2}(k+1) z^k\Theta_{\pm}(z).$$
The first two relations are also true for the Virasoro operators $L^\prime_k$ 
constructed from the $a_n^\prime$'s, by the corollary to Lemma~A. Thus 
$L_k-L_k^\prime$ commutes with $a^\prime_n$ and $V$, and therefore the third idenity is valid with $L^\prime_k$ replacing $L_k$. This holds on any irreducible
submodule for the $a_n^\prime$'s and $V$, in particular those on which $a_0^\prime$ has integer eigenvalues. 
It is easy to see that the eigenvalues of $A_0$ run through all possible integers, 
so that those of $a_0^\prime=m^{-1}A_0$ in particular include all integers. Thus the Fubini--Veneziano relations are satisfied on any irreducible representation of the type discussed in
part (e) of Lemma~A. But these are exactly the
irreducible representations for $U^m$ and $\a$ that occur in ${\cal F}$, 
so the Fubini--Veneziano relations hold on ${\cal F}$.
\vskip .1in
We now pass to ${\cal F}_2={\cal F}\otimes {\cal F}$. The unitary operators 
$V=U\otimes U^*$ and $W=U\otimes U$ act on this space. 
Now that $V$ commutes with the action of $d_n=a_n\otimes I +I\otimes a_n$ and therefore leaves its 
mutliplicity space invariant, since these can be indentified with its space of singular vectors as explained above.
More explicitly invariant under  $b_n=a_n\otimes I - I\otimes a_n$ is invariant under $V$.
Similarly any module invariant under $b_n$ is invariant under $W$.

As a consequence the two level one representations of $\widehat{\s}$ are invariant under the vertex operators
$$\Psi_m(z)=\Phi_m(z)\otimes \Phi_{-m}(z)=V^{-m} z^{mb_0}
\exp(\sum_{n>0} {z^{n} m b_{-n}\over n}) \exp( \sum_{n<0} {m z^n b_{-n}\over n})=\sum \Psi_m(n) z^{-n}.$$
The two systems $(b_n),\, V$ and $(d_n),\, W$ are examples of systems $(A_n),\,U$
$$[A_m,A_n]=2m\delta_{m+n,0}\cdot I,\,\,\,\, UA_mU^*=A_m +2\delta_{m,0} I.$$
Thus if $L_0={1\over 2 } A_0^2 +\sum_{n>0} A_{-n}A_n$, then we have 
$$UL_0U^*  -L_0={1\over 2}[(A_0+2I)^2-A_0^2)=A_0 + 2I.$$
Such a system has two inequivalent positive energy representations with $A_0$ acting
with integer eigenvalues, either all odd or all even. Thus both sets of  operators $H(n),V$ and
$K(n),w$ have only two irreducible 
representations with $H(0)$ and $K(0)$ having half integer
eigenvalues, In one the eigenvalues run through the whole of ${\Bbb Z}$, and in the other through ${1\over 2} +{\Bbb Z}$. It follows that operators $H(n),V$ act irreducibly on each $\H_j$ and  that the $K(n),W$ act irreducibly on the associated multiplicity spaces.
The operator $V$ also has a compatibility with the operators $E(n)$ and $F(n)$ on ${\cal F}^{\otimes 2}$:
$$VE(n)V^*=E(n+2),\,\,\,\, VF(n)V^*=F(n-2).$$
In fact more generally
$$ U_1 E(n) U_1^*=E(n+1),\,\, U_2E(n)U_2^*=E(n-1),\,\, U_1 F(n) U_1^*=F(n-1),\,\, U_2F(n)U_2^*=F(n+1).$$ 
Indeed it suffices to show the result for $E(n)$ since $E(n)^*=F(-n)$. The differences 
$X=U_1E(n)U_1^*-E(n+1)$ and $Y= U_2E(n)U_2^*-E(n-1)$ 
commute with the fermionic operators $v_i$ and $v_i^*$, so that 
$X(\Omega\otimes\Omega)$ and $Y(\Omega\otimes \Omega)$ are proportional to $\Omega\otimes\Omega$. Hence $X$ and $Y$ are scalar operators. Since 
$[H,X]=X$ and $[H,Y]=Y$, we must have $X=0=Y$. 

\vskip .1in
\noindent \bf Example~2 (Frenkel--Kac--Segal). \rm $E(z)=\Psi_1(z)$ and $F(z)=\Psi_{-1}(z)$ on ${\cal F}^{\otimes 2}$ 
and hence on the representations $H_0$ and $H_{1/2}$. Indeed the result follows from part (c) of Lemma~C since, as we have just checked,
both sides of these identities satisfy the same covariance relations and
$$\Psi_1(0)(\Omega\otimes \Omega)=\Omega_1\otimes \Omega_{-1}=E(-1)(\Omega\otimes \Omega),
\,\,\,\,\Psi_{-1}(0)(\Omega\otimes \Omega)=\Omega_{-1}\otimes \Omega_1=F(-1)(\Omega\otimes \Omega).$$ 
\vskip .1in
Graeme Segal found explicit formulas for the Goldstone vectors in the vacuum representation $\H_0$. 
The same method can be applied in representation $\H_{1/2}$. Let $x_i=b_{-i}$ and define operators 
$c_n$ ($n>0$) by the formal power series expansion
$$\Psi_-(z)=\exp \sum_{n>0} {b_{-n}\over n} z^{n} =\sum_{n\ge 0} c_n z^n.$$
Thus $c_0=1$. For $n<0$, we define $c_n=0$. 
Similarly we can define 
$$\Psi_+(z)=\exp \sum_{n<0} {b_{-n}\over n} z^{n},$$
so that 
$$\Psi_m(z)=V^{-m} z^{mb_0} \Psi_-(z)^m\Psi_+(z)^m.$$
Evidently
$$\Psi_+(z)=E_+(z)\otimes E_+(z)^{-1},\,\,\,\, \Psi_-(z)=E_-(z)\otimes E_-(z)^{-1}.$$
For any signature $f_1\ge f_2\ge f_3 \ge 0$, we define 
$$X_f=\det c_{f_i-i+j}.$$
\vskip .1in
\noindent \bf Remark. \rm The proofs below can be understood in terms of the combinatorics of the Weyl character formula
for $U(n)$. The characters of the irreducible representations $\pi_f$ of $U(n)$ are indexed by signatures 
$f_1\ge f_2 \ge \cdots f_n$ with the character given by 
$\chi_f(g)= {\rm Tr}\, \pi_f)(g)=\det z_i^{f_j + n-j}/\det z_i^{n-j}$ if $z_1,z_2,\dots,z_n$ are the eigenvalues of $g$. 
The denominator is the Vandermonde determinant given by $\Delta(z)=\prod_{i<j} (z_i-z_j)$. 
The group $U(n)$ acts irreducibly on $S^k{\Bbb C}^n$: it is the representation 
with signature $(k,0,0,\dots,0)$. The group $U(m)\times U(n)$ acts on $S^k({\otimes C}^m\otimes {\otimes C}^n)$ 
and the character if this representation on an element $(g,h)$ 
is given by $\sum_{f_n\ge 0,|f|=k }\chi_f(g)\chi_f(h)$, where $|f|=f_1+\cdots+f_n$. From this it follows that 
$\chi_f(g) =\det C_{f_i-i+j}$ (the Jacobi-Trudy identity) if $C_k=\chi_{(k,0,\dots,0)}(g)$. 
Moreover the following generating function identity holds:
$$\prod_{i=1}^n f(z_i)=\sum_{f_n\ge 0} \chi_f(z)\cdot \det c_{f_i-i+j},$$
where $f(z)=1+c_1z + c_2 z^2 + \cdots$ (with $c_0=1$ and $c_i=0$ for $i<0$).
Note that the generating function for the $X_j$'s giving the characters has $X_j=\sum z_i^j$. 
\vskip .1in
\noindent \bf Lemma~D. \it If $A$ and $B$ are formal power series of operators in $z$ and $w^{-1}$ such that $C=[A,B]$ is
a multiple of the identity operator, then $e^Ae^B=e^C e^Be^A$.
\vskip .05in
\noindent \bf Proof. \rm By Lemma~B, $[e^A,B]=Ce^A$ so that $e^ABe^{-A}=B+C$. Hence
$$e^Ae^Be^{-A}=\exp e^ABe^{-A}=e^{B+C}=e^C e^B.$$ 
\vskip .1in
\noindent \bf Corollary. \it (a) $E_+(z)E_-(w)=(1-w/z)^{-1} E_-(w)E_+(z)$.

\noindent (b) $E_+(z)^{-1}E_-(w)^{-1}=(1-w/z)^{-1} E_-(w)^{-1}E_+(z)^{-1}$.
\vskip .05in
\noindent \bf Proof. \rm To prove (a), let 
$A=\sum_{n>0} a_n z^{-n}/n$ and $B=\sum_{n<0} a_nw^{-n}/n$. Then
$$C=[A,B]=  \sum_{n>0} z^{-n}w^{-n}/n=-\log (1-{w\over z}),$$
so that 
$$e^C=\exp -\log(1-{w\over z})=(1-{w\over z})^{-1}.$$
(b) can be proved similarly or by taking formal inverses in (a). 
\vskip .1in
We now establish the formulas for Goldstone's vectors. 
As in Section~4, where the Goldstone vectors are defined, let 
$\xi_j$ be a lowest energy unit vector in $\H_j$ with $H\xi_j=j\xi_j$. Then for $m\in {\Bbb Z}$ we may take $\xi_{m+j}=V^{-m}\xi_j$. 
\vskip .1in
\noindent \bf Lemma~E (Goldstone's formulas). \it If $k\ge 0$ is a half integer and $m\ge 0$ is an integer, the Goldstone vector in $ \H_j[- k]$ of energy $(k+m)^2$ 
is proportional to $X_f \xi_{-k}$ where $f_1=2k+m,\dots, f_m=2k+m$ and $f_i=0$ for $i>m$ and the Goldstone vector in $\H_j[k]$ of energy $(k+m)^2$ is proportional to
$X_{f^\prime}\xi_k $ where $f^\prime_1=m,\dots,f^\prime_{2k+m}=m$ and $f^\prime_i=0$ for $i>2k+m$ (the transposed signature).
\vskip .05in 
\noindent \bf Remark. \rm In {\bf [32]}, Segal only treats the case $j=0$ and his
statement is not quite correct since it omits the transpose. Note that the passage from $-k$ to $k$ is equivalent to replacing the variables $X_n$ by $-X_n$ (and hence inverting the generating series $\sum c_n z^n$). It corresponds to the action of the Weyl
group element $s=\pmatrix{0 &1\cr-1 & 0\cr}$ which establishes an 
isomorphism between $\H_j[-k]$ and $\H_j[k]$. Indeed the automorphism 
$b_n\mapsto -b_n$ fixes the Virasoro algebra generators 
and changes the sign of the eigenvalue of $b_0$.   
 
\vskip .05in
\noindent \bf Proof. \rm By the corollary to Lemma~D, 
$$\Psi_+(z)\Psi_-(w)=(1-{w\over z})^{-2}\Psi_-(w)\Psi_+(z).$$
On the other hand the Goldstone vector $\xi_+\in \H_j[k]$ is proportional to 
$E(0)^{m+2k}\xi_{-m-k}$ and the Goldstone vector $\xi_-\in \H_j[-k]$ is proportional
to $E(0)^m\xi_{-m-k}$. Let $p=p_\pm=m +k \pm k$. Thus $\xi_\pm$ is proportional to
$E(0)^{p_\pm}\xi_{-m-k}$. 
Now $\Psi(z)\equiv \Psi^1(z)=\sum E(n) z^{-n-1}$. Hence $\xi_\pm$ is proportional to the 
coefficient of $z_1^{-1}\dots z_p ^{-1}$ in 
$\Psi(z_1)\dots \Psi(z_{p})\xi_{-k-m}$. On the other hand if $C(z)=\sum c_n z^n$, then
$$\eqalign{\Psi(z_1)\dots \Psi(z_p)\xi_{-k-m}
&=\prod_{i<j} (1-{z_i\over z_j})^2 
z_1^{2(p-1)}z_2^{2(p-2)}\cdots z_p^2 (z_1\cdots z_p)^{-2k-2p}\Psi_-(z_1) \cdots \Psi_-(z_p)\xi_{\pm k}\cr
&=\prod_{i<j} (z_1z_2\cdots z_p)^{-2k-2m+p} C(z_1)\dots C(z_p)\prod_{i<j}(z_i-z_j)(z_i^{-1} -z_j^{-1})\xi_{\pm k}.\cr}$$
Now we already noted that
$$C(z_1)\dots C(z_p)=\sum_{g_{p+1}= 0} \chi_g(z)\cdot \det c_{g_i-i+j}.$$
Thus we are looking for the constant coefficient in
$$(z_1z_2\cdots z_p)^{-2k-2m+p} C(z_1)\dots C(z_p)\prod_{i<j}(z_i-z_j)(z_i^{-1} -z_j^{-1})\xi_{\pm k}.$$
By the combinatorics of the Weyl integration formula,
the constant coefficient in
$$\chi_f(z^{-1}) \chi_g(z) \prod_{i<j}(z_i-z_j)(z_i^{-1} -z_j^{-1})$$
equals $p!\delta_{f,g}$ (the usual orthogonality relations). Thus $\xi_-$ is proportional
to $X_f\xi_{-k}$  with $f_1=2k+m,\dots, f_m=2k+m$ and $f_i=0$ for $i>m$; and 
$\xi_+$ is proportional to $X_{f^\prime}\xi_k$ with 
$f^\prime_1=m, \dots, f^\prime_{2k+m}=m$ and $f^\prime_i=0$ for $i>2k+m$. 
Thus $f^\prime$ is the transposed Young diagram. 

Following Segal, we now check this directly 
without using any unproved combinatorial results on symmetric functions. We have 
$$\Delta(z^{-1})=\det z_i^{-(p-j)}$$
and similarly
$$\Delta(z)=\sum_{\sigma\in S_p} \varepsilon(\sigma) 
z_{\sigma(1)}^{p-1} z_{\sigma(2)}^{p-2} \dots z_{\sigma(p-1)}.$$
We are looking for the constant term in
$$C(z_1)\cdots C(z_p)(z_1\cdots z_p)^{-2k-2m+p} \Delta(z^{-1})\sum_{\sigma\in S_p} 
\varepsilon(\sigma)\sigma[z_1^{p-1} \cdots z_{p-1}].$$
Since $\Delta(z^{-1})$ satisfies $\sigma\Delta(z^{-1})=\varepsilon(\sigma)\Delta(z^{-1})$ and the other terms
are invariant under permutation, this is the same as the constant term in
$$\eqalign{p! C(z_1)\cdots C(z_p)& (z_1\cdots z_p)^{-2k-2m+p} \det (z_i^{-(p-j)})\,z_1^{p-1} \cdots z_{p-1}\cr
&=p!C(z_1)\cdots C(z_p) \det z_i^{j-2k  -2m +p -i}\cr}$$
But, if the determinant is written as an alternating sum over the symmetric group $S_p$,
this is clearly proportional to $\det c_{2k+2m-p +i-j}$, as required. 
\vskip .1in

\noindent \bf Remark. \rm The proof above is formally similar 
to the proof of the Weyl character formula
for $U(N)$ (or $SU(N)$). Wallach made this link explicit by noticing that the vertex operator formalism
allows a map to be defined from central functions on $U(n)$ into the oscillator representations. This map
intertwines the action of Virasoro operators with the natural action of the operators  
$\ell_{n}$ ($n>0$) on central functions given by
by $\ell_{n}= {\cal D}_n +\alpha {\rm Tr}(g^n)$ where  ${\cal D}_nf(g)={d/dt}f(g\exp{-tg^n})|_{t=0}$. 
\vskip .1in
We will now exhibit explicitly primary fields by compressing the primary fields
$\Psi_k(z)$ on $\H_j$ between Goldstone vectors. This will reduce to the computation
of the expressions of the form $L_1^{|f|} X_f\xi_{p}$. 
Note that since $X_f$ has total degree $f$, we have
$${1\over |f|!} L_1^{|f|} X_f\xi_{p}=\alpha_f(p)\xi_{p},$$
where $\alpha_f(p)$ is a constant, which we now determine through a product formula.
\vskip .1in
So we consider now the 
action of $L_{1}$ in the oscillator representation with $b_0=2\mu I$, so that$H(0)=\mu I$.
This can be identified with the action on ${\Bbb C}[X]$ with $b_{-n}=\sqrt{2}X_n$ 
and $b_n=\sqrt{s}2n\partial_{X_n}$ for $n>0$. The operator $L_1$ is given by the formula
$$\eqalign{L_1&={1\over 4} \sum_{p+q=1} b_pb_q\cr
&={1\over 2} b_{0}b_{1} + {1\over 2}\sum_{p\ge 1}  b_{-p} b_{p+1}\cr
&=\mu \partial_{X_1} + \sum_{n\ge 2}n X_{n-1}\partial_{X_n}.\cr}$$
This operator can thus be described in terms of two derivations
$D_1$ and $D_2$ on the ring ${\Bbb C}[X]$, namely $D_1X_n=nX_{n-1}$ ($n\ge 2$), $D_1 X_1=0$ and $D_2=\partial_{X_1}$. 
Now if $f(z)=\exp \sum_{n>0} X_nt^n/n=\sum c_nt^n$, then $D_1f=(\sum_{n>0} X_nt^{n+1})f=t^2f^\prime(t)$. It follows
that $D_1c_n= (n-1)c_{n-1}$. Similarly $D_2f=tf$, so that $D_2c_n=c_{n-1}$. Let $D=D_1+\mu D_2$, also a derivation.
Since $X_f$ has degree $|f|$ and $D$ decreases the degree by $1$, 
it follows that $D^{|f|}X_f$ is a scalar.
\vskip .1in
\noindent \bf Lemma~F. \it $(|f|!)^{-1}D^{|f|} \det c_{f_i-i+j}= \det {\mu+f_i-i+j-1\choose f_i -i +j}$.
\vskip .05in
\noindent \bf Proof. \rm  
Consider a term in the determinant $c_{m_1}\cdots c_{m_n}$. The term $c_{m_i}$ has degree $m_i$
so the product has degree $\sum m_i=|f|$. Since $Dc_k=(k-1+\mu)c_{k-1}$, we have
$$D^k c_k =\mu(\mu+1) \cdots (\mu+k-1)=k!\cdot {\mu-1 +k\choose k}.$$
Hence by the Leibniz rule,
$$\eqalign{D^{|f|} c_{m_1}\cdots c_{m_n}&={|f|!\over m_1! \cdots m_n!}D^{m_1}(c_{m_1})D^{m_2}(c_{m_2}) \cdots D^{m_n} (c_{m_n})\cr
  &=|f|! \prod_{i=1}^n {\mu -1 +m_i\choose m_i}.\cr}$$
So applying $(|f|!)^{-1}D^{|f|}$ to the determinant has the effect of replacing each term $c_k$ by ${\mu-1+k\choose k}$, as claimed.
\vskip .1in 
\noindent \bf Corollary. \it  If $p$ is a half--integer, then 
$(|f|!)^{-1}L_1^{|f|}(X_f \xi_{p}) =\alpha_p (f)\xi_{p}$ where 
$\alpha_p(f)= \det {2p+f_i-i+j-1\choose f_i -i +j}$. 
\vskip .1in
\rm We now calculate the above determinant for the signature $f_i=N$ ($i=1,m$), $f_i=0$ ($i>m$).
\vskip .1in
\noindent \bf Lemma~G. \it If $N\ge m$ then,
$$\det {\lambda+N-1+ j-i \choose N+j-i}_{i,j=1,m} 
= {\prod_{j=1}^{N} \prod_{i=1}^m (\lambda-i+j)\over \prod_{j=1}^{N} \prod_{i=1}^m (N+i-j)}.$$   
\vskip .05in
\noindent \bf First proof (Weyl's character and dimension formulas). \rm Recall that, 
if $f_1\ge f_2\ge \cdots \ge f_n\ge 0$ is a signature, then the character of the corresponding irreducible 
representation $\pi_f$ of the unitary group $U(n)$ is given by
$$\chi_f(g)={\rm Tr}\, \pi_f(g)={\det z_j^{f_i +n-i}\over \det z_j^{n-i}},$$
where $z_1,\dots,z_n$ are the eigenvalues of $g$. The signature defines a Young diagram; 
the hook length of any node in the Young diagram is defined to be one plus the sum 
of the number of nodes in the same column or row after that node. If we place $n-i+j$ in at the $j$th
node of the $i$th row, then the dimension $d_f$ of $\pi_f$ is the product of these numbers divided by the 
product of the hooklengths. It is equivalently given by the usual Weyl dimension formula
$$d_g ={\prod_{j<i} (f_i-f_j-i+j)\over (n-1)!!},$$
although this formula is not manifestly written as a polynomial in $n$.
If $f_1=r$ and $f_i=0$ for $i>n$, then $\chi_r$ is the character of the 
symmetric power $S^r {\Bbb C}^n$, which has dimension ${n+r-1\choose r}$.
On the other hand $\chi_f$ can be expressed as a determinant in the $\chi_r$'s:
$$\chi_f=\det \chi_{f_{i}+j-i},$$
where we set $\chi_r=0$ for $r<0$. 
In particular evaluating the character at $1$, we get a formula expressing a deteminant of binomial coefficients
as a product of linear factors in $n$. Taking a Young diagram with $m$ rows and $K+m$ boxes in each row,
we have $f_i=K+m$ ($1\le i\le m$) and $f_i=0$ for $i>m$. Hence
$$\det {n-1+K+m +j-i \choose K+m +j-i}_{i,j=1,\dots,m} 
={\prod_{j=1}^{K+m} \prod_{i=1}^m (n-i+j)\over \prod_{j=1}^{K+m} \prod_{i=1}^m (K+m+i-j)}.$$  
The stated determinantal formula follows by setting $\lambda=n$ and $N=K+m$.
\vskip .05in
\noindent \bf Remark. \rm 
Although we shall not need this result, the dimension formula implies more generally that
$\det {\mu+f_i-i+j-1\choose f_i -i +j}$ is the product of factors $\mu +j-i$, when there is 
a box in the $i$th row and $j$th column of the Young diagram, 
divided by the product of hooklengths of the diagram. 
In other words if $\mu=-p$, with $p\ge 1$ an integer, then
$(|f|!)^{-1}D^{|f|} X_f$ equals the dimension of the representation $\pi_f$ of $U(p)$. It is therefore stictly positive when $f_{p+1}=0$ and vanishes if $f_{p+1}\ne 0$. Since $L_{-1}\xi_0=0$, all these expressions vanish for $p=0$. 
\vskip .1in
\noindent \bf Second proof. \rm We have to show that
$$\det {a+N +i-j\choose N+i-j}_{i,j=1, m} = \prod_{p=1}^N \prod_{q=1}^m 
{a+1+p-q\over N-p+q}.$$ 
We subtract a multiple of the second column from the first column to make the 
last entry of the first column $0$; then we subtract a multiple of the third column from the second column to make the 
last entry of the second column $0$, and so on. After these modification
the $(i,j)$ entry for $i,j<n$ is modified from 
${N+a +i-j\choose a}$ to ${N-1+a +i-j\choose a}(m-i)/(N +m-j)$.
Thus writing $F_m(a,N)=\det {a+N +i-j\choose N+i-j}_{i,j=1, m}$, we have
$$F_m(a,N)={{a+N \choose N} \over {N+m\choose m-1}}\cdot F_{m-1}(a-1,N-1).$$
This may be rewritten as
$$F_m(a,N)= {(a+1)\cdots (a+N) \over m(m+1)\cdots (m+N-1)}\cdot F_{m-1}(a-1,N).$$
Thus 
$$F_m(a,N)=(m-1)!!{\prod_{i=1}^{m} \prod_{j=1}^N {a+1 +j-i}\over \prod_{i=1}^m \prod_{j=1}^N (N+i-j)}.$$
\vskip .1in
\noindent \bf Corollary. \it If $f_i=N$ for $1\le i\le m$ and $f_i=0$ for $i>m$, then 
$(|f|!)^{-1}L_1^{|f|}(X_f \xi_{p}) =\alpha_f(p) \xi_{p}$ 
where 
$$\alpha_p(f)={\prod_{j=1}^{N} \prod_{i=1}^m (p-i+j)\over \prod_{j=1}^{N} \prod_{i=1}^m (N+i-j)}.$$   
\vskip .1in
\noindent \bf Lemma~H. \it For $k$ an integer, $\Psi_k(z)\xi_0=e^{zL_{-1}}\xi_k$.
\vskip .05in
\noindent \bf Proof. \rm Let $f(z)=\Psi_k(z)\xi_0$. We have $[L_{-1},\Psi_k(z)]=d\Psi_k/dz$ and $L_{-1}\xi_0=0$.
Hence $df(z)/dz=L_{-1}f(z)$. Since $f(0)=\xi_k$, the result follows.
\vskip .1in
\noindent \bf End of proof. \rm We shall just use the fact that if $f$ is a signature as above with
$f_{p+1}=0$ with $p\ge 1$, then
$$(L_1^{|f|} X_f\xi_{p},\xi_{p})\ne 0.$$  
We fix an integer $m\ge 1$ and a half integer $k\ge 0$. Since we can take adjoints, it 
will suffice to show that if $i$ is a 
non--negative half integer with $i-k$ an integer and $|m-k|\le i\le m+k$, then there are
are Goldstone vectors $\xi_1,\xi_2$ in $\H_j$ with $\xi_1$ having energy $k^2$ or $i^2$ and
$\xi_2$ energy $i^2$ or $k^2$ such that $(\Psi_m(z)\xi_1,\xi_2)\ne 0$. In this case
if $P_1$ is the projection onto the $\v$--module generated by $\xi_i$, then
$P_2 \Psi_m(z) P_1$ or its adjoint is the 
required primary field of type $m^2$ between $L(1,k^2)$ and $L(1,i^2)$.

Assume first that $k$ is an integer and set $s=|k-m|$. Take $\xi_{k}$ in $H_0$ or $H_{1/2}$ and suppose $i$ satisfies $|k-m|\le i \le k+m$ and $i-k\in {\Bbb Z}$.
We assume further that $m\ge k$ and $i\ge m$ and set
$$A(z)=(\Psi_{m}(z)\xi_{-i},X_f\xi_{-i+m}),$$
where $X_f\xi_{-i+m}$ is a Goldstone vector of energy $k^2$. 
Thus $f_r=k+i-m$ for $r=1,\dots, k+m-i$ and $f_r=0$ for $r>k+m-i$. We must show
that $A(z)\ne 0$. But 
$$\eqalign{A(z)&=(\Psi_{m}(z)V^{i}\xi_0, X_fV^{i-m}\xi_0)\cr
&=z^{-2mi} (\Psi_{m}(z)\xi_0,X_f\xi_{m})\cr
&=z^{-2mi} (e^{zL_{-1}}\xi_{m},X_f\xi_{m})\cr
&=z^{-2mi} (\xi_{m},e^{zL_1}X_f\xi_{m}).\cr}$$
This is non--zero because $f_{m+1}=0$ since $m+1>k+m-i$. 

Now suppose that $m\ge k$ and $i\le m$. In this case we consider
$$B(z)=(\Psi_{m}(z)\xi_{-k}, X_f\xi_{m-k}),$$
where $X_f\xi_{m-k}$ is a Goldstone vector of energy $i^2$. Thus
$f_r=m-k+i$ for $r=1,\dots,i+k-m$ and $f_r=0$ for $r>i+k-m$. 
Then
$$\eqalign{B(z)&=(\Psi_{m}(z)V^{k}\xi_0, X_fV^{-m+k}\xi_0)\cr
&=z^{-2km} (\Psi_{m}(z)\xi_0,X_f\xi_{m})\cr
&=z^{-2km} (e^{zL_{-1}}\xi_{m},X_f\xi_{m})\cr
&=z^{-2km} (\xi_{m},e^{zL_1}X_f\xi_{m}),\cr}$$
which again is non--vanishing since  $m\ge i+k-m$ 
and hence $f_{m+1}=0$.

Thus if $k\le m$ we have shown that there is a non--zero primary
field of type $m^2$ from $L(1,k^2)$ to $L(1,i^2)$ if $V_i\le V_m\otimes V_k$. 
taking adjoints it follows that there is a non--zero primary field of type $m^2$
from $L(1,i^2)$ to $L(1,k^2)$ if $i\ge m$ and $k\le m$ and $V_k\le V_m\otimes V_i$. 
So now assume that $k\ge m$ and $i\ge m$ with $V_i\le V_m\otimes V_k$.  Since may take 
adjoints, we many assume that $i\ge k$. But then the
computation with $A(z)$ above shows that the corresponding primary field exists. 

Now suppose that $k\in {1/2} +{\Bbb Z}$. Let $W$ be 
a shift operator such that $W\xi_0=\xi_{1/2}$. Let $k=p+{1\over 2}$ and $i=q+{1\over 2}$.
If $k<m$ and $i>m$ or $k>m$ and $i>m$, then as before we take
$$A(z)=(\Psi_{m}(z)\xi_{-i},X_f\xi_{-i+m}),$$
with the same choice of $f$ as before. Then we have
$$\eqalign{A(z)&=(\Psi_{m}(z)V^{q}W^{}\xi_0, X_fV^{q-m}W \xi_0)\cr
&=z^{-2im} (\Psi_{m}(z)\xi_0,X_f\xi_{m})\cr
&=z^{-2im} (e^{zL_{-1}}\xi_{m},X_f\xi_{m})\cr
&=z^{-2im} (\xi_{m},e^{zL_1}X_f\xi_{m}),\cr}$$
which as before is non--zero. Thus the primary fields exist in this case. 
If $i<m$ and $k<m$ we take
$$B(z)=(\Psi_{m}(z)\xi_{-k}, X_f\xi_{m-k}),$$
with the same choice of $f$ as before. In this case
$$\eqalign{B(z)&=(\Psi_{m}(z)V^pW\xi_0,X_fV^{p-m}W\xi_0)\cr
&=z^{-2km}(\Psi_{m}(z)\xi_0,X_f\xi_{m})\cr
&=z^{-2km}(e^{zL_{-1}}\xi_{m},X_f\xi_{m})\cr
&=z^{-2km}(\xi_{m},e^{zL_1}X_f\xi_{m}),\cr}$$
which as before is non--zero. It follows that in all cases there is a non--zero 
primary field of type $m^2$ between $L(1,k^2)$ and $L(1,i^2)$ if $V_i\le V_m\otimes V_k$.
 
\vskip .1in
\noindent \bf 9. Proof of the character formula for $L(1,j^2)$ using the Jantzen filtration. \rm
Let $A(x)=\sum_{i\ge 0} A_i x^i$ be an analytic family of non--negative self--adjoint matrices defined 
for $x$ real and small on $V={\Bbb R}^n$ or ${\Bbb C}^n$. We define a filtration $V^{(i)}$ ($i\ge 0$) by
$V^{(0)}=V$ and for $m\ge 1$ 
$$V^{(m)}=\bigcap_{i=0}^{m-1}{\rm ker}\, A_i.$$
In the examples $A(x)$ will be invertible for $x\ne 0$ in a neighbourhood of $0$. 
It follows that $V^{(m)}=(0)$ for $m$ sufficiently large. We call $(V^{(m)})$ 
the {\it Jantzen filtration} associated with $A(x)$. 
\vskip .1in
\noindent \bf Lemma. \it The order of $0$ as a root of $\det A(x)$ equals $k=\sum_{i\ge 1} {\rm dim}\, V^{(i)}$. 
\vskip .05in
\noindent \bf Proof. \rm We set 
$U_0=V^{(0)}=V$, $U_1=V^{(n_1)}$ the first $k>0$ with $V^{(k)}\ne V^{(0)}$, 
$U_2=V^{(n_2)}$ the first $k>n_1$ with $V^{(k)}\ne V^{(n_1)}$, and so on. We set $n_0=0$. 
Let $v_1,\dots,v_{m_1}$ be an orthonormal basis of $U_0\ominus U_1$, $v_{m_1+1}, \dots,v_{m_2}$ 
an orthonormal basis of $U_1\ominus U_2$ and so on. Thus $(A(x)v_i,v_j)=x^{n_s}(A_{n_s}v_i,v_j) +$ higher powers of $x$ for
$i>m_s$. Moreover the first term vanishes if $j>m_{s+1}$. On the other hand the matrix 
$(A_{n_s}v_i,v_j)_{m_s<i,j\le m_{s+1}}$ is invertible. Thus we can divide rows $m_s+1,\dots m_{s+1}$ by $x^{n_s}$ 
and then set $x=0$. The resulting matrix is block triangular with invertible blocks on the diagonal, so is 
itself invertible. Hence the order of $0$ as a root of $\det A(x)$ equals 
$$k=\sum n_s\cdot ({\rm dim}\, U_{s-1} -{\rm dim}\, U_s)= \sum_{i\ge 1} {\rm dim}\, V_i,$$
as claimed.
\vskip .1in
Note that multiplying a bilinear form by (or matrix $A(x)$) by a factor $c_0+c_1x +c_2x^2+\cdots$ with $c_0\ne 0$ does not change the Jantzen filtration.  In this case we will say that the forms are {\it equivalent}.
\vskip .1in
We shall need an additional simple functoriality property of the Janzten filtration:
\vskip .05in
\noindent \bf Lemma. \it Let $V$, $W$ be finite--dimensional inner product spaces 
and for $x\in (-\varepsilon,\varepsilon)$. 
let $A(x)$, $B(x)$ be self--adjoint non--negative matrices depending analytically on $x$.
Let $X(x)\in {\rm Hom}(V,W)$ and $Y(x)\in {\rm Hom}(W,V)$ be analytic functions of $x$ such that
$(X(x)A(x)v,w)=(B(x)Y(x)v,w)$ for all $v\in V, w\in W$. Then $X(0)$ and $Y(0)^*$ induce maps
of the corresponding Jantzen filtrations which are dual to each other on $V^{(i)}/V^{(i+1)}$ abd
$W^{(i)}/W^{(i+1)}$.
\vskip .05in
\noindent \bf Proof. \rm The condition implies that $X(x)A(x)=B(x)Y(x)$. 
If $A(x)=\sum A_n x^n$, $B(x)=\sum B_nx^n$,
$X(x)=\sum X_n x^n$ and $Y(x)=\sum Y_n x^n$, then 
$$\sum_{i+j=n}X_i B_j=\sum_{i+j=n} A_i Y_j.\eqno{(*)}$$
Thus $X_0A_0=B_0Y_0$ on $V$. The equation $(*)$ implies that $X_0A_i=B_iY_0$ on 
$V^{(i)}={\rm ker} A_0\cap \cdots \cap {\rm ker} A_{i-1}$. But then $X_0$ and $Y_0^*$ induce dual maps between
$V^{(i)}/V^{(i+1)}$ and $W^{(i)}/W^{(i+1)}$, as claimed.  
\vskip .1in

Now consider a representation $M(1,j^2)$. The point $(c,h)=(1,j^2)$ lies on the curve $\varphi_{p,1}(c,h)=0$ where
$p=2j+1$. The curve is parametrized by a variable $t$ via $c(t)=13 - 6t - 6t^{-1}$ and $h(t)=(j^2+j)t -j$, 
the point $(1,j^2)$ corresponding to the value $t=1$. If $t\ne 0,\infty$, the Verma module $M_t=M(c(t),h(t))$ is defined. 
It is a direct sum of finite energy spaces $M_t(n)$ on which $L_0$ acts as multiplication by $n+h(t)$. The
space $M_t$ has a canonical invariant bilinear form $B_t(v,w)$ defined on it. Invariance implies that the subspaces
$M_t(n)$ are orthogonal. If we take as basis of $M_t(n)$ elements $v_i=L_{-i_r}\cdots L_{-i_1}v_t$ 
with $i_1\le i_2\le \cdots$ and $\sum s\cdot i_s=n$, then all these spaces can be identified with the same space 
$V={\Bbb R}^{{\cal P}(n)}$. Thus $A^{(n)}_{ij}(t)=B_t(v_i,v_j)$ is a symmetric matrix, the entries of which are polynomials in 
$t$ and $t^{-1}$ with real coefficients. Since the Shapovalov form is invariant under the Virasoro algebra, it follows from the second lemma above that of $M=M(1,j^2)$ then, for $t=1+x$ 
with $x$ small, the Jantzen filtration defines a filtration of $M$ by $\v$--submodules $M^{(i)}$. In particular we can calculate
$\sum_{i\ge 1} {\rm dim}\, M^{(i)}(n)$. 
It is exactly the degree of $t=1$ as a root of the Kac 
determinant at level $n$ which is known explicitly. 
\vskip .1in
We now use the Janzten filtration to compute the character of $M(1,j^2)$. From the coset construction, we know that
$M(1,j^2)$ contains a Verma submodule $M(1,(j+1)^2)$; this is also 
immediate from the Kac determinant formula for $c=1$:
$${\rm det}_N(1,h)=\prod_{p,q;\, pq\le N} (h-{(p-q)^2\over 4})^{\P(N-pq)}.$$ 
Continuing in this way, it has a decreasing chain of Verma
submodules $M(1,(j+k)^2)$ with $k\ge 1$. 
If we fix $h=j^2$ and set $c=1+x$, then as a function of $x$ the determinant is proportional to
$$\eqalign{\Psi_N(x)&=\prod_{1\le r\le N}(j^2+{r^2-1\over 24}\cdot x)\cdot\prod_{1\le s\le r\le N} [(j^2 -{(r-s)^2\over 4})^2 \cr
&\,\,\,\,\,\,\,\,\,\,+((r^2+s^2 -2)j^2 +{1\over 2} (rs+1)(r-s)^2)\cdot{x\over 24} + 
(r^2-1)(s^2-1)({x\over 24})^2].\,\,\,\,\,\,\,\,\,\,\,\,\,\,\,\,\,\,\,\,\,\,\,\,\,\,\,\,
\,\,\,\,\,\,\,\,\,\,\,\,\,\,\,\,\,\,\,\,\,\,\,\,\,\,\,\,\,\,\,\,\,\, (1)\cr}$$  
We take $x$ as the parameter of the Jantzen filtration. 
Since the operators $L_n$ depend polynomially on the parameters $t$ and $t^{-1}$, and hence analytically on $x$,
it is evident that when $h=j^2$, they preserve the filtration $M^{(i)}$ of $M=M(1,j^2)$. 
If $j\ge 0$, then from (1) we get
$$X_{j^2}(q)\equiv\sum_{i\ge 1} {\rm ch}\, M^{(i)}= \sum_{N\ge 0} a(N)q^{N+j^2}$$
where
$$a(N)=\sum_{1\le r\le N,\, } \P(N-r(r+2j)).$$
Thus we obtain
$$X_{j^2}(q)=\sum_{i\ge 1} {\rm ch}\,M^{(i)}= \varphi(q)\cdot \left(\sum_{r\ge 1} q^{(r+j)^2} \right).\eqno{(2)}$$
Let $v_r$ be the singular vector in $M(1,j^2)$ of energy $(j+r)^2$. Let  $n_r>0$ 
be maximal subject to $v_r\in M^{(n_r)}$. Thus
the Verma module generated by $v_r$ lies in $M^{(n_r)}$ and $\beta(w_1,w_2)=(w_1,w_2)_x/x^{n_r}|_{x=0}$ restricts to a 
multiple of the Shapovalov form. 
In particular $n_r$ is the order of the zero of $(v_r,v_r)_x$. (Note that for $x>0$ 
the Shapovalov form is positive definite, so that $(v_r,v_r)_x>0$ for $x>0$.)  
Since the Shapovalov form vanishes on any submodule, we have $\beta(v_{r+1},v_{r+1})=0$ 
and hence $n_r<n_{r+1}$. Moreover $M(1,(j+r)^2)\subseteq M^{(i)}$ 
for $n_{r-1}<i\le n_r$.  Hence
$$X_{j^2}(q)\ge \sum_{N\ge 1} {\rm ch}\, M^{(i)}\ge \sum_{r\ge 1} (n_r-n_{r-1}) {\rm ch}\, M(1,(j+r)^2)
=\varphi(q)\cdot \sum_{r\ge 1} (n_r-n_{r-1}) q^{(j+r)^2},\eqno{(3)}$$
where the inequality is between coefficients of $q^k$ and we set $n_0=0$.
Comparing (3) with (1) and (2), it follows that $n_r=r$ 
for $r\ge 1$ and that equality holds in (3). Hence $M^{(1)}=M(1,(j+1)^2)$ and
$${\rm ch}\, L(1,j^2) =(q^{j^2}-q^{(j+1)^2})\varphi(q),$$
as required. 
\vskip .1in
\noindent \bf 10. Proof of the character formula for the discrete series $0<c<1$ using the Jantzen filtration. \rm 
The proof of the character formula for the discrete series
$$c=1-{6\over m(m+1)},\,\, h=h_{r,s}={(r(m+1)-sm)^2-1\over 4 m(m+1)},\,(1\le s\le r\le m-1)$$
is proved similarly using $x=h-h_{r,s}$ as a parameter. Let
$$\varphi_{r,r}(c,h)=h+{1\over 24}(r^2-1)(c-1)$$
and for $r,s\ge 1$ with $r\ne s$,
$$\varphi_{r,s}(c,h)= \left(h-{(r-s)^2\over 4}\right)^2 +{h\over 24}(r^2+s^2-2)(c-1)
+{1\over 24^2} (r^2-1)(s^2-1)(c-1)^2 +{1\over 48}(c-1)(r-s)^2(rs+1).$$
We shall need the following result, implicit in {\bf [8]}. We give two proofs, the  
second of which is an expanded and slightly corrected version of a result from {\bf [13]}.
\vskip .1in
\noindent \bf Proposition (Fuchs algorithm). \it If $\varphi_{r,s}(c,h)=0$, then $M(c,h)$ has a non--zero 
singular vector at energy level $h+rs$.
\vskip .05in
\noindent \bf Proof~1. \rm A singular vector at energy level $h+rs$  has the form $Pv_{h}$ where 
$$P_d=\sum_{k =0}^{rs} q_k L_{-1}^k.$$
Each term $q_k$ is a sum of monomials $L_{-k}^{n_k} \cdots L_{-2}^{n_2}$ where $2n_2+3n_3+\cdots=d-k$. 
By the uniqueness result of Fuchs, we may assume $q_{d}=1$, where $d=rs$. Set $n_1=d-k$ and 
order these momomials lexigraphically on 
$(n_2,n_3,\dots)$. Thus the first term is $L_{-1}^d$.

Consider a term $\cdots L_{-p-1}^{n_{p+1}}L_{-p}^{n_p}L_{-1}^{d-n_1}v$ with $n_p>0$ 
and suppose that we have already determined the coefficients of all 
previous monomials in the lexicographic
order to be polynomials in $c$ and $h$. 
Let $a$ be the coefficient of $Bv=\cdots L_{-p-1}^{n_{p+1}}L_{-p}^{n_p}L_{-1}^{d-n_1}v$ in 
the singular vector $w$.
Let $k=d-n_1$. Then since fot a singular vector we would require $L_{p-1}w=0$, in this case we just
require that the coefficient of 
$$\cdots L_{-p-1}^{n_{p+1}}L_{-p}^{n_p-1}L_{-1}^{k+1}v\eqno{(*)}$$
in $L_{p-1}w$ must be zero.
We claim that this coefficient equals $n_p(2p-1)a$ plus a sum of lower order coefficients times 
polynomials in $c$ and $h$. In particular $a$ is again a polynomial in $c$ and $h$.

If $A$ is an element of the enveloping algebra of the Virasoro algebra we shall write $Av \sim 0$ 
if $Av$ is a combination of terms lower than $Bv$ in the lexicographic order, 
with coefficients polynomials in $c$ and $h$.
Now we have
$$L_{p-1}w=\sum_{j\ge 0} L_{p-1}q_j L_{-1}^{j}v=\sum_{k\ge 0}[L_{p-1},q_j]L_{-1}^jv+q_jL_{p-1}L_{-1}^jv.$$
We look for terms ending with $L_{-p}^{n-1}L_{-1}^{k+1}$ in this expression. We first note that it 
follows by induction on $k\ge 1$ that if $j\ge 1$ then
$L_jL_{-1}^kv$ lies in ${\rm lin}\, \{AL_{-1}^iv:i<k,\,\, A\in {\cal U}_2\}$. Indeed
$$L_{j}L_{-1} L_{-1}^{k-1} v=(j+1)L_{j-1}L_{-1}^{k-1}v +L_{-1}L_jL_{-1}^{k-1}v,$$
which has the same form since $L_{-1}{\cal U}_2={\cal U}_2 L_1$ and $L_{-1}^{k-1}v$ is an eigenvector of $L_0$.
If $k\ne c$, then the induction hypothesis forces $q_k\sim D$ with $D\in {\cal U}_p$. 
So either $k> c$, in which case any monomial
$q_k$ could have a non--zero contribution from a monomial $L_{-s}^{m_s} \cdots L_{-p}^{m_p}$ with $m_p>n_p$. 
Clearly taking the Lie bracket with $L_{p-1}$ can increase the exponent of $L_{-1}$ by at most one
while diminishing the exponent 
of $L_{-p}$ by at most one. So monomials ending with $L_{-p}^{n-1} L_{-1}^{c+1}$ 
have coefficients that are polynomials in $c$ and $h$.
If $k<c$, there is no way to increase the power of $L_{-1}^k$ to $L_{-1}^{c+1}$ 
by taking the Lie bracket with $L_{p-1}$. For $k=c$ and the terms $[L_{p-1},AL_{-p}^{n}] $ with $A$ a 
monomial in ${\cal U}_{p+1}$, we have
$$[L_{p-1},AL_{-p}^{n}]=[L_{p-1},A]L_{-p}^{n}+A[L_{p-1},L_{-p}^{n_p}].$$
The first term lies in ${\cal U}_2$ while for the second
$$[L_{p-1},L_{-p}^{n}]v=(2p-1) \sum_{a+b=n-1} L_{-p}^a L_{-1} L_{-p}^bv=(2p-1)n L_{-p}^{n-1}L_{-1} + B,$$
where $B\in {\cal U}_2$. There could be several terms in $q_c$ with $n_p=n$; but in this particular case 
on bracketing with $L_{p-1}$ 
and taking the term ending with $L_{-p}^{n-1}L_{-1}^{c+1}$ we obtain $A  L_{-p}^{n-1}L_{-1}^{c+1}$. It follows that 
all these terms are linearly independent. But then there can be no cancellation and 
the coefficient of $ L_{-r}^{n_r} \cdots L_{-p}^{n_p-1}L_{-1}^{c+1} v$ must be $n_p(2p-1)a$ 
plus a sum of lower order coefficients times 
polynomials in $c$ and $h$, as claimed.   

Now given this vector $w=\sum  a_\alpha(c,h) L_\alpha L^{d-|\alpha|}$ where the coefficients $a_\alpha(c,h)$
are polynomials in $c$ and $h$, the condition $L_1w=0=L_2w$ gives a series of polynomials $b_i(c,h)$ which must vanish
for $w$ to be a singular vector. By the coset construction, we know that $L(c_m, h_{r,s}(q))$ has a singular vector
at energy $h_{r,s}(m)+rs$ if $1\le s\le r\le m-1$. Thus $b_i(c_m,h_{r,s}(m))=0$. Thus $b_i(c,h)$ vanishes at a point of
accumulation on the real curve $\varphi_{r,s}(c,h)=0$. It follows that $b_i(c,h)=0$ whenever $\varphi_{r,s}(c,h)=0$. Thus
$w$ defines a non--zero singular vector if $\varphi_{r,s}(c,h)=0$, as required.
\vskip .05in
\noindent \bf Proof~2. \rm A singular vector at energy level $h+rs$  has the form $Pv_{h}$ where 
$$P_d=\sum_{k =0}^{rs} q_k L_{-1}^k.$$
Each term $q_k$ is a sum of monomials $L_{-k}^{n_k} \cdots L_{-2}^{n_2}$ where $2n_2+3n_3+\cdots=d-k$. 
By the uniqueness result of Fuchs, we may assume $q_{d}=1$, 
where $d=rs$.  Set $n_1=d-k$ and order these momomials lexigraphically on 
$(n_1,n_2,n_3,\dots)$. Thus the first term is $L_{-1}^d$.

Consider a term $\cdots L_{-p-1}^{n_{p+1}}L_{-p}^{n_p}L_{-1}^{d-n_1}v$ with $n_p>0$ 
and suppose that we have already determined the coefficients of all previous monomials in the lexicographic
order to be polynomials in $c$ and $h$. 
Let $a$ be the coefficient of $Bv=\cdots L_{-p-1}^{n_{p+1}}L_{-p}^{n_p}L_{-1}^{d-n_1}v$ in the singular vector $w$.
Let $k=d-n_1$. Then since fot a singular vector we would require $L_{p-1}w=0$, in this case we just
require that the coefficient of 
$$\cdots L_{-p-1}^{n_{p+1}}L_{-p}^{n_p-1}L_{-1}^{k+1}v\eqno{(*)}$$
in $L_{p-1}w$ must be zero.
We claim that this coefficient equals $n_p(2p-1)a$ plus a sum of lower order coefficients times 
polynomials in $c$ and $h$. In particular $a$ is again a polynomial in $c$ and $h$.

If $A$ is an element of the enveloping algebra of the Virasoro algebra we shall write $Av \sim 0$ 
if $Av$ is a combination of terms lower than $Bv$ in the lexicographic order, 
with coefficients polynomials in $c$ and $h$. 

We have
$$L_{p-1}w=\sum_{k=0}^e L_{p-1}q_k L_{-1}^{k}v=\sum_{k=0}^e[L_{p-1},q_k]L_{-1}^kv+q_kL_{p-1}L_{-1}^kv.$$
We look for terms ending with $L_{-p}^{n-1}L_{-1}^{e+1}$ in this expression. We first note that it 
follows by induction on $k$ that if $j\ge 1$ then
$L_jL_{-1}^kv$ lies in ${\rm lin}\, \{AL_{-1}^iv:i<k,\,\, A\in {\cal U}_2\}$. Indeed
$$L_{j}L_{-1} L_{-1}^{k-1} v=(j+1)L_{j-1}L_{-1}^{k-1}v +L_{-1}L_jL_{-1}^{k-1}v,$$
which has the same form since $L_{-1}{\cal U}_2={\cal U}_2 L_1$ and $L_{-1}^{k-1}v$ is an eigenvector of $L_0$.

By the inductive hypothesis, all terms with $k>e$ are have uniquely determined polynomial coefficients.
A monomial with $k\le e$ can be written
$u=B\cdot L_{-p}^{n} A \cdot L_{-1}^kv$ 
with $B$ a monomial in $L_{-j}$'s for $j>p$ and $A$ a monomial in the $L_{-j}$'s with
$2\le < p$. But then
$$L_{p-1}u=[L_{p-1},B] \cdot L_{-p}^nA L_{-1}^kv 
+ B\cdot (L_{p-1},L_{-p}^n]\cdot AL_{-1}^kv)+ BL_{-p}^n L_{p-1}AL_{-1}^kv.\eqno{(1)}$$
Since $[L_{p-1},B]$ lies in ${\cal U}_2$, the first term in this expression has no terms ending in $L_{-1}^{e+1}$.
For the second note that
$$[L_{p-1},L_{-p}^{n}]=(2p-1) \sum_{a+b=n-1} L_{-p}^a L_{-1} L_{-p}^b=(2p-1)n L_{-p}^{n-1}L_{-1} + D,$$
where $D\in {\cal U}_2$. Clearly, if $k<e$, then neither $BDL_{-1}^k v$ nor $B L_{-1} AL_{-1}^k v$ contain terms 
ending in $L_{-1}^{e+1}v$. Finally for the last term, note that if $X$ is a monomial in $L_{-j}$ for $1\le j \le t$,
then for $1\le s \le t$ we have $L_s Xv=\sum a_i X_iv$, where $X_i$ has the same form as $X$ with the exponent of
$L_{-1}$ increased by at most one. Thus no terms ending in $L_{-p}^{n_p-1} L_{-1}^{e+1}$ arise this way if $k<e$. 

When $k=e$, then by the inductive hypothesis, we only need to consider terms of the form 
$u=B\cdot L_{-p}^nL_{-1}^ev$ with either $n\ge n_p>0$ or zero. If $n=0$, then
$$L_{p-1}u=L_{p-1}BL_{-1}^ev=[L_{p-1},B]L_{-1}^ev + B L_{p-1} L_{-1}^e v.$$
The first term contains no terms ending in $L_{-1}^{e+1}v$. Nor does the second, since if $kj0$, $L_{j}L_{-1}^ev$ lies in
${\rm lin}\, \{L_{-1}^iv:0\le i<e\}$. If $n>n_p$, the none of the terms in (1) give rise to a monomial ending in
$L_{-1}^{n_p-1} L_{-1}^{e+1}v $. 

Finally there could be several terms $BL_{-p}^{n}L_{-1}^e$ in $q_e$ with $n=n_p$. However taking the term ending with 
$L_{-p}^{n-1}L_{-1}^{e+1}v$ in (1), we get $B  L_{-p}^{n-1}L_{-1}^{e+1}v$. Thus 
all these terms are distinct. But then there can be no  cancellation and
and thus
the coefficient of $ L_{-r}^{n_r} \cdots L_{-p}^{n_p-1}L_{-1}^{c+1} v$ must be $n_p(2p-1)a$ 
plus a sum of lower order coefficients times polynomials in $c$ and $h$.

Now given this vector $w=\sum  a_\alpha(c,h) L_\alpha L_{-1}^{d-|\alpha|}$ where the coefficients $a_\alpha(c,h)$
are polynomials in $c$ and $h$, the condition $L_1w=0=L_2w$ gives a series of polynomials $b_i(c,h)$ which must vanish
for $w$ to be a singular vector. By the coset construction, we know that $L(c_m, h_{r,s}(q))$ has a singular vector
at energy $h_{r,s}(m)+pq$ if $1\le s\le r\le m-1$. Thus $b_i(c_m,h_{r,s}(m))=0$. Thus $b_i(c,h)$ vanishes at a point of
accumulation on the real curve $\varphi_{r,s}(c,h)=0$. It follows that $b_i(c,h)=0$ whenever $\varphi_{r,s}(c,h)=0$. Thus
$w$ defines a non--zero singular vector if $\varphi_{r,s}(c,h)=0$, as required.

\vskip .1in
\bf \noindent Corollary. \it If $r,s\ge 1$ and $c=1-6/m(m+1)$, then $M(c,h_{r,s})$ has a singular vector of energy
$h=h_{r,s}+rs$. 
\vskip .1in
\noindent \bf Proposition. \it Let $c=1-6/m(m+1)$,  $1\le s\le r\le m-1$, $r^\prime=m-r$, $s^\prime=m+1-q$.  Then $M(c,h_{r,s})$ has singular vectors $a_1,b_1,a_2,b_2,\cdots $ where
for $k\ge 1$, $a_k$ has energy $\alpha_k=h_{r^\prime,k(m+1)+s^\prime}=h_{r,s-k(m+1)}$ and $b_k$ has energy
$\beta_k=h_{r,s+k(m+1)}=h_{r^\prime, s^\prime -k(m+1)}$. The energies of the vectors are 
strictly increasing in the above order.
Moereover if $A_k$ and $B_k$ are the Verma modules generated by $a_k$ and $b_k$, then
$A_{k+1}+ B_{k+1}\subseteq A_k\cap B_k$. 
\vskip .05in
\noindent \bf Proof. \rm Since
$$\alpha_{k+1}=\alpha_k+ r^\prime(k(m+1)+s^\prime),\,\, \beta_{k+1}=\alpha_k+r(s-k(m+1))$$
and 
$$\alpha_{k+1}=\beta_k+ r(s+k(m+1)),\,\, \beta_{k+1}=\beta_k+ r^\prime (s^\prime -k(m+1)),$$
existence follows by successive applications of the previous corollary. The uniqueness of singular vectors as a 
particular energy level in a Verma module implies the result on the set theoretic inclusion.
\vskip .1in
For fixed $m$ and $1\le s\le r\le m-1$, take $x=h-h_{r,s,}(m)$.  
As in the previous section let 
$m_k>0$ be the order of the zero of $(a_k,a_k)_x$ and $n_k$ that of $(b_k,b_k)_x$.
For this to make sense, we have to check that, if $\xi=a_k$ or $b_k$, then $(\xi,\xi)_x$ does not vanish identically. But m $(\xi,\xi)_x>0$ for $x$ sufficiently large, 
because the Kac determinant formula and the unitarity of 
Verma modules for $c\ge 1$ together imply that,
at a fixed energy level, the Shapovalov form is positive definite for $c>0$ and $h$ sufficiently large. The same argument as used in the 
Section~9 shows that $m_k<m_{k+1},n_{k+1}$ and $n_k<m_{k+1},n_{k+1}$ since $A_{k+1}+B_{k+1}\subseteq A_k \cap B_k$. 
\vskip .1in
\noindent \bf Proposition. \it $A_1+B_1$ is the maximal submodule of $M=M(c,h_{r,s})$ 
and  $A_{k+1}+ B_{k+1}=A_k\cap B_k$ for $k\ge 1$. Moreover $m_k=n_k=1$ and
$$\sum_{i\ge 1} {\rm ch}\, M^{(i)} = {\rm ch}(A_1+B_1)+\sum {\rm ch A_{2k}}+{\rm ch}\, B_{2k} 
= \sum_{k\ge 0}{\rm ch} \,A_{2k+1}+{\rm ch}\, B_{2k+1}.$$
\vskip .05in
\noindent \bf Proof. \rm From the Kac determinant formula 
$h_{r,s}=h_{r_1,s_1}$ for $(r,s)=(r_1,s_1)+a(m,m+1)$ with
$a\ge 0$ or $(r,s)=-(r_1,s_1)+b(m,m+1)$ with $b\ge 1$. Hence
$$\sum {\rm ch} \, M^{(i)}=\varphi(q)\cdot\sum_{a\in {\Bbb Z}} 
q^{(r+am)(s+a(m+1))}.\eqno{(1)}$$
On the other hand if $M_k=\max (m_k,n_k)$ 
and $N_k=\min (m_k,n_k)$, then $N_k\le M_k < N_{k+1}\le M_{k+1}$.   
On the other hand
$$\sum {\rm ch} \, M^{(i)}\ge \sum_{k\ge 1}(N_{k+1}-M_k) 
{\rm ch}\, (A_k+ B_k)\ge \sum_{k\ge 1}{\rm ch}\, (A_k+ B_k).
\eqno{(2)}$$
Thus
$$\sum {\rm ch} \, M^{(i)}\ge \sum_{k\ge 1}{\rm ch}\, (A_k+ B_k).$$
We wish to prove the same inequality with $A_k\cap B_k$ replacing $A_{k+1}+B_{k+1}$ for all 
$k$. Indeed suppose that $R_k$ is chosen maximal such that $A_k\cap B_k \subseteq M^{(R_k)}$. Then
$M_k+1\le R_k< N_k$. Hence $R_k<R_{k+1}$ and so 
$$\sum {\rm ch} \, M^{(i)}\ge {\rm ch}\, (A_1+B_1) +\sum_{k\ge 1}(R_{k}-R_{k-1}) {\rm ch}\, A_k\cap B_k \ge 
{\rm ch}\,(A_1+B_1)+ \sum_{k\ge 1}{\rm ch}\, A_k\cap B_k.$$
Now we have a short exact sequence
$$0\rightarrow A_k \cap B_k \rightarrow A_k \oplus B_k \rightarrow A_k +B_k\rightarrow 0.$$
Hence
$${\rm ch}\, (A_1+B_1) + {\rm ch}\,A_1\cap B_1 ={\rm ch}\, A_1 + {\rm ch}\, B_1,$$ 
and for $k\ge 1$ 
$${\rm ch}\, A_{2k}\cap B_{2k} +{\rm ch}\, A_{2k+1} \cap B_{2k+1}
\ge 
{\rm ch}\, (A_{2k+1}+ B_{2k+1}) +{\rm ch}\, A_{2k+1} \cap B_{2k+1}={\rm ch}\, A_{2k+1} +{\rm ch}\, B_{2k+1}.$$
Hence
$$\sum {\rm ch} \, M^{(i)}\ge {\rm ch}\, (A_1+B_1) +\sum_{k\ge 1}{\rm ch}\, A_k\cap B_k\ge 
 \sum_{k\ge 0} {\rm ch}\, A_{2k+1} +{\rm ch}\, B_{2k+1}.\eqno{(3)}$$
On the other hand the right hand side equals
$$\varphi(q)\cdot \sum_{a\in {\Bbb Z}} 
q^{(r+am)(s+a(m+1))},$$
so it  coincides with the left hand side. So we have equality in (3), so that
$$A_{2k}\cap B_{2k}= A_{2k+1}+B_{2k+1}$$
for all $k$. Now instead of $M$ we take $M^\prime=A_1$ and 
set $A_k^\prime=B_{k+1}$, $B_k^\prime=A_{k+1}$ for $k\ge 1$. Repeating
the argument above we get
$$\sum_{i\ge 1} {\rm ch}\,(M^\prime)^{(i)}\ge \sum_{k\ge 1}  {\rm ch}\, A_k\cap B_k
\ge \sum_{k\ge 1} {\rm ch}\, A_{2k} +{\rm ch}\, B_{2k}.\eqno{(4)}$$
By (1), the left hand side of (4) equals
$$\varphi(q)\cdot\sum_{a\in {\Bbb Z}} q^{(r_1+am)(s_1+a(m+1)},$$
where $r_1=r$ and $s_1=-s$. But this equals the last term on the right hand side. Hence
equality holds in (4), so that
$$A_{2k-1}\cap B_{2k-1}= A_{2k}+B_{2k}$$
for $k\ge 1$.  
Returning to equation (2), we have equality between the furthest terms and 
hence $N_{k+1}-M_k=1$ for all $k$.  
Evidently $m_1< m_2< \cdots$ and $n_1< n_2<\cdots$. Now we have $M^{(1)}\supseteq A_1+B_1$, so $m_1=1=n_1$. Hence
$M_1=1$. So $N_2=\min(m_2,n_2)=2$. Note that if $N_k<M_k$ then we would have an extra contribution of 
${\rm ch} \, A_k$ or  ${\rm ch} \, B_k$ 
on the right hand side of (2). Since equality holds, it follows that $N_k=M_k$, so that
$n_k=m_k$. Since $m_1=1=n_1$, we finally get $m_k=k=n_k$.  

\vskip .1in

\noindent \bf Corollary. \it ${\rm ch}\, L(c,h_{r,s})=\varphi(q)\cdot \sum_{k\in {\Bbb Z}} (-1)^k q^{h_{r,s+k(m+1)}}$.    
\vskip .05in
\noindent \bf Proof. \rm We have
$$\eqalign{{\rm ch}\, L(c,h_{pq})&= {\rm ch}\, M(c,h_{r,s})/(A_1+B_1)\cr
&={\rm ch}\, M(c,h_{r,s}) - {\rm ch}\, (A_1+B_1)\cr
&={\rm ch}\, M(c,h_{r,s})+\sum_{k\ge 0}{\rm ch} \,A_{2k+1}+{\rm ch}\, B_{2k+1}
  - \sum_{k\ge 0}(\,{\rm ch} \,A_{2k+1}+{\rm ch}\, B_{2k+1})\cr
&=\varphi(q)\cdot \sum_{k\in {\Bbb Z}} (-1)^k q^{h_{r,s+k(m+1)}}.\cr}$$    
\vskip .1in

\noindent \bf APPENDIX~A: Alternative proofs of Fubini--Veneziano relations. \rm 
\vskip .1in
\noindent \bf First indirect proof. \rm Recall that
$\phi(n):{\cal F}\rightarrow {\cal F}$ is a primary field for the system of operators 
$(a_n)$, $U$ and $L_0$ if
$$[a_n,\phi(k)]=m\phi(k+n),\,\,[L_0,\phi(n)]=-(n+\mu)\phi(n),\,\, U\phi(n)U^*=\phi(n+m),$$
for some $m\in{\Bbb Z}$ and $\mu\in {\Bbb R}$. We now fix $k$ and consider
$$\psi_k(n)=[L_k,\phi(n-k)]+(n-k)\phi(n).$$
It is easy to check that $\psi$ satisfies the same conditions as $\phi$ with the same choice of $\mu$. So by uniqueness,
$\psi_k(n)=c(k)\phi(n)$ for some constant $c(k)$:
$$[L_k,\phi(n-k)]+(n-k)\phi(n)=c(k)\phi(n).\eqno{(1)}$$
The Jacobi relation, the relation $[L_a,[L_b,T]]-[L_b,[L_a,T]]=(a-b)[L_{a+b},T]$ and the non--vanishing of $\phi(0)$ 
imply that 
$$(a-b)c(a+b) =ac(a)-bc(b).$$
But this functional equation implies that $c$ is an affine function $c(a)=\alpha a + \beta $. Indeed the $c$ satisfying 
this functional equation form a vector space. So subtracting an affine function from $c$, we may 
assume that $c(0)=0=c(1)$ and must then show that $c\equiv 0$. Taking $a=-b$, we see that $c(a)=-c(-a)$ for all $a$. 
But for $a>1$, $(a-1)c(a+1)=ac(a)$. Hence $c(a)=0$ for all $a>0$ and hence for all $a$.  
So from (1) we have
$$[L_k,\phi(n)]=(-n -\alpha k -\beta)\phi(n+k).\eqno{(2)}$$
Now on the one hand we have $\Phi_m(z)=\sum \phi(n)z^{-n-\delta}$ while on the other, regardless of the labelling of 
the $\phi(n)$'s, we have
$$\Phi_m(z)\Omega|_{z=0}=\Omega_m.$$
It follows that $\delta=0$ and $\phi(n)\Omega=0$ for $n>0$, with 
$\phi(0)\Omega=\Omega_m$. Since $L_0\Omega_m={m^2\over 2}\Omega_m$, we must have $\beta=-m^2/2$. From equation (2) 
we have that
$$[L_k,\Phi_m(z)]=z^{k+1} \Phi_m^\prime(z) + \Delta z^k(k+1)\Phi_m(z),$$
where $\Delta=1-\alpha$ and $\delta=\alpha-\beta-1$. But $\delta=0$. Hence $\alpha=1-m^2/2$ and $\Delta=m^2/2$. This 
proves the Fubini--Veneziano relations.

\vskip .1in
\noindent \bf Second proof by direct verification. \rm We start by noting 
the adjoint relation
$$(\Phi_m(z^{-1}))^*=z^{-m^2} \Phi_{-m}(z)$$
follows immediately because
$$\eqalign{\Phi_m(z^{-1})^*&=\exp(\sum_{n<0} {-z^n m a_n\over n}) \exp( \sum_{n>0} {-m z^n a_n\over n})z^{-ma_0}U^{m}\cr
&=z^{-ma_0} U^{m} \exp(\sum_{n<0} {-z^n m a_n\over n}) \exp( \sum_{n>0} {-m z^n a_n\over n})\cr
&=z^{-m^2}\Phi_{-m}(z).\cr}$$
We shall need the following generalisation of Lemma~B in Section~8. 
\vskip .1in 

\noindent \bf Lemma. \it Let $A$ be a formal power series in $z$ (or $z^{-1}$) with operator coefficients and let 
$D$ be an operator such that $C=[A,B]$ is a multiple of the identity operator, where $B=[D,A]$.
Then $[D,e^A]=(B+C/2)e^A$.
\vskip .05in
\noindent \bf Proof. \rm We have 
$$\eqalign{[D,A^N]&=\sum_{p+q=N-1} A^pBA^q\cr
                 &=\sum_{p+q=N-1} BA^{p+q} +pCA^{p+q-1}\cr
                 &=NBA^{N-1} + {N(N-1)\over 2} CA^{N-2}.\cr}$$
Hence
$$[D,e^A]=(B+C/2)e^A.$$
\vskip .1in
\noindent \bf Corollary. \it 
\noindent (a) $[L_k,E_+(z)]=(-\sum_{n>0} a_{n+k} z^{-n})E_+(z)$, 

and $[L_{-k},E_-(z)]=(-\sum_{n>0} a_{-n-k} z^{-n})E_-(z)$
if $k\ge 0$.

\noindent (b) $[L_{-k},E_+(z)]=(-\sum_{n>0} a_{n-k} z^{-n}+{m^2\over 2}(k-1)z^k)E_+(z)$

and $[L_k,E_-(z)]=(-\sum_{n>0} a_{-n-k} z^{n}-{m^2\over 2}(k+1)z^{-k})E_-(z)$
for $k>0$.  

\vskip .05in
\noindent \bf Proof. \rm These formulas are straightforward consequences of the lemma, 
setting $D=L_k$ and 
$$A=\sum_{\pm n>0} {m\over n} a_nz^{-n}.$$
For example to prove the first formula in (b), we have 
$$B=[D,A]=-\sum_{n>0} ma_{n-k}z^{-n}=-\sum_{i>-k} ma_i z^{-i+k}$$ 
so that
$$\eqalign{[A,B] &=-[\sum_{n>0} {m\over n} a_nz^{-n},\sum_{i>-k} ma_i z^{-i+k}]\cr
&=-\sum_{n=1}^{k-1} [a_n,a_{-n}]{m^2\over n} z^k\cr
&=-m^2(k-1)z^k.\cr}$$
\vskip .1in
\noindent \bf Proof of the Fubini--Veneziano relations. \rm We first check 
the commutation relations with $L_k$ for $k\ne 0$. For $k>0$, we have
$UL_{-k}U^*=L_{-k} + a_{-k}$, so that $U^mL_{-k}U^{-m}=L_{-k} +m a_{-k}$ and hence
$$[L_{-k},U^m]=-mU^ma_{-k}.$$
Thus
$$[L_-k,\Phi_m(z)]=U^m z^{-ma_0}(B_-E_-E_++ E_-B_+E_+),$$
where 
$$B_-=-m a_{-k}-\sum_{n> 0}m a_{-n-k} z^{-n},\,\,B_+=-\sum_{n>0} m a_{n-k} z^{-n}-{m\over 2}(z^k-1)/(z-1).$$
On the other hand
$$z^{-k+1}\Phi^\prime_m(z)=U^mz^{-ma_0} z^{-k}(-ma_{0})E_-E_+ + U^m z^{-ma_0} (C_-E_-E_+ +E_-C_+E_+),$$
where 
$$C_-=-\sum_{n>0} ma_{-n}z^{n-k},\,\,C_+=-\sum_{n>0} ma_n z^{-n-k}.$$
Thus 
$$B_-=C_--A$$
and 
$$B_+=C_+ + A
-{m^2\over 2}(k-1)z^{-k}I$$
where
$$A=\sum_{i=0}^{k-1} ma_i z^{-i-k}.$$
On the other hand
$$\eqalign{[A,E_-(z)]&= \sum_{i=0}^{k-1} m[a_i,E_-(z)]z^{i-k}\cr
& = \sum_{i=1}^{k-1} m^2z^{-k} E_-(z)\cr
&=m^2z^{-k}(-1)k E_-(z).\cr}$$
Hence
$$\eqalign{[L_{-k},\Phi_m(z)]-z^{-k+1}\Phi^\prime_m(z)
&=[-{m^2\over 2}(k-1) + m^2(k-1)]\Phi_m(z)\cr
&={m^2\over 2}(k-1) z^{-k}\Phi_m(z),\cr}$$
as required. The relation $[L_k,\Phi_m(z)]$ can be proved similarly or follows from this one by taking adjoints.
\vskip .1in
\noindent \bf APPENDIX B: Explicit construction of singular vectors 
and asymptotic formulas of Feigin--Fuchs. 
\rm In this appendix we give a slightly simplified account of the approach of 
{\bf [3], [7]} to singular vectors in $M(c(t),h_{r,1}(t))$, 
where $c(t)=13-6t-6t^{-1}$, $r=2j+1$ for $j$ a non--negative half integer and $h_{r,1}(t)\equiv h(t) =(j^2+j)t-j$.
The case of interest in the text has $t=1$, so that $c=1$ and $h=j^2$. 
The treatment of {\bf [8]}, however, uses the more general case, when the singular vector is a polynomial in $t$ 
and a knowledge of the constant term and leading coefficient is required. A proof of 
their formulas in a more general setting was given in {\bf [2]}, 
who commented that in the particular case
the method of {\bf [3]} for determining 
the singular vector could not be used to derive the result. We 
show on the contrary that the formulas are a trivial 
consequence of the algorithm of {\bf [3]}, which provided an 
alternative  approach to explicit formulas for the 
singular vectors of Benoit and St--Aubin {\bf [4]}.  
Let $E,F,H$ be a canonical basis of $\s$ satisfying 
$[E,F]=2H$, $[H,E]=E$ and $[H,F]=-F$ and let
$C=H^2 +(EF+FE)/2$, the Casimir element. Let $V=V_j$ be the irreducible 
representation of $\s$ of spin $j$ with 
(non--orthonormal) basis $v_{-j}, v_{-j+1}, \cdots, v_j$ satisifying
$Hv_k=kv_k$, $F^kv_j=v_{j-k}$ for $k>0$ and $F^{2j+1}v_j=0$. Let $W=\bigoplus_{m\in {\Bbb Z}} W(m)$ be a 
representation of the Virasoro algebra with $L_0w=(h+m)w$ for $w\in W(m)$ and ${\rm dim}\, W(m)<\infty$ 
for all $m$. We make operators $A$ on $V$ and $B$ on $W$ act on $V\otimes W$ as $A\otimes I$ and $I\otimes B$ respectively.
With this convention, we define
$$N=-F +\sum_{m\ge 0} (-tE)^m L_{-m-1},\,\, M=L_0 - H - tC, \,\, L= L_1 +E(t(H-1)+1),\,\, K=L_2 -tE^2(t(H-{3\over 2})+{7\over 4}).$$ 
\vskip .1in
\noindent\bf Lemma~A. \it $[N,M]=aN$, $[N,L]=bEN +c M$ and $[N,K]=d E^2N + e EM +f L$ where $a=-1$, $b=-3t$, $c=-2$, 
$d=-5t^2$, $e=-4t$ and $f=-3$. 
\vskip .05in
\noindent \bf Remark. \rm More generally, if we define ${\cal L}_{-1}=N$, ${\cal L}_0=M$ and 
$${\cal L}_k=L_k -(-E)^kt^{k-1}(t(H-{k+1\over 2})+{3k+1\over 4})$$ 
for $k\ge 1$ (so that ${\cal L}_1=L$ and ${\cal L}_2=K$), then for $p\ge 0$
$$[{\cal L}_{p},{\cal L}_{-1}] =\sum_{q\ge 0} (p+1+q)t^q E^q{\cal L}_{p-1-q}$$
and for $p,q\ge 0$
$$[{\cal L}_p,{\cal L}_q]=(p-q){\cal L}_{p+q}.$$
 
\vskip .05in
\noindent \bf Proof. \rm This is a straightforward verification using the commutation relations.
\vskip .1in
\noindent \bf Lemma~B. \it Let $\xi\in W$ satisfy $L_0\xi=h(t) \xi$. Then there is
a unique vector $w=\sum_{k=0}^{2j+1} v_{-j+k}\otimes \xi_{-j+k}$, $\xi_{-j}=\xi$ with $L_0\xi_k=(h(t) +k)\xi_k$ 
and $Nw=v_j\otimes \eta$.
In this case $\eta=P_j \xi_0$ with $P_j$ a uniquely determined homogenenous polynomial of total degree $r=2j+1$ in 
the universal enveloping algebra ${\cal U}_-$ of the Lie algebra $\v_-$ 
with basis $L_{-k}$ ($k>0$) with coefficients polynomials in $t$ of degree $\le 2j$. 
The constant coefficient is $L_1^r$ and the coefficient of $t^{2j }$ is $((2j)!)^2  L_{-r}$. More generally, 
the coefficient of $L_{-1}^r$ is $1$. 

\vskip .05in
\noindent \bf Remark. \rm We will check below that $P_j$ is non--zero by calculating its action a specific module. 
In fact it can also be seen directly on a 1--dimensional module for $\v_-$. For $z\in {\Bbb C}$ define that
algebra homomorphism $\psi_z:\v_-\rightarrow {\Bbb C}$ by $\psi_z(L_{-1})=z$ and $\psi_z(L_{-k})=0$ for $k>1$. 
Then $\psi_z(P_j)=z^{2j+1}$. Indeed we have to solve $N^\prime w=av_j$ with $w=\sum a_kv_k$, $a_{-j}=1$ and $N^\prime=-F+zI$.
The solution is $w=v_{-j} + z v_{-j+1} +z^2 v_{-j+2} + \cdots$ since we can check that $(zI -F)w=z^{2j+1}v_j$.
Hence $\psi_j(P_j)=z^{2j+1}$. (Alternatively this can be proved using the Lie algebra endomorphism of $\v_-$ defined by 
$\theta(L_{-1})=L_{-1}$ and $\theta(L_k)=0  $ for $k\le -2$.) 
\vskip .05in
\noindent \bf Proof. \rm  The vectors $\xi_k$ are defined inductively by $\xi_{-j}=\xi$ and for $k=-j,\dots, j$
$$\xi_{k+1}\otimes v_k=L_{-1}\xi_{k} \otimes v_k - L_{-2} \xi_{k-1} \otimes tE v_{k-1} 
+ L_{-3} t^2 E^2v_{k-2}+\cdots. \eqno{(*)}$$
This proves uniqueness. By induction $\xi_{-j+k}$ has the form $Q_k\xi$, with $Q_k$ a polynomial in ${\cal U}_-$ 
of total degree $k$. Performing this process for the Verma module $M(1,j^2)\cong {\cal U}_-$, gives a uniquely determined 
polynomial $P_j\in {\cal U}_-$. Since there is a natural homomorphism of $M(c,h)\rightarrow W$ which is compatible 
with the equations defined by $N$, it is clear that $\eta=P_j\xi$.

Since the Lie algebra with basis $L_{-k}$ ($k\ge 0$) is the semidirect product of the algebra 
with basis $L_{-k}$ ($k\ge 2$) and ${\Bbb C}L_{-1}$, it makes sense to talk about the coeffcients of powers of 
$L_{-1}$ in ${\cal U}_-$. Taking the homomorphism $\pi$ sending $L_{-k}$ to zero for 
$k\ge 2$ the recurrence relation becomes:
$$\pi(\xi_{k+1})\otimes v_k=L_{-1}\pi(\xi_{k}) \otimes v_k,$$
so that the coefficient of $L_{-1}^r$ is $1$. Alternatively using the evaluation map $\sigma$ that set $t=0$ 
in the recurrence relation gives
$$\sigma(\xi_{k+1})\otimes v_k=L_{-1}\sigma(\xi_{k}) \otimes v_k,$$
which yields the same result.
To get the leading cofficient in $t$, note that by the recurrence relation $\xi_{-j+k}$ is polynomial in $t$ of degree 
at most $k-1$ for $k\le r$. On the other hand by definition
$$\eta\otimes v_j=L_{-1}\xi_{j} \otimes v_j - L_{-2} \xi_{j-1} \otimes tE v_{j-1} 
+ L_{-3}\xi_{j-2}\otimes t^2 E^2v_{j-2}+ \cdots + L_{-2j-1}\xi_{-j}\otimes t^{2j} E^{2j}v_{-j}.$$
Thus the highest power of $t$ in $\eta$ is $t^{2j}$ and the coefficient is $c L_{-r}\xi$, where $cv_{j}=E^{2j}v_{-j}$.
On the other hand $C=H^2 -H +EF$ so that
$$Ev_{k}=EFv_{k+1}=(j^2+j -(k+1)^2 +k+1)v_{k+1}=(j^2+j-k^2 -k)v_{k+1}=(j-k)(j+k+1)v_{k+1}.$$
Thus $E^{2j}v_{-j}=((2j)!)^2v_j$ and hence $c=((2j)!)^2$.
 
\vskip .1in
\noindent \bf Lemma~C. \it In the Verma module $M(c(t),h(t))$ generated by the vector $\xi$,
there is a non--zero homogeneous polynomial $P_j$ of total degree $(2j+1)$ in 
the universal enveloping algebra ${\cal U}_-$ of the Lie algebra with basis $L_{-k}$ ($k>0$)  
such that $P_j\xi$ is a non--zero singular vector.
\vskip .05in
\noindent \bf Proof. \rm Since the map $A\mapsto A\xi$ gives an isomorphism between ${\cal U}_-$ and the Verma module, 
$P_j\ne 0$ forces $P_j\xi \ne 0$. We check that 
$$L_0\eta=(j+1)^2 \eta,\,\, L_1\eta=0=L_2\eta$$
Since $L_1$ and $L_2$ generate the Lie algebra with basis $L_{-k}$($k>0$), the second identities imply that $L_k\xi_0=0$ 
for all $k>0$, i.e. that $\eta$ is a singular vector.

Let $w$ be the vector in Lemma~B with $\xi_{-j}=\xi$, $Nw=v_j\otimes \eta$ and let $P=EN$. Thus $Pw=0$ and any solution 
$w^\prime \in V\otimes M_j^2$ of $Pw^\prime =0$ must satisfy $Nw^\prime= v_j\otimes \eta^\prime$ 
for some $\eta^\prime \in M_{j^2}$. On the other hand the uniqueness statement in Lemma~B, 
if $w^\prime=\sum v_i\otimes \xi_i^\prime$ satisfies $\xi_{-j}=0$ and $Pw^\prime=0$ then $w^\prime=0$.  
 
We now show that $w^\prime=Aw$ with $A=M,L,K$ satisfies these hypotheses. First we check that the coefficient of 
$v_{-j}$ in $Aw^\prime$ is zero. For $A=M$, we have
$$\eqalign{M(v_{-j}\otimes \xi_{-j})&=(L_0-H-tC)(v_{-j}\otimes \xi_{i})\cr 
&=(L_0-H -t(j^2 +j))(v_{-j}\otimes \xi_{-j})\cr
&=((j^2 +j)t -j +j  -(j^2 +j)t)(v_{-j}\otimes \xi_{-j})\cr
&=0.\cr}$$
(In fact a similar computation shows directly that $M(v_{i}\otimes \xi_{i})=0$.) The result for $A=L$ and $A=K$ 
holds because $L_{1} \xi=0=L_{2}$ and 
$Ev_{i}$ and $E^2 v_i$ are multiples of $v_{i+1}$ and $v_{i+2}$ respectively.    

Now we check that $PAw=0$ for the three choices of $A$. We first observe that, since $EH=(H-I)E$, we have 
$EA=A^\prime E$ with
$$M^\prime=L_0-H-tC +I,\,\, L^\prime=L_1-E(t(H-2)+1) ,\,\, K^\prime=L_2-tE^2(t(H-{5\over 2})+{7\over 4}).$$   
Since $[N,M]=aN$, we have
$NM=MN +a N$
so that
$$ENM=EMN+aEN=M^\prime EN + a EN.$$
Hence
$$PM=M^\prime P + a P.$$
So that $PMw=0$, since $Pw=0$. By uniqueness, $Mw=0$.
Similarly $NL=LN + bEN + cM$ so that
$$ENL=ELN+bE^2N +cEM=L^\prime EN + bE^2N + cEM.$$
Hence $PLw=0$, since $Pw=0$ and $Mw=0$. By uniqueness, $Lw=0$.
Finally $NK=KN +dE^2N + eEM + fL$ so that
$$ENK=K^\prime EN+dE^3N + eE^2M + fEL.$$
Hence $PKw=0$, since $Nw=0$, $Mw=0$ and $Lw=0$. By uniqueness, $Kw=0$. 
\vskip .1in

\noindent \bf APPENDIX~C: Proof of Feigin--Fuchs product formula using explicit formula for singular vectors.
\vskip .1in
\noindent \bf Lemma~A. \it Let $A$ be an upper triangular $n\times n$ matrix, 
$B$ the lower triangular matrix with $1$'s below the diagonal and 
$0$'s elsewhere and suppose $A-B$ is invertible. 
Then the solution $v=\sum v_i e_i$ of $(A-B)v=e_1$ has 
$v_n=\det (A-B)^{-1}$.
\vskip .05in
\noindent \bf Proof. \rm Let $X=A-B$ and $Y=X^{-1}$. Thus 
$v_n=y_{n1}=(-1)^{n-1} \det D/\det(X)$ where $D$ is the $(1,n)$ minor 
obtained by
deleting the first row and last column of $X$. This matrix is upper triangular
with entries $-1$ on the diagonal, so that $\det D =(-1)^{n-1}$. Hence
$v_n=(\det X)^{-1}$ as required.
\vskip .1in
\noindent \bf Lemma~B. \it There is a unique vector $w=\sum w_ie_i$ 
with $w_n=1$ such that  $(A-B)w=ae_1$.  For this solution $a=\det (A-B)$.
\vskip .05in
\noindent \bf Proof. \rm  The solution satisfies the recurrence relations
$w_n=1$ and for $1\le k\le n$
$$w_k=\sum_{i=k}^n a_{ki}w_i$$
and is therefore uniquely determined. If $A-B$ is invertible the result 
follows immediately from Lemma~A. If not then $(A-B)u=0$ 
has a non--zero solution $u=\sum u_i e_i$. If $u_n\ne 0$, then, 
rescaling if necessary, we may assume that $u_n=1$. But then
by uniqueness $a=0=\det (A-B)$, as required. If $u_n=0$, then $w^\prime=w+u$ 
would satisfy $w^\prime_n=1$ and $(A -B)w^\prime=e_1$, so that by uniqueness
$u=0$, a contradiction. The result follows.

\vskip .1in
\noindent \bf Lemma~C.  \it  Let $P_j$ be the operator giving 
the singular vector in the Verma module $M(1,j^2)$ with $j$ a 
non--negative half integer. Let $p\ge 0$ be an integer.
Then for the natural action $\ell_n= -z^{n+1} d/dz -p^2(n+1)z^n$
on ${\Bbb C}[z,z^{-1}]z^{\mu}$ with $c=0$, we have $P_jz^{\mu+p^2}=az^{-d+\mu+p^2}$ where
$a=(-1)^d\det (-F +(I+E)^{-1}(-\mu I -p^2 I +(2p+1)H+jI))$.
\vskip .05in
\noindent \bf Remark. \rm Only the case $p=1$ is needed in our applications.

\vskip .05in
\noindent \bf Proof. \rm Let $\lambda=p^2$ and $\nu=\mu+\lambda$. We set
$$u(z)= \sum_{i=0}^{2j} a_iz^{-i +\nu}v_{-j+i},\eqno{(1)}$$
with $a_0=1$.  The coefficient $a$ is determined as the unique solution of
$$av_{j}z^{\nu -2j -2}= -Fu(z) -\sum_{m\ge 0} (-E)^m z^{-m}u^\prime(z) -\sum_{m\ge 0} \lambda m (-E)^m z^{-m+1}u(z).$$
Thus 
$$av_j z^{\nu -2j-2}=-Fu(z) +\sum_{m\ge 0} (-E^m)z^{-m}(\mu z^{\nu-1}v_{-j} +\sum_{i=1}^{2j} 
(-i+\nu) z^{-i-1+\nu}v_{-j+i}) +\lambda E(I+Ez^{-1})^{-2} u(z).\eqno{(2)}$$
Let $w(z)=z^{-\nu} u(z)$. Thus multiplying by $z^{-\nu}$, we get 
$$av_j z^{-2j-2} =-Fw(z) +\sum_{m\ge 0} (-E)^m z^{-m}(-\nu z^{-1} v_{-j} +\sum_{i=1}^{2j+} 
(-i+\nu) z^{-i-1}v_{-j+i} +\lambda E(I+Ez^{-1})^{-2} w(z).\eqno{(3)}$$ 
Let $w=w(1)=v_{-j} +a_1v_{-j+1} + a_2v_{-j+2} +\cdots$ and set $z=1$ in (3). Setting $a_0=1$, this yields
$$av_j=-Fw +\sum_{m\ge 0} (-E)^m \sum_{i=0}^{2j} (-i+\mu) a_i v_{-j+i} +\lambda E(I+E)^{-2} w.$$   
Now
$$Hv_{-j+i}=(-j+i)v_{j-i},$$
so that
$$(-i-\mu)v_{-j+i} =(-H-j-\mu)v_{-j+i}.$$
Hence we get
$$\eqalign{av_j &= -Fw +\sum_{m\ge 0} (-E)^m (H+j)w +\lambda E(I+E)^{-2}w\cr
&=[-F+(I+E)^{-1}(\mu I - H -jI)]w+\lambda E(I+E)^{-2}w.\cr}$$
It follows from Lemma~B that
$$a=\det (-F +(I+E)^{-1}(\nu I - H -jI)) +\lambda E(I+E)^{-2}).\eqno{(4)}$$
When $\lambda=0=p$, this immediately gives the result for $p=0$. When $\lambda=1=p$,
we have $HE=EH +E$, so that $H(I+E)=(E+I)H +E$ and hence
$$(I+E)^{-1}H(I+E)=H+(I+E)^{-1}E.\eqno{(5)}$$
Similarly $FE=EF-2H$, so that $F(I+E)=(I+E)F-2H$ and hence
$$(I+E)^{-1}F(I+E)=F-2(I+E)^{-1}H.\eqno{(6)}$$
So in this case
$$-F +(I+E)^{-1}(-\nu I +H +jI)) +E(I+E)^{-2})=(I+E)^{-1}(-F+(I+E)^{-1}(-\nu I +3H +jI)(I+E),$$
and hence 
$$a=\det (-F+(I+E)^{-1}(-\nu I +3H +jI).$$
It follows by induction from (5) and (6) that for $p\ge 1$
$$(I+E)^{-p} H(I+E)^p=H+p(I+E)^{-1}E,\,\, (I+E)^{-p} F(I+E)^p=F-2p(I+E)^{-1}H -(p^2-p) (I+E)^{-2}E.$$
But then
$$-F +(I+E)^{-1}(H+\alpha I) +p^2 (I+E)^{-2}E= (I+E)^{-p}(-F +(I+E)^{-1}(H+\alpha I))(I+E)^p.$$
Hence
$$a=\det (-F +(I+E)^{-1}(H+jI -\nu I),$$
as required.
\vskip .1in
\noindent \bf Lemma~D. \it Let $P_j$ be the operator giving 
the singular vector in the Verma module $M(1,j^2)$ with $j$ a 
non--negative half integer. Let $p\ge 0$. Then for the natural action $\ell_n = -z^{n+1} d/dz -p^2(n+1)z^n$ 
on ${\Bbb C}[z,z^{-1}]z^{\mu +p^2}$ with $c=0$, we have $Pz^{\mu+p^2}=az^{-2j-1+\mu+p^2}$ where
$$a=(-1)^d\prod_{k\in S} (\mu -j^2 +k^2),$$
and $S=\{-j,-j+1,\cdots,j-1,j\}$ is the set of eigenvalues of $H$ on $V_j$.
\vskip .05in
\noindent \bf Proof. \rm By Lemmas B and C, we have
$$a=\det(-F +(I+E)^{-1}(-\mu-p^2 +j +(2p+1)H)).$$
Since $\det(I+E)=1$, it follows that
$$a=\det (-F - EF -\mu + j +H).$$
On the other hand
$$C={1\over 2}(EF+FE+2H^2)=EF-H+H^2$$
is a central operator acting on $v_{-j}$ as $j^2+j$. Hence
$$-EF=-j-j^2 -H +H^2,$$
so that
$$\eqalign{a&=\det (-F   -\mu-p^2 -j -H  -j^2 +H^2 + j +(2p+1)H)\cr
&=\det(j^2 -2pH -\mu-p^2 -H^2)\cr
&=\det( -\mu - (H+p)^2 +j^2)\cr
& =(-1)^d\prod_{k\in S} (\mu -j^2+(k+p)^2),\cr}$$
as required. 
\vskip .1in

\noindent \bf APPENDIX D: Holomorphic vector bundles and flat connections. \rm
Given a holomorphic vector bundle on the Riemann sphere ${\Bbb C}\cup\{\infty\}$, 
we can restrict it to the covering by two discs
$\{z:|z|<2R\}$ and $\{z:|z|>r/2\}$ where $R>1>r$ which can be refined to smaller discs
$\{z:|z|<R\}$ and $\{z:|z|>r\}$. We may identify each disc with a disc $D^\prime=\{z:|z|<1+\delta\}$ 
and the smaller disc with the unit disc $D=\{z:|z|<1\}$. Because it arises by restriction, 
there is a finite covering of $D^\prime$ by opens $U_i$ and holomorphic maps 
$g_{ij}:U_i\cap U_j \rightarrow GL_n({\Bbb C})$ with $g_{ij}g_{jk}=g_{ik}$ on $U_i\cap U_j\cap U_k$. Taking a 
connection on the corresponding complex vector bundle, 
parallel transport along lines radiating from the origin trivialises the vector bundle. Thus there are smooth 
maps $h_i:U_i\rightarrow GL_n({\Bbb C})$ such that $g_{ij}=h_i h_j^{-1}$. Thus 
$g_{ij}\cdot \partial_{\overline{z}} h_j=\partial_{\overline{z}}h_i$ and hence
$$h_i^{-1}\partial_{\overline{z}}h_i=h_j^{-1}\partial_{\overline{z}}h_j$$
on $U_i\cap U_j$. It follows that there is is a $C^\infty$ map $A:D^\prime\rightarrow M_n({\Bbb C})$ such that 
$A= h_i^{-1}\partial_{\overline{z}}h_i$ on $U_i$. We claim that there is a smooth map 
$f:D^\prime\rightarrow GL_n({\Bbb C})$ with 
$$\partial_{\overline{z}}f=-Af\eqno{(*)}$$ 
For then 
$$\partial_{\overline{z}}(h_if)=\partial_{\overline{z}}(h_i)f+h_i\partial_{\overline{z}}(f)
=h_i(Af+\partial_{\overline{z}}(f))=0.$$
Thus $k_i=h_if$ is holomorphic on $U_i$ with $g_{ij}=k_ik_j^{-1}$, as required. We need the following 
generalisation of Dolbeault's lemma.

\vskip .1in
\noindent \bf Lemma~A. \it Let $A\in C^\infty(D_1,M_s)$ and $h\in  C^\infty(D_1,{\Bbb C}^s)$.
Then we can find a solution $f\in  C^\infty(D,{\Bbb C}^s)$ of the homogeneous equation 
$\partial_{\overline{z}}f=-Af$ on $D$. The inhomogenous equation  $\partial_{\overline{z}}f=-Af +h$ on $D_r$ 
is solvable if the dual equation has the {\rm analytic continuation property}: with $R$ fixed in $(1,1+\delta)$, 
any solution 
of $\partial_{z}f=A^*f$ in $D_R=\{z:|z|<R\}$ 
vanishing on some open $U$ in $D_R$ is identically zero in $D_R$. 

\vskip .05in
\noindent \bf Remark. \rm In the example above the solutions of $\partial_{\overline{z}}f=-Af$ in $U\subset D$ 
correspond to holomorphic sections of the holomorphic vector bundle over $U$. Indeed $\xi_i(z)=h_i(z)f(z)$ satisfies
$g_{ij}\xi_j=\xi_i$, so it is a section, and $\partial_{\overline{z}}\xi_i=0$. So it is holomorphic. Conversely, 
if $(\xi_i)$ is a holomorphic section, then $h_i^{-1}\xi_i$ fit together to give a function $f(z)$ which evidently
satisfies $\partial_{\overline{z}}f=-Af$ in $U\subset D$. Similarly if we define dual bundle by 
$g_{ij}^\prime(z)=(g_{ij}(z)^{t})^{-1}$
then $h_i^\prime(z)=(h_i(z)^t)^{-1}$ and 
$$A^\prime(z)= h_i^t (\partial_{\overline{z}}h_i^{-1})^t =- A(z)^t.$$
Similarly the corresponding quantities for the conjugate antiholomorphic bundle $g^c_{ij}(z)=(g_{ij}(z)*)^{-1}$ are
$h_{i}^c(z)=(h_{i}(z)*)^{-1}$ and $A^c(z)=-A(z)^*$. Thus any solution of $\partial_{z}f=A^*f$ in $U$ yields an 
antiholomorphic section of the antiholomorphic vector bundle over $U$.
 
\vskip .05in
\noindent \bf Proof. \rm Let $\psi\in C_c^\infty(D)$ be a bump function equal to $1$ on $D_r$ and $0$ off $D_\rho$ for
some $r<\rho<1$. Replacing $A$ by $\psi A$ and $h$ by $\psi h$, we may assume that $X$ lies in $C^\infty_c(D, M_s)$
and $h$ lies in $C_c^\infty(D,{\Bbb C}^s)$ and both therefore extend to the whole of ${\Bbb C}$. 
To solve the equation $\partial_{\overline{z}}f=-Af+h$ on $D_r$, we include the disc $D$ in 
a large square  $F=[-R,R]\times [-R,R]$. By identifying 
opposite sides, doubly periodic functions on $F$ can be identified with functions on a torus $T={\Bbb T}^2$. The operator
${\cal D}=\partial_{\overline{z}}=\partial_x +i\partial_y$. Let $H_k(T)$ be the L${}^2$ Sobolev spaces for $T$ 
constructed using the Laplacian operator $\Delta={\cal D}^*{\cal D}={\cal D}{\cal D}^*=-\partial_x^2 -\partial_y^2$
(see {\bf [5]}, {\bf [18]}, {\bf [35]} or {\bf [38]}). The operator $D$ defines
a Fredholm operator of index $0$ from $H_k(T)$ to $H_{k-1}(T)$, since it is diagonalised in the natural basis. 
Its kernel consists of the constant functions and its image is the orthogonal complement of the constant functions.
The operator ${\cal D}+A$ is therefore also Fredholm of index zero from $H_k(T)\otimes M_s$ to $H_{k-1}(T)\otimes M_s$. Thus the equation
$$({\cal D}+A)f=h+g\eqno{(1)}$$ 
is soluble provided $(g+h,w_i)_{(k-1)}=0$ for finitely many vectors $w_1,\dots ,w_p\in H_{k-1}(T)\otimes M_n$. 
But if $U$ is an of 
open set $[-R,R]\times [-R,R]$ with $\overline{U} \cap \overline{D_r}=\emptyset$, then $C_c^\infty(U)$ 
is embedded in $H_{k-1}(T)$. We need to justify why there will be a vector $g$ in $C_c^\infty(U)$  satisfying 
$(g+h,w_i)_{(k-1)}=0$ ($i=1,\dots,p$). In the homogeneous case, when $h=0$, this is clear because the image of $C_c^\infty(U)$ 
contains a subspace of dimension $p+1$. 

In the inhomogeneous case, the $w_i$ will be the smooth functions in the finite--dimensional kernel of 
${\cal D}^*+A*=-\partial_{z}+A(z)^*$. By assumption their restrictions to some $U$ with 
$\overline{U}\cap \overline{D_r}=\emptyset$ are
linearly dependent, i.e. regarded as elements of $L^2(U,{\Bbb C}^s)\subset H_0(T)$. 
Since  $C^\infty_c(U,{\Bbb C}^s)$ is dense in $L^2(U,{\Bbb C}^s)$,  
we can find a solution $g\in C^\infty_c(U,{\Bbb C}^s)$ of $(g+h,w_i)=0$ ($i=1,\dots,p$). 
Hence the conditions are satisfied with $k=1$.
\vskip .1in
\noindent \bf Corollary. \it The inhomogeneous equation is always 
solvable in the example coming from a holomorphic vector bundle.
\vskip .05in
\noindent \bf Proof. \rm Linearly independent antiholomorphic 
sections of the conjugate bundle stay linearly independent on any open $U$. Indeed in this case if $w$ is in the kernel of ${\cal D}^*+A*=-\partial_{z}+A(z)^*$ and vanishes outside $\overline{D_r}$ then
it vanishes everywhere (by antiholomorphicity inside $D_R$ for $r<R<1$).
 
\vskip .1in
\noindent \bf Lemma~B. \it If $A\in C^\infty(D,M_n)$ and $D_r=\{z:|z|<r\}$ for $0<r<1$, then we can find
$B\in C_c^\infty(D,M_n)$ such that $\partial_{\overline{z}} +A$ and $\partial_z+B$ commute on $D_r$.
\vskip .05in
\noindent \bf Proof. \rm The commutativity condition on $D_r$ is equivalent to
$$\partial_{\overline{z}}B=\partial_{z}A-{\rm ad}(A(z))\cdot B,$$
which can be solved by Lemma~A with $s=n^2$, $f=B$, $h=\partial_{z}A$ and the operator $A(z)$ given by $-{\rm ad}(A(z))$.

\vskip .1in 
\noindent \bf Lemma~C. \it If $A\in C^\infty(D,M_n)$ and $D_r=\{z:|z|<r\}$ for $0<r<1$, then we can find
$f\in  C^\infty(D,M_n)$ such that $\partial_{\overline{z}}f=-Af$ on $D_r$ with $f(z)$ invertible on 
$D_r$.
\vskip .05in
\noindent \bf Proof. \rm Taking $B(z)$ as in Lemma~B, we get
operators $\partial_{\overline{z}} +A$ and $\partial_z+B$ commuting on $D_r$.
These define a flat connection so that parallel transport along any path from $0$ results in an smooth map 
$f$ of $D_r$ into invertible matrices such that $f(0)=I$ and 
$$\partial_{\overline{z}}f =-Af,\,\,\partial_zf=-Bf$$.
\vskip .1in

\noindent \bf Remarks. 1. \rm Standard approximation techniques 
by polynomials (see {\bf  [9]}) can be 
applied to deduce that any holomorphic vector bundle on the unit disc is trivial. 

\noindent \bf 2. \rm The same techniques can be used to prove that a holomorphic vector bundle on an open ball $B$ in
${\Bbb C}^n$ is trivial. In that case the we have commuting operators $\partial_{\overline{z}_i} +A_i$ and 
must find commuting operators $\partial_{z_i} +B_i$. This can be accomplished 
by successively
using operators on L${}^2$ Sobolev spaces of ${\Bbb T}^{2n}$ corresponding to the de Rham complex 
in the first $2k$ real variables in ${\Bbb T}^{2n}$ and the Dolbeault complex
in the last $2n-2k$ real variables (regarded as $n-k$ complex variables).

\vskip .2in

\noindent \bf REFERENCES. \rm
\vskip .1in
\item{1.} A. Astashkevich, {\it On the structure of Verma modules over Virasoro and Neveu-Schwarz algebras,} 
Comm. Math. Phys. {\bf 186} (1997), 531–-562.
\item{2.} A. Astashkevich and D. Fuchs, {\it Asymptotics for singular vectors in 
Verma modules over the Virasoro algebra,}  Pacific J. Math.  {\bf 177}  (1997), 201–-209.
\item{3.} M. Bauer, P. Di Francesco, C. Itzykson and J.--B. Zuber, {\it Covariant differential equations and 
singular vectors in Virasoro representations,}  Nuclear Phys. B  {\bf 362}  (1991), 515–-562.
\item{4.} L. Benoit and Y. St--Aubin, {\it A proof of the explicit formula for certain singular 
vectors of the Virasoro algebra} in  ``Lie theory, differential equations and representation theory,''  77–-84, 
Univ. Montr\'eal, 1990. 
\item{5.} L. Bers, F. John and M. Schechter, ``Partial differential equations,'' 
Lectures in Applied Mathematics, A.M.S., 1979.
\item{6.} R. E. Borcherds, {\it Vertex algebras, Kac-Moody algebras, and the Monster,} 
Proc. Nat. Acad. Sci. U.S.A. {\bf 83} (1986), 3068–-3071.
\item{7.} P. Di Francesco, P. Mathieu and D. S\'en\'echal, ``Conformal field theory,'' 
Graduate Texts in Contemporary Physics, Springer-Verlag, 1997.
\item{8.} B. Feigin and D. Fuchs, {\it Representations of the Virasoro algebra} 
in  ``Representations of Lie groups and related topics,''  465–-554, 
Adv. Stud. Contemp. Math. {\bf 7}, Gordon and Breach, 1990. 
\item{9.} O. Forster, ``Lectures on Riemann Surfaces,'' Graduate Texts in Mathematics, {\bf 81}, 
Springer-Verlag, 1981.
\item{10.} I. Frenkel and V. Kac, {\it Basic representations of affine Lie algebras and dual resonance models,}  
Invent. Math. {\bf 62} (1981), 23–-66. 
\item{11.} I. Frenkel, J. Lepowsky and A. Meurman, ``Vertex operator algebras and the Monster,'' 
Pure and Applied Mathematics, {\bf 134},  Academic Press, 1988. 
\item{12.} D. Friedan, Z. Qiu and S. Shenker, {\it Details of the nonunitarity proof for highest 
weight representations of the Virasoro algebra,}  Comm. Math. Phys.  {\bf 107}  (1986), 535–-542.
\item{13.} D. Fuchs, {\it Singular vectors over the Virasoro algebra and extended Verma modules} in  
``Unconventional Lie algebras'',  65–74, Adv. Soviet Math. {\bf 17}, A.M.S., 1993.

\item{14.} P. Goddard and D. Olive,  {\it Kac--Moody and Virasoro algebras in relation to quantum physics}, 
Internat. J. Modern Phys. A  {\bf 1}  (1986),  303-–414.
\item{15.} P. Goddard, A. Kent and D. Olive, {\it Unitary representations of the Virasoro and super--Virasoro  algebras,}
Comm. Math. Phys.  {\bf 103}  (1986), 105–-119.
\item{16.} I. Gohberg and M. G. Krein, ``Theory and applications of Volterra operators in Hilbert space,'' Translations of Mathematical Monographs {\bf 24}, A.M.S., 1970. 
\item{17.} J. F. Gomes, {\it The triviality of representations of the Virasoro algebra with vanishing 
central element and $L_0$ positive,}  Phys. Lett. B  {\bf 171}  (1986) 75–-76.
\item{18.} P. Griffiths and J. Harris, ``Principles of algebraic geometry,'' Wiley-Interscience, 1978
\item{19.} A. Grothendieck, {\it Sur la classification des fibr\'es holomorphes sur la sph\`ere de Riemann,}  (French)
Amer. J. Math. {\bf 79} (1957), 121–-138. 
\item{20.} J. C. Jantzen, ``Moduln mit einem höchsten Gewicht,'' Lecture Notes in Mathematics {\bf 750}, Springer, 1979.
\item{21.} V. Kac, {\it Contravariant form for infinite-dimensional Lie algebras and superalgebras,} 
in  ``Group Theoretical Methods in Physics',' Lecture Notes in Physics, {\bf 94}, Springer, 1979, 441--445.
\item{22.} V. Kac and A. K. Raina, ``Bombay lectures on highest weight representations of 
infinite-dimensional Lie algebras,'' Advanced Series in Mathematical Physics, {\bf 2}, World Scientific, 1987.
\item{23.} V. Kac, ``Infinite--dimensional Lie algebras'' (third edition), Cambridge University Press, 1990.
\item{24.} V. Kac, ``Vertex algebras for beginners,'' University Lecture Series, {\bf 10}, A.M.S., 1997.
\item{25.} J.--L. Koszul and B. Malgrange, Koszul, {\it Sur certaines structures fibr\'ees complexes,} 
Arch. Math. {\bf 9} 1958 102–-109. 
\item{26.} G. D. Landweber,
{\it Multiplets of representations and Kostant's Dirac operator for equal rank loop groups,}
Duke Math. J. {\bf 110} (2001), 121--160. 
\item{27.} D. E. Littlewood, ``The Theory of Group Characters and Matrix Representations of Groups,'' 
Oxford University Press, 1940.
\item{28.} T. Loke, ``Operator algebras and conformal field theory of the discrete series
representations of Diff(S${}^1$)'', Ph.D. dissertation, University of Cambridge, 1994.
\item{29.} I. G. Macdonald, ``Symmetric functions and Hall polynomials,'' Oxford University Press, 1979
\item{30.} B. Malgrange, ``Lectures on the theory of functions of several complex variables,''
 Tata Institute of Fundamental Research 
Lectures on Mathematics and Physics {\bf 13}, Springer-Verlag, Berlin, 1984.
\item{31.} A. Pressley and G. B. Segal,  ``Loop groups'', Oxford University Press, 1986.
\item{32.} G. B. Segal, {\it Unitary representations of some infinite-dimensional groups}.  
Comm. Math. Phys.  {\bf 80}  (1981), 301–-342. 
\item{33.} M. Wakimoto and H. Yamada, {\it The Fock representations of the Virasoro algebra 
and the Hirota equations of the modified KP hierarchies},  Hiroshima Math. J.  {\bf 16}  (1986), 427–-441.
\item{34.} N. R. Wallach, {\it Classical invariant theory and the Virasoro algebra} in  
``Vertex operators in mathematics and physics (Berkeley, Calif., 1983)'',  475–-482, M.S.R.I. Publ. 
{\bf 3}, Springer, 1985.
\item{35.} F. W. Warner, ``Foundations of differentiable manifolds and Lie groups'', 
Graduate Texts in Mathematics, {\bf 94}, Springer-Verlag, 1983.
\item{36.} A. J. Wassermann, ``Kac-Moody and Virasoro algebras'', Part III lecture notes, Cambridge, 1998. 
http://arXiv:1004.1287v2
\item{37.} A. J. Wassermann, {\it Operator algebras and conformal field theory. III. 
Fusion of positive energy representations of ${\rm LSU}(N)$ using bounded operators.}  Invent. Math.  
{\bf 133}  (1998), 467–-538.
\item{38.} A. J. Wassermann, ``Analysis of operators,'' Part III lecture notes, Cambridge, 2006. \break 
http://www.dpmms.cam.ac.uk/$\sim$ajw 
\item{39.} A. J. Wassermann, Primary fields, operator product expansions and Connes fusion, in preparation.
\item{40.} H. Weyl, ``The Classical Groups: Their Invariants and Representations,'' Princeton University Press, 1939. 
\vfill\eject
\end

APPENDIX B: A differential geometric proof of the triviality 
of holomorphic vector bundles on the disc. \rm
Let $D$ be the open unit disc with a finite covering by opens $U_i$ and holomorphic maps 
$g_{ij}:U_i\cap U_j \rightarrow GL_n({\Bbb C})$ with $g_{ij}g_{jk}=g_{ik}$ on $U_i\cap U_j\cap U_k$. Taking a 
connection on the corresponding complex vector bundle, 
parallel transport along lines radiating from the origin trivialises the vector bundle. Thus there are smooth 
maps $h_i:U_i\rightarrow GL_n({\Bbb C})$ such that $g_{ij}=h_i h_j^{-1}$. This 
$g_{ij}\cdot \partial_{\overline{z}} h_j=\partial_{\overline{z}}h_i$ and hence
$$h_i^{-1}\partial_{\overline{z}}h_i=h_j^{-1}\partial_{\overline{z}}h_j$$
on $U_i\cap U_j$. It follows that there is is a $C^\infty$ map $A:U\rightarrow M_n({\Bbb C})$ such that 
$A= h_i^{-1}\partial_{\overline{z}}h_i$ on $U_i$. We claim that there is a smooth map 
$f:U\rightarrow GL_n({\Bbb C})$ with 
$$\partial_{\overline{z}}f=-Af\eqno{(*)}$$ 
For then 
$$\partial_{\overline{z}}(h_if)=\partial_{\overline{z}}(h_i)f+h_i\partial_{\overline{z}}(f)
=h_i(Af+\partial_{\overline{z}}(f))=0.$$
Thus $k_i=h_if$ is holomorphic on $U_i$ with $g_{ij}=k_ik_j^{-1}$, as required. We need the following 
generalisation of Dolbeault's lemma.

\vskip .1in
\noindent \bf Lemma~A. \it Let $A\in C^\infty(D,M_s({\Bbb C}))$, $h\in  C^\infty(D,{\Bbb C}^s)$  and 
$D_r=\{z:|z|<r$ for $0<r<1$. 
Then we can find a solution $f\in  C^\infty(D,{\Bbb C}^s)$ of the homogeneous equation 
$\partial_{\overline{z}}f=-Af$ on $D_r$. The inhomogenous equation  $\partial_{\overline{z}}f=-Af +h$ on $D_r$ 
is solvable if the dual equation has the {\rm analytic continuation property}: any solution 
$\partial_{z}f=A^*f$ in $D_R$ with $r<R<1$ vanishing on some open $U$ in $D_R$ is identically zero in $D_R$. 

\vskip .05in
\noindent \bf Remark. \rm In the example above the solutions of $\partial_{\overline{z}}f=-Af$ in $U\subset D$ 
correspond to holomorphic sections of the holomorphic vector bundle over $U$. Indeed $\xi_i(z)=h_i(z)f(z)$ satisfies
$g_{ij}\xi_j=\xi_i$, so it is a section, and $\partial_{\overline{z}}\xi_i=0$. So it is holomorphic. Conversely, 
if $(\xi_i)$ is a holomorphic section, then $h_i^{-1}\xi_i$ fit together to give a function $f(z)$ which evidently
satisfies $\partial_{\overline{z}}f=-Af$ in $U\subset D$. Similarly if we define dual bundle by 
$g_{ij}^\prime(z)=(g_{ij}(z)^{t})^{-1}$
then $h_i^\prime(z)=(h_i(z)^t)^{-1}$ and 
$$A^\prime(z)= h_i^t (\partial_{\overline{z}}h_i^{-1})^t =- A(z)^t.$$
Similarly the corresponding quantities for the conjugate antiholomorphic bundle $g^c_{ij}(z)=(g_{ij}(z)*)^{-1}$ are
$h_{i}^c(z)=(h_{i}(z)*)^{-1}$ and $A^c(z)=-A(z)^*$. Thus any solution of $\partial_{z}f=A^*f$ in $U$ yields an 
antiholomorphic section of the antiholomorphic vector bundle over $U$.
 
\vksip .05in
\noindent \bf Proof. \rm Let $\psi\in C_c^\infty(D)$ be a bump function equal to $1$ on $D_r$ and $0$ off $D_s$ for
some $r<s<1$. Replacing $A$ by $\psi A$ and $h$ by $\psi h$, we may assume that $X$ lies in $C^\infty_c(D, M_n({\Bbb C}))$
and $h$ lies in $C_c^\infty(D,{\Bbb C}^s)$ and both therefore extend to the whole of ${\Bbb C}$. 
To solve the equation $\partial_{\overline{z}}f=-Af+h$ on $D_r$, we include the disc $D$ in 
a large square  $F=[-R,R]\times [-R,R]$. By identifying 
opposite sides, doubly periodic functions on $F$ can be identified with functions on a torus $T$. The operator
$D=\partial_{\overline{z}}=\partial_x +i\partial_y$. Let $H_s(T)$ be the Sobolev spaces for $T$ 
constructed using the Laplacian operator $\Delta=D^*D=DD^*=-\partial_x^2 -\partial_y^2$. The operator $D$ defines
a Fredholm operator of index $0$ from $H_s(T)$ to $H_{s-1}(T)$, since it is diagonalised in the natural basis. 
Its kernel consists of the constant functions and its image is the orthogonal complement of the constant functions.
The operator $D+A$ is therefore also Fredholm of index zero from $H_s(T,M_n)$ to $H_{s+1}(T,M_n)$. Thus the equation
$$(D+A)f=h+g\eqno{(1)}$$ 
is soluble provided $(g+h,h_i)_{(s+1)}=0$ for finitely many vectors $h_1,\dots ,h_k\in H_{s+1}(T,M_n)$. 
But if $U$ is an of 
open set $[-R,R]\times [-R,R]$ with $\overline{U} \cap \overline{D_r}=\emptyset$, then $C_c\infty(U)$ 
is embedded in $H_{s+1}(T)$. We need to justify why there will be a vector $g$ in $C_\^\infty(U)$  satisfying 
$(g+h,w_i)_{(s+1)}=0$ ($i=1,k$). In the homogeneous case, when $h=0$, this is clear because the image of $C_c\infty(U)$ 
contains a subspace of dimension $k+1$. 

In the inhomogeneous case, the $w_i$ will be the smooth functions in the finite dimensional kernel of 
$D^*+A*=-\partial_{z}+A(z)^*$. By assumption their restrictions to some $U$ with 
$\overline{U}\cap \overline{D_r}=\emptyset$ are
linearly dependent, i.e. regarded as elements of $C^\infty(U,{\Bbb C}^s)\subset L^2((U,{\Bbb C}^s)$. 
The pairing with $C^\infty_c(U,{\Bbb C}^s)$ is non--degenerate, so, taking $s=1$, 
we can find a solution $g\in C^\infty_c(U,{\Bbb C}^s)$ of $(g+h,w_i)=0$ ($i=1,k$).
\vskip .1in
\noindent \bf Corollary~1. \it The inhomogeneous equation is always 
solvable in the example coming from a holomorphic vector bundle.
\vskip .05in
\noindent \bf Proof. \rm Linearly independent antihomolomorphic 
sections of the conjugate bundle stay linearly independent on any open $U$. A little more care is needed here.
Indeed in this case if $w$ is in the kernel of $D^*+A*=-\partial_{z}+A(z)^*$ and vanishes outside $\overline{D_r}$ then
it vanishes everywhere (by antiholomorphicity inside $D_R$ for $r<R<1$).
\noindent \bf Corollary~2 (Dolbeault's Lemma). \it The inhomogeneous equation 
$\partial_{\overline{z}}f=g$ is always solvable on $D$.
\noindent \bf Proof. \rm Take discs $D_n={z:|z|<r_n$ with $r_n\uparrow 1$. By Corollary~1, 
we can find $f_n\in C^\infty(D,M_n{(\Bbb C}))$ such that
$$\partial_{\overline{z}}f_n=g$$
on $D_n$. 
Hence $h(z)=f_{n+1}(z)-f_n(z)$ is holomorphic on $D_n$. Note that we can modify $f_{n+1}$ to $f_{n+1}+ g$ 
with $g$ holomorphic. We claim that we can choose $g$ with polynomial entries
such that $\sup_{D_{n-1}}\|h(z)-g(z)\|$ is arbitrarily small,
which can be accomplished by using the truncation of the Taylor 
expansion of $h(z)$. So we can assume that
$$\sup_{D_{n-1}} \|f_{n+1}(z)-f_n(z)\|<\delta_n,$$
with $\delta_n$ arbitrarily small. 
It follows that
$f_{n+p}$ forms a Cauchy sequence on $\overline{D_{n-1}}$ and hence has a uniform limit $f$. 
But  $f_{n+p}-f_n$ is holomorphic in $D_{n-1}$ so the limit $f-f_n$ is also holomorphic in $D_{n-1}$. 
In particular $f$ is smooth in $D_{n-1}$ and $\partial_{\overline{z}}f=\partial_{\overline{z}}f_n=g$
on $D_{n-1}$ and hence everywhere.
 
\vskip .1in
\noindent \bf Lemma~B. \it If $A\in C^\infty(D,M_n({\Bbb C}))$ and $D_s=\{z:|z|<r\}$ for $0<r<1$, then we can find
$B\in C_c^\infty(D,M_n({\Bbb C}))$ such that $\partial_{\overline{z}} +A$ and $\partial_z+B$ commute on $D_r$.
\vksip .05in
\noindent \bf Proof. \rm The commutativity condition on $D_r$ is equivalent to
$$\partial_{\overline{z}}B=\partial_{z}A-{\rm ad}(A(z))\cdot B,$$
which can be solved by Lemma~A with $s=n^2$, $f=B$, $h=\partial_{z}A$ and the operator $A(z)$ given by $-{\rm ad}(A(z))$.

\vskip .1in 
\noindent \bf Lemma~C. \it If $A\in C^\infty(D,M_n({\Bbb C}))$ and $D_r=\{z:|z|<r$ for $0<r<1$, then we can find
$f\in  C^\infty(D,M_n({\Bbb C}))$ such that $\partial_{\overline{z}}f=Af$ on $D_r$ with $f(z)$ invertible on 
$D_r$.
\vskip .05in
\noindent \bf Proof. \rm Taking $B(z)$ as in Lemma~B, we get
commuting operators $\partial_{\overline{z}} +A$ and $\partial_z+B$ commute on $D_s$.
These define a flat connection so that parallel transport along any path from $0$ results in an smooth map of 
$f$ on $D_r$ into invertible matrices such that $f(0)=I$ 
$$\partial_{\overline{z}}U =-AU,\,\,\partial_zU=-UB$$ commute on $D_r$.
\vskip .1in
This result can be extended by the technique of Corollary~2 to the whole disc as follows:
\vksip .1in
\noindent \bf Lemma~D. \it If $A\in C^\infty(D,M_n({\Bbb C}))$, then we can $f\in  C^\infty(D,GL_n({\Bbb C}))$
such that $(\partial_{\overline{z}}f=Af$. 
\vskip .05in
\noindent \bf Proof. \rm  Take discs $D_n={z:|z|<r_n$ with $r_n\uparrow 1$. By Lemma~C, 
we can find $f_n\in C^\infty(D,M_n{(\Bbb C}))$ such that
$$(\partial_{\overline{z}} +A)f_n=0$$
on $D_n$ with $f_n$ invertible on $D_r$. Note that 
$$\partial_{\overline{z}}(f_n)f_n^{-1}=\partial_{\overline{z}}(f_{n+1})f_{n+1}^{-1}$$
on $D_n$, so that
$$\partial_{\overline{z}}(f_n^{-1}f_{n+1})=-f_{n}^{-1}[\partial_{\overline{z}}(f_{n})f_n^{-1}]f_{n+1} + 
f_n^{-1}[\partial_{\overline{z}}(f_{n+1})f_{n+1}^{-1}] f_{n+1}=0.$$
Hence $h(z)=f_n(z)^{-1}f_{n+1}(z)$ is holomorphic on $D_n$. Note that we can modify $f_{n+1}$ to $f_{n+1} g$ 
with $g$ holomorphic and invertible on $D_{n+1}$. We claim that we can choose $g$ with polynomial entries
such that $\sup_{D_{n-1}}\|h(z)g(z)-I\|$ is arbitrarily small. But this is equivalent to choosing 
$\sup_{D_{n-1}} \|g(z)-h(z)^{-1}\|<\delta$ arbitrarily small, 
which can be accomplished by using the truncation of the Taylor 
expansion of $h(z)^{-1}$. Moreover, since $GL_n({\Bbb C})$ is open in $M_n({\Bbb C})$, 
$g$ will be invertible on $\overline{D_{n-1}}$ for $\delta$ sufficiently small. So we can assume that
$$\sup_{D_{n-1}} \|I-f_n(z)^{-1}f_{n+1}(z)\|<\delta_n,$$
with $\delta_n$ arbitrarily small. 
It follows that
$f_{n+p}$ forms a Cauchy sequence on $\overine{D_{n-1}}$ and hence has a uniform limit $f$. 
But  $f_n^{-1}f_{n+p}$ is holomorphic in $D_{n-1}$ so the limit $f_{n}^{-1}f$ is also holomorphic in $D_{n-1}$. 
In particular $f$ is smooth in $D_{n-1}$.   
Moreover $\partial_{\overline{z}}(f_{n}^{-1} f)=0$ so that
$$\partial_{\overline{z}}(f)=\partial_{\overline{z}}(f_n(f_{n}^{-1} f))=\partial_{\overline{z}}(f_n)(f_{n}^{-1} f)
=Af_n(f_{n}^{-1} f)=Af.$$
Since this holds for each $D_{n-1}$, $f$ is smooth in $D$ and $\partial_{\overline{z}}(f)=Af$ as required.
Finally note 
that
$$\sup _{D_{n-1}} 
\|I-f_n(z)^{-1}f_{n+p}(z)\|\le
\delta_n +(1+\delta_n)\delta_{n+1} + \cdots +(1+\delta_n)\dots (1+\delta_{n+p-1})\delta_{n+p)\le
1/2,$$ 
so that passing to the limit
$$\sup _{D_{n-1}} 
\|I-f_n(z)^{-1}f(z)\|\le 1/2.$$
Hence $f$ is invertible on $D_{n-1}$ and hence everywhere.    
\vskip .1in
\noindent \bf *Lemma E. \it Let $v_1(z),\dots,v_m(z)$ be holomorphic maps of $D$ into $V={\Bbb C}^n$ which are linearly
independent at each point. If $(e_i)$ is a basis of $V$, there exists a holomorphic map $g:D\rightarrow GL_n((\Bbb C})$ 
such that $g(z)v_s(z)=e_s$ for $s=1,dots,m$. Indeed given vectors $u_1,\dots,u_s$ 
\vskip .05in
\noindent \bf Proof~1. \rm The $C^\infty$ problem can easily be solved,  
so we can find $f$ smooth and invertible such that $f(z)v_s(z)=e_s$. In fact, orthonomalising
the vectors $v_1,\dots,v_m$, we get an orthonormal set $u_1,\dots, u_m$. It follows that the orthogonal projection
$P(z)$ onto the subspace spanned by $v_1(z), \dots,v_m(z)$ depends smoothly on $z$. Fix $\zeta$ with $|\zeta|=1$ 
and set $F(t)=I-2P(\zeta t)$. Let $A(t)={1\over
2}\dot{F}(t)F(t)^{-1}$ and let $g(t)$ be the solution of $\dot{g}(t)=A(t)g(t)$ 
with $g(0)=I$ and let $h(t)$ the solution of $\dot{h{=-h(t)A(t)$ with $h(0)=I$. Thus ${d\over dt}(hg)=0$ that
$hg=I$ and $g$ is invertible. Moreover
$${d\over dt}(g^{-1}Fg)=-g^{-1}\dot{g} g^{-1} F g +g^{-1}\dot{F}g + g^{-1}F\dot{g}=g^{-1}(-AF +2AF + FA)g=0.$$
Since $\dot{F}F=-F\dot{F}$ so that $AF+FA=0$. Since solutions of ordinary differential equations 
depend smoothly on initial
data, this gives a smooth family $g(z)$ such that $F(z)=g(z)F(0)g(z)^{-1}$. 
But then $\partial_{\overline{z}}f(z) v(z)=0$ and hence $\partial_{\overline{z}}(f)f^{-1}e_s=0$. Thus,
setting $A(z)=\partial_{\overline{z}}(f)f^{-1}$, we have
$$A(z)=\pmatrix{0 & B\cr 0 & D\cr}.$$
We look for a solution of $\partial_{\overline{z}}h(z)=A(z)h(z)$ with
$$h(z)=\pmatrix{I& X(z) \cr 0 & Y(z)\cr},$$
with $Y$ invertible. Thus we require
$$\partial_{\overline{z}}Y(z)=Y(z).$$
Lemma~D implies this can always be solved. We then use Corollary~2 to solve the inhomogeneous equation
$$\partial_{\overline{z}}X(z)=B(z)Y(z).$$
Thus $h(z)$ fixes $e_1,\dots e_m$ and $\partial_{\overline{{z}}(h)h^{-1}=\partial_{\overline{z}}(f)f^{-1}=A$. Hence
$$\partial_{\overline{z}} (h^{-1}f)=-h^{-1}Af+h^{-1}Af=0,$$
so that $g=h^{-1}f$ is the required map.

\vskip .05in
\noindent \bf Proof~2. \rm The problem can be solved in an sufficiently small open neighbourhood of every point, 
so we may assume there is a countable cover of $D$ by open discs $U_i$ 
(so that all non--empty intersections are connected), 
and holomorphic maps $h_i:U_i\rightarrow GL_n({\Bbb C})$ 
such that $h_i(z)v_s(z)=e_j$ on $U_i$. Thus $g_{ij}(z)=h_i(z)h_j(z)^{-1}$ ($z\in U_i\cap U_j$) 
satisfies $g_{ij}g_{jk}=g_{ik}$ on $U_i\cap U_j\cap U_k$  and $g_{ij}(z)e_s=e_s$. Thus 
$$g_{ij}=\pmatrix {I & B_{ij} \cr 0 & G_{ij}},$$
where $G_{ij}(z)$ takes values in $GL_{n-p}({\Bbb C})$. It also satisfies $G_{ij}G_{jk}=G_{ik}$ on 
$U_i\cap U_j\cap U_k$, and hence we can find $H_i:U_i\rightarrow GL_{n-p}({\Bbb C})$ holomorphic such that
$G_{ij}(z)=H_i(z)^{-1}H_j(z)$. Thus replacing $h_i(z)$ by 
$$\pmatrix{I & 0 \cr 0 & H_i(z)\cr}\cdot h_i(z),$$
we may assume that $G_{ij}(z)\equiv I$. But then $B_{ij}(z)+B_{jk}(z)=B_{ik}(z)$ on $U_i\cap U_j\cap U_k$. Each coordinate
of $B_{ij}(z)$ defines a holomorphic complex function $f_{ij}(z)$ on 
$U_i\cap U_I$ satisfying the same condition. Thus $\gamma_{ij}(z)=\exp f_{ij}(z)$ defines a line bundle, 
so we can find holomorphic functions $k_i$ on $U_i$ such that $\gamma_{ij}=\exp k_i(z)-k_j(z)$ on $U_i\cap U_j$.
Hence on $U_i\cap U_j$
$$f_{ij}(z)-k_i(z) + k_j(z)=2\pi i m_{ij}$$
with $m_{ij}$ an integer (since $U_i\cap U_j$ is connected). Clearly $m_{ij}+m_{jk}=m_{ik}$ if 
$U_i\cap U_j\cap U_k\ne \emptyset$. Now by elementary topology we can find integers $n_i$ such that $m_{ij}=n_i-n_j$. 
(Indeed placing vertices at the centre of each disc and joining the vertices by an edge if the disc intersect, 
we get a triangulation of $D$, each triangle corresponding to a non--empty triple intersection. If we choose one
vertex $a$ and assign the value $n_a=0$ there, 
then the value $n_b$ at any other vertex $b$ is obtained by summing the $m_{ij}$ along 
edges of a path between the vertices; the condition on triangles implies that this value is uniquely determined.)
So replacing $k_i(z)$ by $K_i(z)=k_i(z) +2\pi i n_i$, we have $f_{ij}=K_i-K_j$. Since we can do this for each coordinate, 
it follows that there are holomorphic $p\times (n-p)$ matrix--valued function $C_i(z)$ on $U_i$ such that
$B_{ij}=C_j-C_i$ on $U_i\cap U_j$. Replacing $h_i(z)$ by
$$\pmatrix{I & C_i(z) \cr 0 & I\cr}\cdot h_i(z),$$
we may get $g_{ij}(z)\equiv I$. Thus $h_i(z)=h_j(z)$ on $U_i\cap U_j$ and so they define a homolorphic map 
$h$ of $D$ into $GL_n({\Bbb C})$ such that $h(z)v_s(z)=e_s$. 
\vskip .1in
\noindent \bf Remark. \rm The same method can be used to prove that a holomorphic vector bundle on an open ball $B$ in
${\Bbb C]^n$ is trivial. In this case there are $n$ commuting operators $\overline{\partial}_i+A_i(z)$ on $B$, 
where $\overline{\partial}_i=\partial_{\overline{z_i}}=\partial_{x_i} + i \partial_{y_i}$. This time $B$ can be 
embedded in a torus $T={\Bbb T}^{2n}$. Let $\Lambda$ be the exterior algebra in the variables $d\overline{z_i}$ and let
$$\overline{\partial}
:H^s(T)\otimes \Lambda^{{\rm even}} \rightarrow H^{s-1}(T)\otimes \Lambda^{{\rm odd}}$$ 
be the operator
$$\overline{\partial}\omega = \sum \overline{\partial}_i d{\overline z}_i\omega.$$
This $\overline{\partial}$--operator is Fredholm of index $0$ and hence its extension to vector valued functions 
$$\overline{\partial}
:H^1(T,{\Bbb C}^s)\otimes \Lambda^{{\rm even}} \rightarrow H^{0}(T,{\Bbb C}^s)\otimes \Lambda^{{\rm odd}}$$
has index zero. As before operators of type $\overline{\partial} +\sum A_i(z) d\overline{z}_i$ 
are also Fredholm of index $0$. This means that given $g_1,\dots,g_n\in H^s(T,{\Bbb C}^s)$, the equation
$(\overline{\partial}_i+A_i)f=g_i$ is solvable provided the compatibility conditions 
$(\overline{\partial}_i +A_i)g_j=(\overline{\partial}_j+A_j)g_i$ are satisfied and
and $(g_j,h_{k,j})=0$ for a finite set of $h_{k,j}$'s. Again in the example the $h_{k,j}$'s will be antiholomorphic in a
ball concentric with $B$. Thus if we look for commuting operators $\partial_i + B_i$ which also commute with the operators
$\overline{\partial}_i +A_i$, we have to solve
$$\overline{\partial}_i B_j = \partial_j A_i -{\rm ad}( A_i) \cdot B_j.$$
This can be done step by step, first solving for $B_1$. To solve for  $B_2$, the extra condition
$$\partial_1B_2  +{\rm ad}(B_1)B_2=\partial_2 B_1$$
has to be added and the differential $dz_1$ added to the elliptic complex. Once $B_2$ has been determined, 
we can solve for $B_3$ by adding the equations
$$\partial_iB_j  +{\rm ad}(B_i)B_j=\partial_j B_i$$
with $j=3$ and $i=1,2$. Again $dz_1$ and $dz_2$ must be added to the elliptic complex.
Proceeding in this way, we can successively solve for $B_1,\dots,B_n$. 
The determination of $B_{k+1}$ involves the product of the de Rham complex on the first $2k$ factors of ${\Bbb T}^{2n}$ 
and the Dolbeault complex on the last $(2n-2k)$ factors. Having produced a flat connection, the proof proceeds as in the
case $n=1$, with the holomorphic trivialisation again achieved by parallel transport.
\vskip .1in

\vskip .1in
Now consider a representation $M(1,j^2)$. The point $(c,h)=(1,j^2)$ lies on the curve $\varphi_{p,1}(c,h)=0$ where
$p=2j+1$. The curve is parametrized by a variable $t$ via $c(t)=13 - 6t - 6t^{-1}$ and $h(t)=(j^2+j)t -j$, 
the point $(1,j^2)$ corresponding to the value $t=1$. If $t\ne 0,\infty$, the Verma module $M_t=M(c(t),h(t))$ is defined. 
It is a direct sum of finite energy spaces $M_t(n)$ on which $L_0$ acts as multiplication by $n+h(t)$. The
space $M_t$ has a canonical invariant bilinear form $B_t(v,w)$ defined on it. Invariance implies that the subspaces
$M_t(n)$ are orthogonal. If we take as basis of $M_t(n)$ elements $v_i=L_{-i_r}\cdots L_{-i_1}v_t$ 
with $i_1\le i_2\le \cdots$ and $\sum s\cdot i_s=n$, then all these spaces can be identified with the same space 
$V={\Bbb R}^{p(n)}$. Thus $A^{(n)}_{ij}(t)=B_t(v_i,v_j)$ is a symmetric matrix, the entries of which are polynomials in 
$t$ and $t^{-1}$ with real coefficients. Now for each $t$, Fuchs' uniqueness theorem shows that there is a unique
polynomial
$$Q_p=L_{-1}^p + a_{p-1}(t) L_{-1}^{p-1} + \cdots a_0$$
where $a_i$ is a ${\Bbb R}[t,t^{-1}]$--linear 
combination of monomials  $L_{-i_r}\cdots L_{-i_1}$ with $2\le i_1\le i_2\le \cdots$ and $\sum si_s=p-i$
such that $w_t=Q_p(t)v_t\in M_{j^2}(p)$ is a singular vector, i.e. $L_k w=0$ for $k\ge 0$. If $N$ is the submodule of $M$ 
generated by $w$, the $N_t(n)=M_t(n)\cap N_t$ has basis $L_{-i_r}\cdots L_{-i_1}w$ with $1\le i_1\le i_2\le \cdots$ and
$\sum r\cdot i_r=n-p$. 
This submodule lies in the annihilator of $B_t$. Thus $B_t$ passes to a bilinear form $\overline{B}_t$ 
on $M(n)/N_t(n)$.  Generically $\overline{B}_t$ is non--degenerate. 
Now recall that if $X$ is a matrix and $P$ a rank one projection then
$$\det (X+aP)= a\det(I-P)X(I-P) + \det X.$$
We claim that the order of $0$ as a root of $\det A(t)$ equals $k=\sum_{i\ge i} {\rm dim}\, V_i$ where if $A=\sum A_it^i$,
$V_i=\cap{j<i} {\rm ker}\, A_j$ ($i\ge 1$). We set $U_0=V_0=V$, $U_1=V_{n_1}$ the first $k>0$ with $V_k\ne V_0$, 
$U_2=V_{n_2}$ the first $k>n_1$ with $V_k\ne V_{n_1}$, and so on. We set $n_0=0$. 
Let $v_1,\dots,v_{m_1}$ be an orthonormal basis of $U_0\ominus U_1$, $v_{m_1+1}, \dots,v_{m_2}$ 
an orthonormal basis of $U_1\ominus U_2$ and so on. Thus $(A(t)v_i,v_j)=t^{n_s}(A_{n_s}v_i,v_j) +$ higher powers of $t$ for
$i>m_s$. Moreover the first term vanishes if $j>m_{s+1}$. On the other hand the matrix 
$(A_{n_s}v_i,v_j)_{m_s<i,j\le m_{s+1}$ is invertible. Thus we can divide rows $m_s+1,\dots m_{s+1}$ by $t^{n_s}$ 
and then set $t=0$. The resulting matrix is block triangular with invertible blocks on the diagonal, so is 
itself invertible. Hence the order of $0$ as a root of $\det A(t)$ equals 
$$k=\sum n_s\cdot ({\rm dim}\, U_{s-1} -{\rm dim}\, U_s)= \sum_{i\ge 1} {\rm dim}\, V_i,$$
as claimed. Since the Shapovalov form is invariant under the Virasoro algebra, it follows that the 
Jantzen filtration defines a filtration of $V_t=V(c(t),h(t))$ by submodules $V_t^{(i)}$. In particular we can calculate
$\sum_{i\ge 1} {\rm dim}\, V_t^{(i)}(n)$. It is exactly the degree of $h(t)$ as a root of the Kac 
determinant at level $n$ which is known explicitly.  Using only this information, Astashkevich 
gave an elementary proof of the character $L(c,h)$ for the discrete series $c=1-6/m(m+1)$ and 
$h=[(p(m+1)-qm)^2-1]/4m(m+1)$ with $1\le q\le p\le m-1$ 
for $m\ge 2$ and the Jantzen filtrations on $M(c,h)$ corresponding to 
arbitrary smooth curves through $(c,h)$. (This is the case III$_-$ of Feigin--Fuchs.) 
\noindent \bf Remark. \rm The level one representations of $\widehat{\s}$ can be constructed explicitly 
using vertex operators (see section~8). Not only does this give another method of obtaining the character 
formula for $\H_j$, but it also yields explicit formulas for 
$E(0)^k$ and $F(0)^k$ which can be used to deduce Jeffrey Goldstone's 
formulas for the singular vectors in the oscillator modules, as shown by Graeme Segal.        

\vskip .1in

\noindent \bf Proof~2. \rm Suppose that $q_d=0$. We prove that $P_d=0$. 
Each $q_k$ can be written as a linear combination of monomials $L_{-s}^{n_s} \cdots L_{-2}^{n_2}$ with $n_i\ge 0$. 
Set $n_1=d-k$. We take a lexicographic ordering on monomials according to $(n_1,n_2,\dots)$. Thus the first term is
$L_{-1}^d$ corresponding to $(0,0,0,\dots)$. 

Let $w=P_dv$ with $v$ the cyclic lowest energy 
vector in Verma module. 
Amongst all monomials with non--zero coefficients, there is a smallest one in the lexicographic order
$$\cdots L_{-p}^{n_p} L_{-1}^ev,$$
where $e=d-n_1\ge 0$. Since $w$ is a singular vector, we must have $L_{p-1}w=0$. On the other hand
$$L_{p-1}w=\sum_{k=0}^e L_{p-1}q_k L_{-1}^{k}v=\sum_{k=0}^e[L_{p-1},q_k]L_{-1}^kv+q_kL_{p-1}L_{-1}^kv.$$
We look for terms ending with $L_{-p}^{n-1}L_{-1}^{e+1}$ in this expression. We first note that it 
follows by induction on $k$ that if $j\ge 1$ then
$L_jL_{-1}^kv$ lies in ${\rm lin}\, \{AL_{-1}^iv:i<k,\,\, A\in {\cal U}_2\}$. Indeed
$$L_{j}L_{-1} L_{-1}^{k-1} v=(j+1)L_{j-1}L_{-1}^{k-1}v +L_{-1}L_jL_{-1}^{k-1}v,$$
which has the same form since $L_{-1}{\cal U}_2={\cal U}_2 L_1$ and $L_{-1}^{k-1}v$ is an eigenvector of $L_0$.

A monomial with $k\le e$ can be written
$u=B\cdot L_{-p}^{n} A \cdot L_{-1}^kv$ with $B$ a monomial in $L_{-j}$'s for $j>p$ and $A$ a monomial in the $L_{-j}$'s with
$2\le < p$. But then
$$L_{p-1}u=[L_{p-1},B] \cdot L_{-p}^nA L_{-1}^kv 
+ B\cdot (L_{p-1},L_{-p}^n]\cdot AL_{-1}^kv)+ BL_{-p}^n L_{p-1}AL_{-1}^kv.\eqno{(1)}$$
Since $[L_{p-1},B]$ lies in ${\cal U}_2$, the first term in this expression has no terms ending in $L_{-1}^{e+1}$.
For the second note that
$$[L_{p-1},L_{-p}^{n}]=(2p-1) \sum_{a+b=n-1} L_{-p}^a L_{-1} L_{-p}^b=(2p-1)n L_{-p}^{n-1}L_{-1} + D,$$
where $D\in {\cal U}_2$. Clearly, if $k<e$, neither $BDL_{-1}^k v$ nor $B L_{-1} AL_{-1}^k v$ contain terms 
ending in $L_{-1}^{e+1}v$. Finally for the last term, note that if $X$ is a monomial in $L_{-j}$ for $1\le j \le t$,
then for $1\le s \le t$ we have $L_s Xv=\sum a_i X_iv$, where $X_i$ has the same form as $X$ with the exponent of
$L_{-1}$ increased by at most one. Thus no terms ending in $L_{-p}^{n_p-1} L_{-1}^{e+1}$ arise this way if $k<e$. 

When $k=e$, terms have the form $u=B\cdot L_{-p}^nL_{-1}^ev$ with either $n\ge n_p>0$ or zero. If $n=0$, then
$$L_{p-1}u=L_{p-1}BL_{-1}^ev=[L_{p-1},B]L_{-1}^ev + B L_{p-1} L_{-1}^e v.$$
The first term contains no terms ending in $L_{-1}^{e+1}v$. Nor does the second, since if $kj0$, $L_{j}L_{-1}^ev$ lies in
${\rm lin}\, \{L_{-1}^iv:0\le i<e\}$. If $n>n_p$, the none of the terms in (1) give rise to a monomial ending in
$L_{-1}^{n_p-1} L_{-1}^{e+1}v $. 

Finally there could be several terms $BL_{-p}^{n}L_{-1}^e$ in $q_e$ with $n=n_p$. However taking the term ending with 
$L_{-p}^{n-1}L_{-1}^{e+1}v$ in (1), we get $B  L_{-p}^{n-1}L_{-1}^{e+1}v$. Thus 
all these terms are distinct. But then there can be no  cancellation and 
the coefficient of $ L_{-r}^{n_r} \cdots L_{-p}^{n_p-1}L_{-1}^{e+1} v$ must therefore be non--zero, a contradiction. 

For $p\ge 1$ 
we similarly see that  
$((|f|!)^{-1}L_1^{|f|} X_f\xi_{p},\xi_{p})$ equals the dimension of the 
representation $\pi_{f^\prime}$ of $U(p)$, where $f^\prime$ is the transposed signature (or Young diagram). It vanishes if $f^\prime_{p+1}$ is greater than zero and otherwise is strictly positive, 

Now for each $t$, Fuchs' uniqueness theorem shows that there is a unique
polynomial
$$Q_p=L_{-1}^p + a_{p-1}(t) L_{-1}^{p-1} + \cdots a_0$$
where $a_i$ is a ${\Bbb R}[t,t^{-1}]$--linear 
combination of monomials  $L_{-i_r}\cdots L_{-i_1}$ with $2\le i_1\le i_2\le \cdots$ and $\sum si_s=p-i$
such that $w_t=Q_p(t)v_t\in M_{j^2}(p)$ is a singular vector, i.e. $L_k w=0$ for $k\ge 0$. If $N$ is the submodule of $M$ 
generated by $w$, the $N_t(n)=M_t(n)\cap N_t$ has basis $L_{-i_r}\cdots L_{-i_1}w$ with $1\le i_1\le i_2\le \cdots$ and
$\sum r\cdot i_r=n-p$. 
This submodule lies in the annihilator of $B_t$. Thus $B_t$ passes to a bilinear form $\overline{B}_t$ 
on $M(n)/N_t(n)$.  Generically $\overline{B}_t$ is non--degenerate. 
Now recall that if $X$ is a matrix and $P$ a rank one projection then
$$\det (X+aP)= a\det(I-P)X(I-P) + \det X.$$